\documentclass[12pt]{article}
\usepackage[hmargin={1.75truecm,2truecm},vmargin={2truecm,2truecm}]{geometry}
\usepackage[english]{babel}
\usepackage[utf8]{inputenc}
\usepackage{geometry}
\usepackage[normalem]{ulem}
\usepackage{multicol,float}
\usepackage{color,cancel}
\usepackage{multirow}
\usepackage{array}
\usepackage{comment}
\usepackage{amsmath}
\usepackage{amsthm}
\usepackage{amsfonts}
\usepackage{amssymb}
\usepackage{graphicx}
\usepackage{xcolor}
\usepackage[noend]{algpseudocode}
\usepackage{algorithm}
\usepackage{makecell}
\usepackage{ulem}
\usepackage{cancel}
\usepackage{arydshln}
\usepackage{tikz}
\usepackage{pgfplots}
\pgfplotsset{compat=1.18}
\usepackage{subfigure}
\usepackage{adjustbox}
\usepackage[affil-it]{authblk}
\usepackage{natbib}

\newtheorem{example}{Example}

\newtheorem{rem}{\bf Remark}[section]
\newtheorem{thm}{Theorem}[section]

\newtheorem{prop}[thm]{Proposition}
\newtheorem{definition}[thm]{Definition}

\newcommand{\rojooscuro}[1]{\textcolor[rgb]{0.75,0.00,0.00}{#1}}

\newcommand{\SDOD}{MCFOD}
\title{The Multi-Commodity Flow Problem with Outsourcing Decisions}

\author[1]{Elena Fern\'{a}ndez}
\author[2]{Ivana Ljubi\'{c}}
\author[1]{Nicolás Zerega}
\affil[1]{\small{Department of Statistics and Operations Research, Universidad de Cádiz, Puerto Real, Spain}}
\affil[2]{\small{Department of Information Systems, Decision Sciences and Statistics, ESSEC Business School, France}}
\date{}

\begin{document}
\maketitle
\begin{abstract}
    We address a new prize-collecting problem of routing commodities in a given network with hub and non-hub nodes, in which the service of the non-hub nodes will be outsourced to third-party carriers.   The problem is modeled as a Stackelberg game: there is a major firm (leader) that decides to serve a subset of commodities. The leader aims to outsource first and third legs of transportation services to smaller carriers (who act as followers) by allocating at most one carrier to each non-hub node. The carriers try to maximize their own profits, which are influenced by the leader's offers. The goal of the leader is to determine the optimal outsourcing fees, along with the allocation of carriers to the non-hub nodes, so that the profit from the routed commodities is maximized. The optimal response of the followers must be taken into account, as the followers might refuse to serve some legs in case they are negative or do not maximize their profit.  We also study two alternative settings: one in which the outsourcing fees are fixed, and the other one in which the carriers accept any offer, as long as the resulting profit is non-negative. We prove that the set of possible outsourcing fees can be discretized and formulate the problem as a single-level mixed-integer nonlinear program. For all considered problem variants, we prove NP-hardness and propose and computationally investigate several MIP formulations. We study the computational scalability of these MIP formulations and analyze solutions obtained by varying the reservation prices of the carriers. Finally, by comparing the introduced problem variants, we derive some interesting managerial insights.
\end{abstract}

\section{Introduction}
    In recent years, there has been a significant increase in the outsourcing of various practices. This trend has been particularly prominent in the logistics and freight transportation sectors, where it encompasses activities like last-mile delivery and full integration with external operators, including 3-PL logistics partners \citep{Dimitriou2021}. Similarly, the airline industry has also experienced outsourcing in processes such as check-in, luggage management, and even cabin crew \citep{Commission2019}, among others. Notably, there are cases in which flights are outsourced to third-party airlines, by taking advantage of the existing transfer system in major airports \citep{Commission2021}. The outsourcing of these processes offers several advantages, including enhanced flexibility and a reduced dependency on hiring and training specialized staff.

    Given the increasing prevalence of this trend, it is crucial to study and model this type of situations to gain a better understanding of how the outsourcing decisions contribute to the overall net profit of a major firm when it faces such a process. It is equally important to consider the viewpoint of carriers (who will perform the outsourced tasks) as they too aim to optimize their own net profits.

    The objective of this research is to introduce a novel prize-collecting multi-commodity flow problem for the transportation of demand between different origins and destinations (commodities) when a major firm, referred to as the leader, already possesses a backbone network induced by a set of major hubs and chooses to outsource the service for the demands coming from or going to the remaining non-hub locations using third-party companies (carriers) in order to maximize its overall profit. The incoming trips, from a non-hub to a hub node are referred to as first legs, and the outgoing trips, from a hub to a non-hub node, as third legs.

    The whole process can be seen as a Stackelberg game with a single-leader and  multiple independent followers. The leader  has to make tactical decisions regarding the outsourcing of the services, namely:
    \begin{itemize}
        \item how to allocate carriers to non-hub nodes, and
        \item what should be the outsourcing fees for routing commodities through hub nodes.
    \end{itemize}
    Once they receive the offer and allocation provided by the leader, carriers decide which commodities they accept to serve, taking into account their reservation prices. We assume that the leader has perfect information regarding the carriers' reservation prices, and hence, by anticipating the optimal response of the followers, the leader also decides which commodities  will ultimately be served and how the routing of served commodities will be done through the hub network. We call this problem the \emph{Multi-Commodity Flow Problem with Outsourcing Decisions} (\SDOD). 

    \paragraph{Our Contribution:}  We initially model the \SDOD\  as a bilevel Mixed Integer Non-Linear Program (MINLP).
    We  consider an additional \emph{fixed} setting, in which the outsourcing fees are fixed (e.g., by some external authorities) and cannot be modified. We also study two \emph{relaxed} problems, in which the followers are more flexible, and are willing to accept any offer that meets their reservation price (so any offer resulting in a non-negative profit will be accepted). That way, we obtain three additional problem variants, denoted by \SDOD$_F$ , $r$\SDOD\ and $r$\SDOD$_F$, respectively. We prove that all introduced variants are NP-hard.

    Leveraging the inherent properties associated with the independency of allocations and costs of each carrier, we express the model as a single-level MINLP and discretize the outsourcing fee decisions for the \SDOD\ problem. By deriving necessary adaptations for the other three problem variants, we prove that \SDOD\ and $r$\SDOD\ result in the same optimal solution. By exploiting structural properties of the optimal solutions, we introduce several MILP formulations for the \SDOD\ and show how to adapt them for the other problem variants. Some of these MILP approaches \emph{explicitly} identify the routing of the flows for the served commodities, possibly aggregating them, whereas an alternative MILP approach uses \emph{implicit path} variables to determine the routing of commodities based on carrier allocations.

    Our computational study is based on established benchmark sets from the literature. We first study the scalability of the proposed MILP formulations. The obtained results demonstrate the superiority of the implicit paths formulation which is able to solve instances with up to 200 nodes and 6 carriers to optimality within one hour. We also provide some interesting managerial insights by analyzing how the reservation prices, the two different strategies of the follower, and fixing of outsourcing fees affect the profit of the leader and the overall service level. These findings provide motivation for studying more complex systems in the future.

    The article is organized as follows. {We close the Introduction with a review of some related literature. The \SDOD\, and the \SDOD$_F$ and their ``relaxed'' counterparts, $r$\SDOD\ and $r$\SDOD$_F$, are formally defined in Section \ref{sec:Definition}, where the notation is introduced and an illustrative example is presented. Bilevel formulations and single-level MINLP reformulations are developed in Section \ref{sec:bilevelModel}, whereas solution properties and the problem complexity are studied in Section \ref{sec:properties}. Based on these properties, in Section \ref{sec:single-levelformu} we introduce single-level MILP reformulations for the studied models, and provide an extensive computational analysis in Section \ref{sec:compu}. The paper closes in Section \ref{sec:conclu} with some conclusions and avenues for future research.}
    
    \paragraph{Literature Overview.}
    Hierarchical decision making problems in which a major company decides to outsource some activities to third parties and needs to anticipate the optimal responses of the followers, are typically modeled as Stackelberg games using the tools of bilevel optimization. For recent surveys on mixed integer optimization approaches in bilevel optimization, see \citet{Kleinert-et-al:2021,Beck-et-al:2023}; for a more general introduction to bilevel optimization see the book by \citet{Dempe:2020}.   
     
    In the context of last-mile delivery, for example, a logistics peer-to-peer platform acts as the leader and decides to outsource some of the delivery activities to individual carriers (e.g., occasional drivers, who act as followers) who receive a compensation for each parcel delivered. In \citet{mofidi2019}, the platform makes personalized offers to the carriers, who then select a subset of parcels from this list. The problem is modeled as a Stackelberg game, using a deterministic bilevel model. To deal with uncertainties regarding drivers' preferences,  \citet{horner2021} and \citet{ausseil2022} propose stochastic bilevel models. 

    While in all these studies the outsourcing fees are assumed to be given, \citet{hong2019} propose another Stackelberg game model in which the outsourcing fees are optimized along with the offers made to the carriers.  In \citet{gdowska2018} and \citet{barbosa2023} both a professional delivery fleet and a set of occasional drivers are taken into account, and it is assumed that each delivery request has a fixed probability of being rejected by the driver.

    In all of these models, routing costs for the drivers are not taken explicitly into account. This is done for the first time in a recent work by \citet{Cerulli-et-al:2024} in which the platform decides on the fees to be paid to individual drivers, while anticipating that they may refuse to deliver some parcels if the profit of the underlying route is not maximized.
    
    Contrary to last-mile deliveries, in which the whole service (i.e., delivery of a parcel) is outsourced to a single carrier, in this study we assume that the first leg of the service can be performed by one carrier, the middle leg is performed by the major company, and the third leg can be performed by another carrier. Such settings are common in freight transportation, for example.     
    In a recent study, \citet{DelleDonne:2023} propose to use a public transportation service (PTS) to route commodities in the middle leg, whereas the first leg is performed by 3-PL who transports commodities from a  distribution center to the entry points of the PTS network, and in the third leg individual carriers deliver the parcels to the final destination. This concept is known as Freight-on-Transit (FOT), and \citet{DelleDonne:2023} study the problem from a perspective of a single decision maker who makes strategic and tactical decisions regarding the design of the FOT network. However, extensions like the ones proposed in this paper would allow to model the operational decisions concerning the service of the first and last leg of FOT in a more precise manner.
    
    Our problem is closely related to the  Multi-Commodity Flow Problem (MCFP), in which we are given a set of commodities with their origins, destinations and demands, and arc-routing costs and the goal is to find routing paths for each commodity, while minimizing the total routing costs. In its basic version, the problem is polynomially solvable, see \citet{ahuja1995network}. Moreover, when there are no arc capacities, the optimal route of each commodity corresponds to the shortest path in the given network. Many NP-hard variants (that include additional network-design, scheduling or capacity decisions) have been intensively studied in the literature (see, e.g., \citet{CrainicLaporte:97,Crainic:2000}, as well as the more recent results in \citet{Legault-et-al:23,GENDRON2019203,ChoumanCG17,GendronG17} and further references therein). Our problem introduces outsourcing decisions that need to be applied to the first and last leg of the routing paths of commodities. As we will see later, these additional decisions also render the problem NP-hard, and they distinguish the problem from the known MCFP variants studied so far.  
    
    We can also find relationship with prize-collecting hub location problems, where the goal is to locate the set of hubs and choose a subset of commodities to route so as to maximize the profit, see, e.g., \citet{Alibeyg2016,Alibeyg2018,Taherkhani2019,Taherkhani2020}. Contrary to these hub location problems, in our setting, the backbone network induced by the set of open hubs is already given. 
\section{Problem definition}\label{sec:Definition}
    \paragraph{Input Data:}
    Let $G=(V,\, A)$ be a given directed graph, with $n=|V|$ nodes. Let $H\subset V$ be the set of \emph{hub nodes} and $V\setminus H$ the remaining nodes, referred to as \emph{non-hub nodes}. The backbone network, consisting of the arcs between hub nodes, is defined as $A_H= \{(i,\,j):\,i,\,j\in H\}$ forming a complete network. By definition, $A_H$ considers the existence of \emph{loop} type connections, i.e. connections with same end nodes. The arc set is defined by the union of the backbone network and the arcs between non-hub and hub nodes as:
    $$A=A_H\cup \{(i,\,j):\,i\in V\setminus H,\,j\in H\}\cup \{(i,\,j):\,i\in H,\,j\in V\setminus H \}.$$
    There is a service demand expressed by means of a set of commodities, indexed in a set $R$.
    \begin{itemize}
         \item Each commodity $r\in R$ is associated with a triplet $(o(r),\, d(r),\, w^r)$, where  $o(r),\, d(r)\in V$  denote its origin and destination nodes, respectively, and $w^r >0$ the flow that must be sent from $o(r)$ to $d(r)$. Commodities  need not be served necessarily; if commodity $r\in R$ is served, it produces a revenue $b^r > 0$.
    \end{itemize}
    We are also given a set of external carriers (i.e., followers) $K=\{1, \dots, |K|\}$. In addition:
    \begin{itemize}
         \item We are given $\overline c_{ij}^{rk} >0 $ representing the \emph{reservation price} of carrier $k\in K$ for serving commodity $r\in R$ through \emph{access} arcs $\{ (i,j) : i=o(r) \not \in H, j \in H \}$ or \emph{distribution} arcs $\{ (i,j) : i \in H, j=d(r) \not \in H \}$. This value accounts for the cost of carrier $k$ for routing the demand $w^r$   through arc $(i,\, j)$ plus an additional profit margin that carrier $k$ charges the leader to perform the service. This reservation price $\overline c_{ij}^{rk}$ may depend not only on the carrier $k$, the arc $(i,\, j)$, and the amount of demand $w^r$, but also on the commodity $r$ itself.
        \item  The cost of the leader for routing the demand $w^r$ through inter-hub arc $(i,\, j)\in A_H$ is given as $c_{ij}^{r} >0$. This cost may depend not only on the arc $(i,\, j)$ and on the amount of demand of commodity $r$, $w^r$, but also on the commodity itself. For each $r \in R$, the vector $c^r$ satisfies the triangle inequality. 
    \end{itemize}
    Table \ref{tab:NotSummary} in Appendix \ref{Apdx:Notation} summarizes the notation used in this paper. 
    
    \paragraph{Decisions of the leader:}
    The leader must decide on a set of commodities to serve together with an origin/destination route for each served commodity.    Service routes must be of the form $o(r)-i(r)-j(r)-d(r)$, where $i(r),\, j(r)\in H $ are  hub nodes (possibly $i(r)=j(r)$)  decided by the leader.  The \emph{intermediate leg} of each served commodity $r\in R$,  $(i(r),\, j(r))$, will be handled by the leader, incurring his own routing costs $c_{ij}^{r}$, whereas service of the \emph{first} and \emph{third legs},  $(o(r),\, i(r))$ and  $(j(r),\, d(r))$, respectively, will be outsourced to some external carrier(s). We assume a \emph{multiple allocation} strategy so commodities with the same non-hub origin but different destinations can be connected to the backbone network through different hubs, and flows may arrive to non-hub destinations from different hub nodes.

    The leader allocates every non-hub node $i\in V\setminus H$ to at most one carrier; this means that the first and third legs of all the routed commodities with origin or destination at node $i$ will be handled by that carrier. Furthermore, for each commodity $r\in R$ such that $o(r)\notin H$ the leader offers an \emph{outsourcing fee} $p^r_i$ for routing the first leg of $r$ through access arc $(o(r),\, i)$, where $i\in H$. Similarly, for each $r\in R$ such that $d(r)\notin H$ the leader offers an \emph{outsourcing fee} $q^r_i$ for routing the third leg of $r$ through distribution arc $(i,\, d(r))$, where $i\in H$. In each case, the involved carrier will accept or deny the service based on his own reservation price (see further description below).

    The leader eventually decides on the commodities that will be routed. When end-nodes of the commodities are non-hubs, the commodity cannot be routed unless some outsourcing fee for its first, respectively third leg has been accepted by the involved carriers. This means that the commodities entailing both a first and a third leg cannot be routed unless both a first and a third leg fee have been accepted. Therefore, the allocation of carriers to non-hub nodes becomes a major decision for the leader since,  due to their different reservation prices, different carriers could give different responses to given outsourcing fees $p^r_i$, $q^r_i$, $i\in H$. These responses condition the commodities that the leader may decide to serve, and thus the overall revenue.

    Summarizing, the decisions of the leader are the following:
    \begin{itemize}
        \item Allocate each non-hub node $i\in  V\setminus H$ to at most one carrier. Let $a: V\setminus H \mapsto K \cup \{0\}$ be a mapping that maps each non-hub node $i \in V\setminus H$  to a carrier $k \in K$, and let $a(i) = 0$ if $i$ is not allocated to any carrier.
        \item For each $r \in R$ and each  $i\in H$, determine outsourcing fees for access arcs $p^r_i$ (if $o(r)\notin H$) and distribution arcs  $q^r_i$ (if $d(r)\notin H$).
        \item Identify the set of commodities to be eventually served, $R^*\subseteq R$, among the ones with accepted outsourcing fees for their first and third legs and determine the service route for each of them.
     \end{itemize}
    Without loss of generality, we assume that all commodities $r\in R$ such that $o(r), d(r) \in H$ are preprocessed. First, the service for these commodities will not be outsourced. Second,  the profit of the leader for each of these commodities is $b_r - c_{o(r)d(r)}^r$, so the commodity will be served  if and only if $b_r - c_{o(r)d(r)}^r > 0$. 

    For each $r \in R$, $A^r\subset A_H$ denotes the set of potential interhub arcs for the service route of commodity $r$: 
        $$A^r=\begin{cases}
        \{(o(r), j): j\in H\} & \text{if } o(r)\in H, d(r)\notin H\\
        \{(i, d(r)): i\in H\} & \text{if } o(r)\notin H, d(r)\in H\\
        \{(i, j): i,j\in H\} & \text{if } o(r)\notin H, d(r)\notin H.\\
        \end{cases}$$
    
    \paragraph{Decisions of the Followers:}
    {We consider two alternative settings for the decision of the followers.}
    \begin{enumerate}
        \item \textbf{Profit-Maximizing Carriers:} For each $r\in R$ such that $o(r) \in V\setminus H$, if $o(r)$ is allocated to carrier $k\in K$, the carrier observes the offered outsourcing fees $p^r_i$ for all $i \in H$ {and accepts the most profitable one (provided that it results in a non-negative profit), or refuses to serve commodity $r$, otherwise.}

        Similarly, if $d(r)$ is allocated to $k\in K$, the carrier observes the 
        outsourcing fees $q^r_i$, $i \in H$  and accepts the most profitable one (if it produces a non-negative profit), or refuses to serve commodity $r$, otherwise.
        
        The profit functions for routing the first, respectively third leg, for a given $r \in R$, $k \in K$ and outsourcing fees $p$ and $q$ are, therefore, defined as follows:
        \begin{equation}\label{optimal:followers}
            \Pi^k(r,p) = \left[\max_{i \in H} \left\{ p^r_i - \bar{c}^{rk}_{o(r)i}\right\}\right]^+ \qquad \text{and}
            \qquad
            \Gamma^k(r,q) = \left[\max_{i \in H} \left\{ q^r_i - \bar{c}^{rk}_{id(r)} \right\} \right]^+,
        \end{equation}
        where $[\alpha]^+=\max\{0, \alpha\}$.
        
        In the above two functions, when the value of the inner maximum is non-negative, the respective sets of optimal solutions (hubs) are denoted by $I_{k}(r,p)$ and $J_{k}(r,q)$. If the inner maximum in $\Pi^k(r,p)$ or $\Gamma^k(r,q)$ is negative, we define $I_{k}(r,p)=\emptyset$ and $I_{k}(r,p)=\emptyset$, respectively. That is
        \begin{equation}\nonumber
          I_k(r,p) = \begin{cases}
        \arg\max_{i \in H}\left\{p^r_i - \bar{c}^{rk}_{o(r)i}\right\} & \text{if }\max \left\{p^r_i - \bar{c}^{rk}_{o(r)i}: i \in H \right\}\geq 0\\ \emptyset & \text{otherwise,}
        \end{cases} 
        \end{equation}
        \begin{equation}\nonumber
              J_k(r,q) = \begin{cases}
        \arg\max_{i \in H}\left\{q^r_i - \bar{c}^{rk}_{id(r)}\right\} & \text{if }\max \left\{ q^r_i - \bar{c}^{rk}_{id(r)}: i \in H\right\}\geq 0 \\ \emptyset & \text{otherwise}.
        \end{cases}
        \end{equation}

        The problem can now be formally stated as follows:   
        \begin{definition}[\small{Multi-Commodity Flow Problem with Outsourcing Decisions, \SDOD}]\label{def:bilevel}
            The Multi-Commodity Flow Problem with Outsourcing Decisions can be expressed as the following bilevel optimization problem:
            \begin{subequations}
                \begin{align}
                    \qquad&  \max_{\substack{R^*\subseteq R,\,  a:V\setminus H \mapsto  K\cup\{0\} \\
                    p\geq 0,\, q\geq 0}}\sum_{r \in R^*} \left[ b_r- C_r(a, p, q)\right]&& \nonumber
                \end{align}
            \end{subequations}
            \noindent where, for each commodity $r\in R$ the routing costs $C_r(a, p, q)$ are calculated as:
            \begin{equation} \label{eq:functionC}
                C_r(a, p, q)= \begin{cases}
                        \min \{ p_{i}^r+c_{{i}{j}}^r+q_{j}^r:   i\in I_{a(o(r))}(r,p), j\in J_{a(d(r))}(r,q)\} \;   & o(r) \notin  H , \; d(r) \notin  H\\
                       \min \{ p_{i}^r+c_{id(r)}^r:   i\in I_{a(o(r))}(r,p) \} & o(r) \notin  H, \;d(r) \in H\\
                       \min \{  c_{{o(r)}{j}}^r + q_{j}^r: j\in J_{a(d(r))}(r,q) \}  & d(r) \notin  H,  \;o(r) \in H,
                    \end{cases}
            \end{equation}
            and $a(o(r))$ and $a(d(r))$ refer to the carrier the origin, respectively destination of commodity $r$ is allocated to.
        \end{definition}
        In the definition of $C_r(a, p, q)$, the conditions $i\in I_{a(o(r))}(r,p)$ and $j\in J_{a(d(r))}(r,q)$ reveal the bilevel nature of the \SDOD: they enforce that $i$ and $j$ are  \emph{optimal followers' responses} for routing the first, respectively third leg, given the allocation and outsourcing fee decisions of the leader.
        
        The allocated carrier refuses to serve the first, respectively third leg of commodity $r$  if $I_{a(o(r))}(r,p)=\emptyset$ or $I_{a(d(r))}(r,q)=\emptyset$, and hence the overall routing cost must be set to $\infty$.

        We consider an \emph{optimistic} bilevel optimization setting in which the leader chooses the most profitable routing path 
        in case there are  multiple optimal responses of the followers. This is embedded in the definition of $C_r(a, p, q)$. If $|I_{a(o(r))}(r,p)|>1$ or $|J_{a(d(r))}(r,q)|>1$ so the respective carriers can route commodity $r$ through multiple alternative hubs $i$ or $j$ while receiving the maximum profit, the leader will choose to route $r$ through inter-hub arc $(i(r), j(r))\in A^r$ such that the overall routing cost plus outsourcing fees for commodity $r$ is minimized, i.e. $(i(r),j(r)) \in \arg \min \{p_{i}^r+c_{{i}{j}}^r+q_{j}^r: i\in I_{a(o(r))}(r,p), \, j\in J_{a(d(r))}(r,q)\}$.
        
        \item \textbf{Reservation-Price-Oriented Carriers:}  Alternatively, we may assume that the followers react to the leader's offers following a more flexible policy, according to which they are willing to accept any offered outsourcing fee, provided that it meets their reservation price, even if it is not the most profitable one. This policy can be modeled by redefining the set of acceptable connecting hubs for first and third leg offers to:
        \begin{equation}\nonumber
            rI_k(r, {p}) = \left\{i \in H: {p}^r_i - \bar{c}^{rk}_{o(r)i}\geq 0\right\}
        \text{ and }
            rJ_k(r, {q}) = \left\{i \in H: {q}^r_i - \bar{c}^{rk}_{id(r)}\geq 0\right\},
        \end{equation}
        \noindent respectively.  So, for each commodity $r\in R$ the routing costs, now denoted as $rC_r(a, p, q)$, become:
        \begin{equation} \label{eq:functionrC}
            rC_r(a, p, q)= \begin{cases}
                \min \{ p_{i}^r+c_{{i}{j}}^r+q_{j}^r :   i\in rI_{a(o(r))}(r,p), j\in rJ_{a(d(r))}(r,q)\} \;   &
                o(r) \notin  H , \; d(r) \notin  H \\
               \min \{  p_{i}^r+c_{{i}{d(r)}}^r :   i\in rI_{a(o(r))}(r,p) \} & o(r) \notin  H, \;d(r) \in H\\
               \min \{   c_{{o(r)}{j}}^r + q_{j}^r : j\in rJ_{a(d(r))}(r,q) \}  & d(r) \notin  H,  \;o(r) \in H.
            \end{cases}
        \end{equation}
      \noindent This policy leads to the \emph{relaxed Multi-Commodity Flow Problem with Outsourcing Decisions} defined next:        
        \begin{definition}[\small{\emph{relaxed} Multi-Commodity Flow Problem with Outsourcing Decisions, $r$\SDOD}]\label{def:r-bilevel}
            The $r$\SDOD\ can be expressed as the following optimization problem:
                \begin{subequations}
                    \begin{align}
                        \qquad&  \max_{\substack{R^*\subseteq R,\,  a:V\setminus H \mapsto  K \cup\{0\}\\
                        p\geq 0,\, q\geq 0}}\sum_{r \in R^*} \left[ b_r- rC_r(a, p, q)\right],&& \nonumber
                    \end{align}
                \end{subequations}
            where, for each commodity $r\in R$ the routing costs $rC_r(a, p, q)$ are calculated as in \eqref{eq:functionrC}.
        \end{definition}
    \end{enumerate}
    
    The definitions of the above two problems show that the leader faces the following trade-off for finding suitable outsourcing fees.          
    On the one hand, it can be convenient to increase outsourcing fees so as to get more positive responses, which may produce a higher total revenue for the served commodities. On the other hand, high outsourcing fees reduce the leader's net profit, possibly resulting in fewer served commodities.

    In this paper we also introduce variants of \SDOD\ and $r$\SDOD\ in which the outsourcing fees are fixed. Fees can be fixed by an external (regulating) authority, so that the leader only has to choose the allocation of non-hubs to carriers, the commodities to be served, and their routing paths.  The problems are defined as:
    \begin{definition}[\small{\SDOD\ and $r$\SDOD\ with Fixed Outsourcing Fees, \SDOD$_F$ and $r$\SDOD$_F$}]\label{def:bilevel-fixed}
        Let $\bar p^r_i$ and $\bar q^r_i$ be the given outsourcing fees for each $r \in R$, $i \in H$.
        \begin{itemize}
            \item The \SDOD$_F$ can be expressed as the following bilevel optimization problem:
            \begin{subequations}
                \begin{align}
                    \qquad&  \max_{\substack{R^*\subseteq R,\,  a:V\setminus H \mapsto  K\cup\{0\}}} \sum_{r \in R^*} \left[ b_r- C_r(a, \bar p, \bar q)\right]&& \nonumber
                \end{align}
            \end{subequations}
            where the function $C_r(a, \bar p, \bar q)$ is defined as in \eqref{eq:functionC}.
            \item The $r$\SDOD$_F$ can be expressed as the following optimization problem:
            \begin{subequations}
                \begin{align}
                    \qquad&  \max_{\substack{R^*\subseteq R,\,  a:V\setminus H \mapsto  K \cup\{0\}}} \sum_{r \in R^*} \left[ b_r- rC_r(a, \bar p, \bar q)\right]&& \nonumber
                \end{align}
            \end{subequations}
            where the function $rC_r((a, \bar p, \bar q)$ is defined as in \eqref{eq:functionrC}.
        \end{itemize} 
    \end{definition}
    We observe that in both the \SDOD$_F$ and the $r$\SDOD$_F$,  the leader still has to decide on the optimal allocation of non-hub nodes to carriers, the commodities to be served and their routing paths. We also observe that the \SDOD$_F$ is still modeled as an optimistic bilevel problem.
    
    The following properties apply to any of the problems defined above.
    \begin{rem}
        \begin{enumerate}\item[] 
            \item By associating sufficiently large profits to each commodity, we obtain the problem in which all demand must be satisfied, while minimizing the sum of routing and outsourcing cost.
            \item The assumption that the hub network induced by the nodes $H$ is a complete graph with routing costs that satisfy triangle inequality, can be done without loss of generality. Indeed, for each commodity $r \in R$ and each arc $(i,j)$, $i,j \in H$, $i \neq j$, the routing cost $c^r_{ij}$ will correspond to the cost of the shortest $i$-$j$ path in the hub network (see also the example below).
          \item There is an optimal solution in which the routing of each commodity is not split among multiple paths.
        \end{enumerate}
    \end{rem}

    \begin{example}\label{example:1}
    Consider the \SDOD$_F$ example of Figure \ref{fig:1}(a) on a  network with four hub nodes, $H=\{3,\,4,\,6,\,7\}$, and three non-hubs $V\setminus H=\{1,\,2,\,5\}$, with symmetric routing costs $c^r_{ij}$,  the same for all commodities (displayed next to each edge). This backbone network is not complete, so we first make it complete by adding links $(i,j)$ for all distinct pairs $i,j \in H$ with the cost of the shortest $i$-$j$ path, and then we remove \emph{redundant hubs} 6 and 7, i.e., those hubs which are not adjacent to $V\setminus H$ (see Figure \ref{fig:1}(b)). 
    
    We furthermore assume that there are two carriers $K=\{1,\,2\}$, and four commodities with respective origins and destinations $(o(r),\, d(r))$ given by $(1,\, 5)$, $(2,\, 5)$, $(3,\, 5)$, $(4,\, 5)$, all of them with demand of one unit and associated revenue of 100, i.e., $w^r=1$, $b^r=100$, $r\in R=\{1, \dots, 4\}$.
    
    The routing costs for the leader and the reservation prices for the carriers are given in Table \ref{tab:ex1}. For ease of presentation, we assume in this example that these values  only depend on the arc but not on the commodity itself. As indicated above this will not be the case in general.

    \begin{figure}[htbp]%
        \centering
        \subfigure[Original network, considering non-complete hub network with possible redundant hubs.]{{\includegraphics[page=2, width=.40\textwidth, keepaspectratio]{example1bis2.pdf}}}%
        \qquad
        \subfigure[Network after making the hub network complete and removing redundant hubs.]{{\includegraphics[page=1, width=.40\textwidth, keepaspectratio]{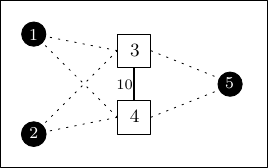}}}%
        \caption{Network for Example \ref{example:1} with set of commodities $\{(1, 5), (2, 5), (3, 5), (4, 5)\}$.}\label{fig:1}
        \footnotesize{Squares indicate hub nodes from $H$, solid lines the interhub links with their corresponding leader routing costs, and dotted lines potential access or distribution arcs to/from the backbone network.}
    \end{figure}
    
    \begin{table}[htbp]
        \caption{Routing costs and reservation prices.\label{tab:ex1}}
        {\begin{tabular}{ccccccccc}
            \multicolumn{2}{c}{} & $(1, 3)$ & $(1, 4)$ & $(2, 3)$ & $(2, 4)$ & $(3, 4)$ & $(3, 5)$ & $(4, 5)$ \\ \hline \hline
            \multicolumn{2}{c} {Routing costs $c^r_{ij}$ (Leader)} & -- & -- & -- & -- & 10 & -- & -- \\
            \hdashline        \multirow{2}{*}{ Reservation prices $\bar c^{rk}_{ij}$
(Carriers)} & k=1 & 10 & 40 & 25 & 40 & -- & 40 & 40 \\
             & k=2 & 30 & 40 & 20 & 40 & -- & 45 & 20 \\ \hline
        \end{tabular}}\\
        \footnotesize{The same for all commodities $r\in R$}
    \end{table}%
    Table \ref{tab:ex1_2} shows fixed outsourcing fees, $\overline p_i^r$ and $\overline q_i^r$,  for all $r \in R,\, i \in H$. For $r=4$, $(o(4), d(4))=(4, 5)$, since its origin is a hub, only the outsourcing fees for the third leg, namely $\overline q_i^r$, are required.
    
    Table \ref{tab:ex1_2} also shows the carriers' responses for the \SDOD\ to the offers  $\overline p_i^r$ and $\overline q_i^r$. There are two rows associated with each carrier. One with its response (YES or NO) to the fixed outsourcing fees $\overline p^r_i$ and $\overline q^r_i$, and another one with the carriers reservation price of the arc against which the fees should be contrasted.

    For example, the first block of columns in Table \ref{tab:ex1_2} shows the followers' responses for the first leg if fixed outsourcing fees are $\overline{p}^1_3=30$, $\overline{p}^1_4=50$. That is, the carrier allocated to non-hub node 1 would receive 30 or 50 euros depending on whether the commodity is routed through access arc $(1, 3)$ or $(1, 4)$, respectively.     

    Even if the offer is higher for the latter connection, only carrier 2 would accept it. The reason is that,  by accepting the offer of 30 euros for access arc $(1, 3)$, carrier 1 would obtain a higher net profit namely $30-10=20$ than the one it would obtain if it accepted the 50 euros offer for access arc $(1, 4)$, which would be of $50-40=10$ euros. Offers that are not accepted because, even if they are profitable, they do not produce the highest net profit to the carrier are marked with an asterisk.
    
    Similarly, for $\overline{q}^1_3=45$, $\overline{q}^1_4=40$, carrier 2 would accept to route the third leg of commodity $r=1$ through hub $4$, achieving the net profit of $40 - 20 = 20$. It can also be observed that some offers would not be accepted neither by carrier 1 nor 2; for instance,  $\overline{q}^3_3=30$. 
    \begin{table}[H]
    \caption{Carriers' responses to fixed outsourcing fees offers (\SDOD$_F$).\label{tab:ex1_2}}
    {\begin{tabular}{ccccccccccccc}
        \multicolumn{2}{c}{$r$} & \multicolumn{4}{c}{$(1, 5)$} & \multicolumn{4}{c}{$(2, 5)$} & \multicolumn{2}{c}{$(3, 5)$} & $(4, 5)$ \\ \hline \hline
        \multicolumn{2}{c}{} & \multicolumn{2}{c}{$p^r_i$} & \multicolumn{2}{c}{$q^r_i$} & \multicolumn{2}{c}{$p^r_i$} & \multicolumn{2}{c}{$q^r_i$} & \multicolumn{2}{c}{$q^r_i$} & $\overline{q}^r_i$ \\
        \multicolumn{2}{c}{$i(r)/j(r)$} & 3 & 4 & 3 & 4 & 3 & 4 & 3 & 4 & 3 & 4 & 4 \\ \hline
        \multicolumn{2}{c}{Fixed outsourcing fee} & 30 & 50 & 45 & 40 & 15 & 35 & 40 & 30 & 30 & 20 & 40 \\ \hline
        \multirow{4}{*}{Response} & \multirow{2}{*}{$k=1$} & YES & NO$^*$ & YES & NO$^*$ & NO & NO & YES & NO & NO & NO & YES \\
         &  & (10) & (40) & (40) & (40) & (25) & (40) & (40) & (40) & (40) & (40) & (40) \\ \cline{2-13} 
         & \multirow{2}{*}{$k=2$} & NO$^*$ & YES & NO$^*$ & YES & NO & NO & NO & YES & NO & YES & YES \\
         &  & (30) & (40) & (45) & (20) & (20) & (45) & (45) & (20) & (45) & (20) & (20) \\ \hline
    \end{tabular}}%
    \end{table}%

    Note that, in this example, the carrier allocation of node 2 is irrelevant because no offer related to routing the first leg of commodity $(2, 5)$, which is the only one involved in that node, is accepted. For the fixed outsourcing fees as depicted in  Table \ref{tab:ex1_2}, Table \ref{tab:ex1_23} analyzes the possible  \SDOD\  outcomes depending on the allocation of carriers to non-hub nodes 1 and 5.
    
    The table shows the routing paths and associated outsourcing plus routing costs for the cases where the involved carriers would accept the corresponding offers. This is why the entries of commodity $(3, 5)$ when $a(5)=1$ are empty. We observe, that for the allocation $a(1)=2$ and $a(5)=1$ commodity $(1,\,5)$ would not be served, because even if both carriers would accept the offers of the leader, the overall outsourcing fee plus routing cost (105) exceeds the revenue (100) for serving it.
    
    Indeed, if the outsourcing fees are fixed as in Table \ref{tab:ex1_2} the best \SDOD\  solution for the given fixed outsourcing fees is to allocate carrier 2 to node 5 and carrier 1 to node 1 (i.e. $a(1)= 1,\, a(5)=2$). Then, the leader is able to route three out of the four commodities with a resulting net profit of $150$ euros.
    \begin{table}[htbp]
    \caption{Potential allocation of carriers to non-hub nodes.\label{tab:ex1_23}}
    \resizebox{\textwidth}{!}{
    \begin{tabular}{cccccccccc}
         &  & \multicolumn{2}{c}{$r=1$} & \multicolumn{2}{c}{$r=3$} & \multicolumn{2}{c}{$r=4$} & \multirow{3}{*}{$R^*$} & \multirow{2}{*}{Net Profit} \\
          &  & \multicolumn{2}{c}{$(o(r), d(r))=(1, 5)$} & \multicolumn{2}{c}{$(o(r), d(r))=(3, 5)$} & \multicolumn{2}{c}{$(o(r), d(r))=(4, 5)$} &  &  \\            
        $a(5)$ & $a(1)$ & path & $C_r(a,p,q)$ & path & $C_r(a,p,q)$ & path & $C_r(a,p,q)$ &  & $\sum_{r \in R^*} \left[ b_r- C_r(a, p, q)\right]$ \\ \hline \hline
        \multirow{2}{*}{1} & 1 & 1-3-5 & 30+45 &  &  & \multirow{2}{*}{4-5} & \multirow{2}{*}{40} & $\{1,\,4\}$ & 200-(75 + 40)=85 \\
         & 2 & 1-4-3-5 & 50+10+45 $>$100 &  &  &  &  & $\{4\}$ & 100-40= 60 \\
        \hdashline            \multirow{2}{*}{2} & 1 & 1-3-4-5 & 30+10+40 & \multirow{2}{*}{3-4-5} & \multirow{2}{*}{10 + 20} & \multirow{2}{*}{4-5} & \multirow{2}{*}{40} & \multirow{2}{*}{$\{1,\,3,\,4\}$} & 300-(80+30+40)=150 \\
         & 2 & 1-4-5 & 50 +40 &  &  &  &  &  & 300-(90+30+40)=140\\
         \hline
    \end{tabular}}\\
    \footnotesize{Service paths and net profits for the solutions with carriers' positive responses for both the first and third legs for fixed outsourcing fees as in Table \ref{tab:ex1_2}.}
    \end{table}%

    This example also illustrates the bilevel nature of  \SDOD$_F$. {If we consider $r$\SDOD$_F$, so that} the condition that \textit{accepted offers must be optimal} for the carriers is substituted by the condition that \textit{accepted offers must result in a non-negative profit} for the carriers, then routes involving profitable but non optimal offers (marked with an asterisk in Table \ref{tab:ex1_2}) should also be considered. In that case, nothing would change for commodities $(3, 5)$ and $(4, 5)$. However, for commodity $(1, 5)$ the allocation $a(1)=a(5)=2$ would produce positive answers of carrier 2 for $\overline{p}^{(1,\,5)}_3$ and $\overline{q}^{(1,\,5)}_4$. Then, the path $1-3-5$ would be accepted, since its outsourcing fees plus routing cost of $75=30+45$ is smaller than the total cost of the original $1-4-5$ route $(90=50+40)$. Again, the same three commodities are served but the route of commodity $(1, 5)$ would be $1-3-5$ with a resulting net profit of $155=300-(75+30+40)$.
        
    Hence, $r$\SDOD$_F$ produces an overestimation of the leader's profit in \SDOD$_F$,  given that any offer acceptable for \SDOD$_F$ would also be acceptable for $r$\SDOD$_F$.
    
    However, if the leader could choose the outsourcing fees, the solution would be very different.  Table \ref{tab:ex1_opt} shows outsourcing fees for an optimal solution for the \SDOD\ problem whose routing cost and reservation prices are provided in Table \ref{tab:ex1}. Specifically, we show the optimal outsourcing fees and the respective carrier responses that the leader would anticipate before making the allocation and routing decisions.
    \begin{table}[H]
    \caption{Optimal outsourcing fees and carriers responses for \SDOD.\label{tab:ex1_opt}}
    \resizebox{\textwidth}{!}{\begin{tabular}{ccccccccccccc}
        \multicolumn{2}{c}{$r$} & \multicolumn{4}{c}{$(1, 5)$} & \multicolumn{4}{c}{$(2, 5)$} & \multicolumn{2}{c}{$(3, 5)$} & $(4, 5)$ \\ \hline \hline
        \multicolumn{2}{c}{} & \multicolumn{2}{c}{$p^r_i$} & \multicolumn{2}{c}{$q^r_i$} & \multicolumn{2}{c}{$p^r_i$} & \multicolumn{2}{c}{$q^r_i$} & \multicolumn{2}{c}{$q^r_i$} & $q^r_i$ \\
        \multicolumn{2}{c}{$i(r)/j(r)$} & 3 & 4 & 3 & 4 & 3 & 4 & 3 & 4 & 3 & 4 & 4 \\ \hline
        \multicolumn{2}{c}{Outsourcing fee} & 10 & 0 & 0 & 20 & 20 & 0 & 0 & 20 & 0 & 20 & 20 \\ \hline
        \multirow{4}{*}{Response} & \multirow{2}{*}{$k=1$} & YES & NO & NO & NO & NO & NO & NO & NO & NO & NO & NO \\
         &  & (10) & (40) & (40) & (40) & (25) & (40) & (40) & (40) & (40) & (40) & (40) \\ \cline{2-13} 
         & \multirow{2}{*}{$k=2$} & NO & NO & NO & YES & YES & NO & NO & YES & NO & YES & YES \\
         &  & (25) & (40) & (45) & (20) & (20) & (45) & (45) & (20) & (45) & (20) & (20) \\ \hline
    \end{tabular}}\\
    \footnotesize{Considering data provided in Table \ref{tab:ex1}}
    \end{table}%
    
    The optimal solution for this example shows that by considering the outsourcing fees as decision variables, we are able to serve all four commodities and maximize the leader's profit, while maintaining positive responses from the carriers. The optimal allocations are defined with $a(1) = 1$ and $a(2) = a(5) = 2$. The total leader's profit increases from 150 to $260$ euros.
   \end{example}

\section{Bilevel Formulation and a Single-Level MINLP Reformulation} \label{sec:bilevelModel}
    In this section we first provide a bilevel problem formulation for the more general \SDOD\ and we derive a single-level MINLP reformulation by exploiting strong duality conditions applied to the lower level problems. As we will see, for \SDOD$_F$ the resulting reformulation is a MILP. {We also show how to adapt these formulations for $r$\SDOD\ and $r$\SDOD$_F$}.
    
    \subsection{Bilevel MILP formulation for the \SDOD}
        We define the following sets of decision variables:
        \begin{itemize}
            \item For the upper level problem, associated with the decisions of the leader:
            \begin{itemize}
                \item[] $s^r\in\{0,\, 1\}$, $r\in R$, for the commodities that are served.
                \item[] $a_{ik}\in\{0,\, 1\}$, $i\in V\setminus H$, $k\in K$, for the allocation of non-hub nodes to carriers.
                \item[] $x_{ij}^r\in\{0,\, 1\}$, $(i, j)\in A^r$, $r\in R$, for the commodities that are served using the interhub arc $(i,j)$.
                \item[] $p_i^r\geq 0$: outsourcing fee for access arc $(o(r), i)$ for commodity $r\in R$.
                \item[] $q_i^r\geq 0$: outsourcing fee for distribution arc $(i, d(r))$ for commodity $r\in R$.
            \end{itemize}
            \item For the lower level problem, associated to the carrier $k\in K$:
            \begin{itemize}
                \item[] $f_{i}^{rk}\in\{0,\, 1\}$, $i\in H$, $r\in R$. $f_{i}^{rk}$ takes the value 1 if and only if $o(r) \in V \setminus H$ is allocated to carrier $k\in K$ and $k$ accepts the offer for routing the first leg of commodity $r$ through access arc $(o(r), i)$. That is, $f_{i}^{rk}=1
                \Rightarrow i\in I_{k}(r, p)$ for the outsourcing fee $p_i^r$.
                \item[] $t_{i}^{rk}\in\{0,\, 1\}$, $i\in H$, $r\in R$. $t_{i}^{rk}$ takes the value 1 if  and only if $d(r) \in V \setminus H$ is allocated to carrier $k\in K$ and $k$  accepts the offer for routing the third leg of commodity $r$ through distribution arc $(i, d(r))$. That is, $t_{i}^{rk}=1
                \Rightarrow i\in J_{k}(r, q)$ for the outsourcing fee $q_i^r$.
            \end{itemize}
        \end{itemize}
        The formulation of the {multicommodity flow problem with outsourcing decision} for the leader is:
        \begin{subequations}\label{mod:leader}
            \begin{align}
                \quad\quad  \max & \sum_{r\in R}b^rs^r- \sum_{r\in R}\sum_{(i,\,j)\in A^r} \left [p_i^r+q_j^r+ c_{ij}^r\right]x_{ij}^r && \label{of}\\  
                 \mbox{s.t. }
                & \sum_{k\in K}a_{ik} \leq 1 && i \in  V\setminus H  \label{const:allocate}\\
                & \sum_{(i,\,j)\in A^r}x_{ij}^r \leq \sum_{k\in K}f_{i}^{rk}\, &&     {r\in R}, i\in H, \text{ s.t. } o(r)\ne i \label{const:serve-origin}\\
                & \sum_{(j,\,i)\in A^r}x_{ji}^r \leq  \sum_{k\in K}t_{i}^{rk}\, &&     {r\in R}, i\in H \text{ s.t. } d(r)\ne i\label{const:serve-destination}\\
                & \sum_{(i,\,j)\in A^r} x^r_{ij} = s^r && r\in R \label{one-interhub}\\
                & {f^{rk}\in\arg\max F^{rk}(a_{o(r)k},\, p)} && r\in R, k\in K \text{ s.t. } o(r)\notin H\label{a-first}\\
                & {t^{rk}\in\arg\max T^{rk}(a_{d(r)k},\, q)} && r\in R, k\in K \text{ s.t. } d(r)\notin H\label{a-third}\\
                & p_i^r, q_i^r\geq 0 \, && i\in H,\, r\in R  \\
                & x_{ij}^r\in \{0,1\} \, &&r\in R, \, (i,j)\in A^r\label{const:binary}\\
                & s^r\in \{0,1\}\, && r\in R \\
                & a_{ik}\in\{0, 1\} && i\in V\setminus H,\, k\in K,
                \label{const:domain}
            \end{align}
        \end{subequations}
        where $F^{rk}(a_{o(r)k}, p)$ and $T^{rk}(a_{d(r)k}, q)$ denote the subproblems for carrier $k\in K$ for the first and third legs of commodity $r\in R$, defined by \eqref{subeq:F} and \eqref{subeq:T}, respectively.
        
        Constraints \eqref{const:allocate} impose that only non-hub nodes can be allocated to carriers, and that each non-hub node can be allocated to at most one carrier. Constraints \eqref{const:serve-origin} and  \eqref{const:serve-destination} relate the interhub arcs (possibly loops) used in the routing of served commodities with their first and third legs, respectively.  Constraints \eqref{one-interhub}  identify served commodities by associating them to a routing arc (possibly a loop) of the backbone network, which, by  Constraints \eqref{const:serve-origin}-\eqref{const:serve-destination} will be the one connecting their first and third legs. By Constraints \eqref{a-first}-\eqref{a-third}  the first and third legs that are accepted by the carriers correspond to optimal decisions of the carriers to the leader's outsourcing fees and to the allocation to carriers to non-hub nodes to carriers.

        The first term of the objective function is the total revenue for the served commodities, whereas the second term refers to routing costs and outsourcing fees. Note that this second term is bilinear, hence, the overall formulation is a bilevel MINLP with non-convex objective function and lower-level problems being linear programs as shown below. 

        \paragraph{The carriers' subproblems.}\label{scr:carriers}
        Given a first-level solution, with outsourcing fees $\overline p$, $\overline q$, and allocation of non-hub nodes $\overline a$,  each carrier $k\in K$ must identify the commodities' first and third leg offers it will accept, among the ones starting at or arriving to the nodes allocated to, as stated in \eqref{optimal:followers}. The decisions on the different commodities are independent from each other. To determine the most profitable access arc for the first leg of carrier $k$ for commodity $r$ we solve the following problem:
        \begin{subequations}
        \label{subeq:F}
            \begin{align}
                 F^{rk}(\overline a_{o(r)k}, \overline p)\quad = \max & \sum_{i\in H} \rho^{rk}_i f_{i}^{rk}  \label{of_sub_f}\\         
                 \mbox{s.t. }
                &  \sum_{i\in H}f_{i}^{rk} \leq \overline a_{o(r)k} &&   \label{const:single allocation-f}\\
                & f_{i}^{rk} \in\{0, 1\} && i\in H,
            \end{align}
        \end{subequations}
        where $\rho_{i}^{rk}=\overline p_i^r- \overline c_{o(r)i}^{rk}$, $i\in H$.
        
        We notice that the optimal solution value of problem \eqref{subeq:F}  is always non--negative. If $\rho_{i}^{rk}<0$ for all $i\in H$ then $f^{rk}_i=0$ for all $i \in H$. Indeed, $F^{rk}(\overline a_{o(r)k}, \overline p)=\Pi^k(r, \overline p)$ when $\overline a_{o(r)k}=1$ and 0, otherwise.
    
        Analogously, the third leg subproblem of carrier $k$ for commodity $r$ is:
        \begin{subequations}
        \label{subeq:T}
            \begin{align}
                T^{rk}(\overline a_{d(r)k}, \overline q)\quad = \max & \sum_{i\in H} \gamma_{i}^{rk} t_{i}^{rk}  \label{of_sub_t}\\
                 \mbox{s.t. }
                &  \sum_{i\in H}t_{i}^{rk} \leq \overline a_{d(r)k} && \label{const:single allocation-t}\\
                & t_{i}^{rk} \in\{0, 1\} && i\in H,
            \end{align}
        \end{subequations}
        where $\gamma_{i}^{rk}=\overline q_i^r- \overline c_{id(r)}^{rk}$, $i\in H$. 
        We also have that $T^{rk}(\overline a_{d(r)k}, \overline q)=\Gamma^k(r, \overline q)$ when $\overline a_{o(r)k}=1$ and 0, otherwise.

        \begin{rem}\label{rem:dual}
            It is easy to see that the integrality conditions on variables $t_{i}^{rk}$ and $q_{i}^{rk}$ can be relaxed  and replaced by non-negativity constraints.  Then, the optimal values of the carriers subproblems, $F^{rk}(\overline a_{o(r)k}, \overline q)$ and  $T^{rk}(\overline a_{d(r)k}, \overline q)$, can be obtained alternatively \emph{via} their respective dual problems namely    
            \begin{subequations}\label{subeq:DF}
                \begin{align}
                    F^{rk}(\overline a_{o(r)k}, \overline p)\quad = \min \quad & \overline a_{o(r)k} u_{rk}  \label{of_sub_f-dual}\\
                    \mbox{s.t. }
                    &  u_{rk} \geq \rho_{i}^{rk} \quad &&  i\in H  \label{const:-f-dual}\\
                    & u_{rk} \geq 0 && i\in H,
                \end{align}
            \end{subequations} 
            \begin{subequations}\label{subeq:DT}
                \begin{align}
                    T^{rk}(\overline a_{d(r)k}, \overline q)\quad = \min \quad & \overline a_{d(r)k} v_{rk}  \label{of_sub_t-dual}\\
                    \mbox{s.t. }
                    &  v_{rk} \geq \gamma_{i}^{rk} \quad &&  i\in H  \label{const:-t-dual}\\
                    & v_{rk} \geq 0 && i\in H,
                \end{align}
            \end{subequations}
            whose optimal solutions are $u^*_{rk}=\left[\max_{i\in H}\rho_{i}^{rk}\right]^+$ and $v^*_{rk}=\left[\max_{i\in H}\gamma_{i}^{rk}\right]^+$, respectively. Hence, we have $F^{rk}(\overline a_{o(r)k}, \overline p)=\overline a_{o(r)k}\left[\max_{i\in H}\rho_{i}^{rk}\right]^+$ and $T^{rk}(\overline a_{o(r)k}, \overline p)=\overline a_{d(r)k}\left[\max_{i\in H}\gamma_{i}^{rk}\right]^+$.  
        \end{rem}
    \subsection{Single-level MI(N)LP Reformulation}
        In order to guarantee that first and third legs correspond to optimal carriers' decisions, we substitute \eqref{a-first}-\eqref{a-third} by  primal-dual optimality conditions on the carriers' decisions. For this, as explained in Remark \ref{rem:dual}, we introduce additional dual variables $u_{rk}$, $v_{rk}$, $r\in R$, $k\in K$ associated with the constraints \eqref{const:single allocation-f} and \eqref{const:single allocation-t} of the carriers subproblems, respectively. Then, the single-level reformulation of \SDOD\ reads as follows:
        \begin{subequations} \label{subeq:MINLPreformulation}
            \begin{align}
              \qquad \qquad \max &  \sum_{r\in R} b^r s^r-\sum_{r\in R}\sum_{(i,\,j)\in A^r}[ c_{ij}^r +   p^{r}_{i} +   q^{r}_{j}]  x_{ij}^r  && \label{of:MINLP}\\                 \mbox{s.t. }
                & \sum_{k\in K}a_{ik} \leq 1 && i\in V\setminus H\label{const:allocate-fixed}\\
                &  \sum_{k\in K}f_{i}^{rk} {\ge} \sum_{(i,\,j)\in A^r}x^r_{ij}&& r\in R, i\in H \text{ s.t. } o(r)\ne i\label{access-route-fixed1}\\
                &  \sum_{k\in K}t_{i}^{rk} {\ge}  \sum_{(j,\,i)\in A^r}x^r_{ji}&& r\in R, i\in H \text{ s.t. } d(r)\ne i\label{distribution-route-fixed1}\\
                &\sum_{(i,\,j)\in A^r} x^r_{ij}= s^r && r\in R \label{one-interhub-fixed} \\
                 &  \sum_{i\in H}f_{i}^{rk} \leq a_{o(r)k} && r\in R, k\in K \label{const:or_allocation_MINLP}\\
                &  \sum_{j\in H}t_{j}^{rk} \leq a_{d(r)k} && r\in R, k\in K \label{const:dr_allocation_MINLP}\\
                & u_{rk}\geq  \rho^{rk}_i&& i\in H,\, r\in R,\, k\in K \label{dual-feas-first}\\
                & v_{rk}\geq \gamma^{rk}_i && i\in H,\, r\in R,\, k\in K \label{dual-feas-third}\\
                & u_{rk}a_{o(r)k} \le  \sum_{i\in H} \rho^{rk}_i\,f_{i}^{rk}\qquad && r\in R,\, k\in K\label{primal-dual-first}\\
                & v_{rk}a_{d(r)k} \le   \sum_{i\in H} \gamma^{rk}_i\,t_{i}^{rk}\qquad && r\in R,\, k\in K \label{primal-dual-third}\\     
                & s^r\in\{0, 1\} &&r\in R, \\
                  & p_i^r, q_i^r\geq 0 \, && i\in H,\, r\in R  \\
                & f_{i}^{rk}, t_{i}^{rk} \ge 0  && i\in H, r\in R, k \in K\\
                & a_{ik} \in\{0, 1\}&& i\in V\setminus H, k\in K\\
                & x_{ij}^r\geq 0&& (i,\,j)\in A^r, r\in R\\
                & u_{rk}, v_{rk}\geq 0 &&  r\in R, k\in K. &&
            \end{align}
        \end{subequations}
        where,  $\rho^{rk}_i= p_i^r- \overline c_{o(r)i}^{rk}$, and $\gamma^{rk}_i=  q_i^r- \overline c_{id(r)}^{rk}$, for all $i\in H$, $r\in R$, $k\in K$.
        
        Constraints\eqref{const:or_allocation_MINLP} and \eqref{const:dr_allocation_MINLP} guarantee that, for all $r\in R$, $k\in K$, variables $f_{i}^{rk}$, and $t_{i}^{rk}$ determine feasible solutions for $F^{rk}(a_{o(r)k},   p)$ and $T^{rk}(a_{d(r)k},   q)$, respectively, whereas Constraints \eqref{dual-feas-first}-\eqref{dual-feas-third} guarantee the feasibility of the dual variables $u_{rk}$ and $v_{rk}$. Finally, Constraints \eqref{primal-dual-first} and \eqref{primal-dual-third} impose that the objective function values of $F^{rk}(a_{o(r)k},   p)$ and $T^{rk}(a_{d(r)k},   q)$ coincide with those of their respective duals, thus guaranteeing the optimality of both the primal and dual solutions.
        \begin{rem}\label{rem:reformulation}
            For \SDOD$_F$, the model \eqref{subeq:MINLPreformulation} is an MILP reformulation. Indeed, taking into account Remark \ref{rem:dual}, we can overcome the bilinear terms in the left-hand-sides of \eqref{primal-dual-first}-\eqref{primal-dual-third} by projecting out $u$ and $v$ variables and  substituting constraints \eqref{dual-feas-first}-\eqref{primal-dual-third} by the following:
            \begin{subequations} \label{subeq:primaldual}
                \begin{align}
                    & a_{o(r)k}\left[\max_{i\in H}\rho_{i}^{rk}\right]^+ \le 
                    \sum_{i\in H} \rho^{rk}_i\,f_{i}^{rk}\qquad && r\in R,\, k\in K\label{primal-dual-first-lin}\\
                    & a_{d(r)k}\left[\max_{i\in H}\gamma_{i}^{rk}\right]^+ \le  \sum_{i\in H} \gamma^{rk}_i\,t_{i}^{rk}\qquad && r\in R,\, k\in K. \label{primal-dual-third-lin} \end{align}
            \end{subequations}
            where, for the fixed outsourcing fees $\overline p$ and $\overline q$, the values of $\rho$ and $\gamma$ can be precalculated, and hence  \eqref{subeq:primaldual} can be turned into linear constraints.  
        \end{rem}
        The following criteria can be applied to eliminate decision variables from model \eqref{subeq:MINLPreformulation} for \SDOD$_F$:
        \begin{itemize}
            \item 
            $x_{ij}^r=0$ for all $r\in R$, $(i,j)\in A^r$ such that $\overline p^r_i+\overline q^r_j+c^r_{ij}\geq b^r$.
            \item 
            $f_{i}^{rk}=0$ for all $i\in H$, $r\in R$, $k\in K$ such that $\rho^{rk}_i<0$. \item $t_{i}^{rk}=0$ for all $i\in H$, $r\in R$, $k\in K$ such that $\gamma^{rk}_i<0$.
        \end{itemize}
        \begin{rem}\label{rem:relaxation}
            {Single-level formulations for $r$\SDOD\ and $r$\SDOD$_F$ can be obtained from those of \SDOD\ and \SDOD$_F$, respectively, by replacing constraints \eqref{subeq:primaldual} by the following ones, where it is assumed that outsourcing fees offered by the leader will be accepted, as long as they result in a non-negative profit for the follower. Depending on the case, the outsourcing fees would be decision variables or fixed:}
            \begin{subequations} \label{subeq:primaldual-relax}
                \begin{align}
                    & \sum_{i\in H} \rho_{i}^{rk}f_{i}^{rk}\geq 0&& r\in R, k\in K \\
                    &\sum_{i\in H} \gamma^{rk}_it_{i}^{rk}\geq 0&& r\in R, k\in K,
                \end{align}
            \end{subequations}  
            As mentioned, and also illustrated in Example \ref{example:1}, these are relaxations of their respective  bilevel problems.                 
        \end{rem}

\section{Solution Properties and Problem Complexity}\label{sec:properties}
    Below we discuss some properties of optimal solutions and show that {all four problems, namely, \SDOD, $r$\SDOD, \SDOD$_F$, and $r$\SDOD$_F$} are NP-hard.     Contrary to the intuition that \SDOD\ should be more difficult to solve than \SDOD$_F$ (since there are less decisions to be made for \SDOD$_F$),
    we will see that (thanks to certain properties of optimal solutions) the \SDOD\ is a \textit{soft} bilevel problem, as it will be possible to find an optimal solution by solving an auxiliary single-level problem, in which outsourcing fees are not explicitly stated. Surprisingly, even though problems $r$\SDOD\ and $r$\SDOD$_F$ are relaxations of the original setting, they remain in the same complexity class as \SDOD.

    \subsection{Optimal outsourcing fees for given non-hub allocations, and first and third legs}
        Even if, in principle, outsourcing fees, $p^r_i$, $q^r_i$, may take continuous non-negative values, we next show that there is an optimal solution where the values of these variables can be restricted to a finite set of ${\cal{O}}(|R|\cdot |H|\cdot |K|)$ elements. Moreover, we will see that an optimal set of outsourcing fees, $\overline p$, $\overline q$ can be easily found, for any solution when its allocation along with its first and third legs are given.
        
        Let $\bar f$, $\bar t$, and $\bar a$ satisfying Constraints \eqref{const:allocate-fixed}-\eqref{const:dr_allocation_MINLP} be given. 
        For each $r\in R$, let $k(r), l(r) \in K$ be the carriers such that $\bar a_{o(r)k(r)}=1$ and $\bar a_{o(r)l(r)}=1$, respectively. When $\sum_{i\in H}\overline f^{r k(r)}_i= 1$, let $i(r)\in H$ be the hub such such that $\overline f^{rk(r)}_{i(r)}= 1$. 
        Similarly, when $\sum_{i\in H}\overline t^{r l(r)}_i= 1$, let $j(r)\in H$ be the hub such that $\overline t^{rk(r)}_{j(r)}= 1$.

        Consider the following outsourcing fees:
        For each $r\in R$, such that $o(r)\notin H$, let:
        \begin{equation}
            \overline p_{i}^{r}= \begin{cases}
                              \overline c^{rk(r)}_{o(r)i} & \text{ if } i=i(r)\\
                                0 & otherwise
                            \end{cases},\qquad i\in H.
            \label{best_p}
        \end{equation}
        
        Similarly, for each $r\in R$, such that $d(r)\notin H$, let:
        \begin{equation}
            \overline q_{i}^{r}= \begin{cases}
                              \overline c^{rl(r)}_{id(r)} & \text{ if } j=j(r)\\ 
                                0 & otherwise
                            \end{cases},\qquad i\in H.
            \label{best_q}
        \end{equation}
        
        With the above fees, for all $r\in R$ we have:
        \begin{equation*}
            {\overline \rho_{i}^{rk}=}\overline p_i^r- \overline c_{o(r)i}^{rk}=\begin{cases}0 & k=k(r),\, i=i(r)\\
            - \overline c_{o(r)i}^{rk} & otherwise \end{cases}, \qquad k\in K,\, i\in H,
        \end{equation*}
        \begin{equation*}
            {\overline\gamma_{i}^{rk}=}\overline q_i^r- \overline c_{id(r)}^{rk}=\begin{cases}0 & k=l(r),\, i =j(r)\\
            - \overline c_{id(r)}^{rk} & otherwise \end{cases}, \qquad k\in K,\, i\in H.
        \end{equation*}
        
        Hence, for the given $\overline f$, $\overline t$, $\overline a$, the outsourcing fees obtained according to \eqref{best_p} and \eqref{best_q}, not only satisfy the non-negativity conditions \eqref{subeq:primaldual-relax}, but they also satisfy the bilevel optimality constraints \eqref{subeq:primaldual}. Moreover, these fees take the smallest possible values that, in each case, produce positive responses from the involved carriers for both $r$\SDOD\ and \SDOD. As a consequence, we have the following result: 
        \begin{prop}\label{propo2}
            Let $\overline f$, $\overline t$, and $\overline a$ be given first and third legs, and non-hub allocations 
            satisfying \eqref{const:allocate-fixed}-\eqref{const:dr_allocation_MINLP}. Then, the outsourcing fees $\overline p$ and $\overline q$ obtained according to \eqref{best_p} and \eqref{best_q} are optimal for  $\overline f$, $\overline t$, and $\overline a$, for both $r$\SDOD\ and \SDOD.
        \end{prop}

        Proposition \ref{propo2} has two main consequences. The first one is that, for both \SDOD\ and $r$\SDOD, there is an optimal solution with outsourcing fees $p^r_i\in P^r_i : = \{\overline c^{rk}_{o(r)i}: k\in K\}\cup\{0\}$ and $q^r_i\in Q^r_i : =  \{\overline c^{rk}_{id(r)}: k\in K\}\cup\{0\}$.
        Therefore, in formulation \eqref{mod:leader} the conditions $p^r_i\geq 0$ and $q^r_i\geq 0$ can be substituted by $p^r_i\in P^r_i$ and $q^r_i\in Q^r_i$, respectively.
         
        The second consequence of Proposition \ref{propo2} is that optimal outsourcing fees can be determined according to \eqref{best_p} and \eqref{best_q}, from optimal non-hub allocations, and first and third leg decisions. 

        If an optimal \SDOD\ solution routes commodity $r\in R$ using as first and third legs the ones associated with $\overline f^{rk(r)}_{i(r)}$ and $\overline t^{rl(r)}_{j(r)}$, respectively, the overall outsourcing plus routing cost incurred by the leader will be:
        \begin{equation*}
            C^{i(r)j(r)}_{k(r)l(r)} = \overline p^{r}_i+c^{r}_{i(r),j(r)}+\overline q^r_j=\overline c_{o(r)i(r)}^{rk(r)}+c^{r}_{i(r)j(r)}+\overline c_{j(r)d(r)}^{rl(r)},
        \end{equation*}
        whereas the outsourcing fees associated with any other first leg $\overline f^{rk}_{ik}$ or third leg $\overline t^{rk}_{ik}$ can be set to zero. Indeed, the above comment applies to $r$\SDOD\ as well. 

        The above observation indicates that, even if, in principle, outsourcing fees are independent on the carriers who would route the first and third legs, optimal outsourcing fees will, in fact, depend on who are the carriers in charge of access (first) and distribution (third) legs. By setting the outsourcing fees as in \eqref{best_p} and \eqref{best_q},  the followers will be incentivized to accept the offer which is the most profitable for the leader. This is quite intuitive for the problem at hand and allows us to model the \SDOD\ as a problem that considers only the leader costs along with allocation decisions, where the costs (outsourcing fees) of the first and third legs are precisely the carriers reservation prices of the corresponding legs.

        Furthermore, if the carriers $k(r)$ and $l(r)$ are known for commodity $r\in R$, then the best interhub arc for routing it, $(i(r), j(r))$, can be determined \textit{a priori} by finding:
        \begin{equation*}
            (i(r), j(r))\in\arg\min\{\overline c^{rk(r)}_{o(r)i} +\overline c^{rl(r)}_{jd(r)}+ c_{ij}^r: i,j\in H\}.
        \end{equation*}
        
        For any pair $(k,l)$ of carriers potentially serving the first, respectively third leg of commodity $r$, we can precompute in $O(|H|^2)$ time the total routing cost for the leader through path $o(r)-i(r)-j(r)-d(r)$ as:
        \begin{equation}
            C^r_{kl}=\min\{\overline c^{rk}_{o(r)i} +\overline c^{rl}_{jd(r)}+ c_{ij}^r: i,j\in H\}. \label{Crkl1}
        \end{equation}

        The above expression is also valid for routes that do not have a first or a third leg. Specifically, the route of a given commodity $r\in R$ will have a first leg only if $o(r)\notin H$, and it will have a third leg  only if $d(r)\notin H$. If $o(r)\in H$, then $\overline c^{rk}_{o(r)i} +\overline c^{rl}_{jd(r)}+ c_{ij}^r=\overline c^{rl}_{jd(r)}+ c_{ij}^r$ so, for a given $l\in K$, $C^r_{kl}$ takes the same value for all $k\in K$. Likewise,  if $d(r)\in H$, then $\overline c^{rk}_{o(r)i} +\overline c^{rl}_{jd(r)}+ c_{ij}^r=\overline c^{rl}_{io(r)}+ c_{ij}^r$ so, for a given $k\in K$, $C^r_{kl}$ takes the same value for all $l\in K$. The routing costs defined in \eqref{Crkl1} determine, in each case, the best overall routing costs for the commodities. From the definition of \eqref{Crkl1}, we observe that the set $R$ can be reduced \emph{a priori}. 
        
        \begin{rem}
            If $\min_{k,\,l\in K}C^r_{kl} \geq b^r$ then commodity $r$ will not be profitable and can be removed from $R$.
        \end{rem}
        
        We can provide now an alternative problem definition obtained after preprocessing routing costs.

        \begin{prop}[\SDOD\ with preprocessed routing costs]\label{prop:defSD}
            Let $C^r_{kl}$ be the costs for routing commodity $r \in R$ if its origin is served by carrier $k$ and its destination is served by carrier $l$, defined as in \eqref{Crkl1}. Then, an optimal \SDOD\ solution is obtained by allocating each non-hub node to at most one carrier  and finding a subset of commodities $R^* \subseteq R$,  that maximize \[\sum_{r \in R^*} \left[ b^r - C^r_{k(r)l(r)} \right]\] where  for each $r \in R$, its origin and destination are allocated to carriers $k(r)$  and $l(r)$, respectively.
        \end{prop}
        In the above definition the allocation of non-hubs to carriers becomes the major decision, as the hub network is totally preprocessed. This alternative point of view will allow us to prove the problem complexity below.

        We have a similar result for the \SDOD$_F$ when outsourcing fees, $\overline p$ and $\overline q$, are fixed. Now, for a given commodity $r\in R$, if the carriers $k(r)$ and $l(r)$ are known, then the best interhub arc for routing $r$, $(i(r),j(r))$, can be determined \textit{a priori} by finding the best interhub arc $(i, j)\in A^r$ among where hubs $i$ and $j$ are optimal responses of $k(r)$ and $l(r)$ for outsourcing fees $\overline p$ and $\overline q$, respectively:
        \begin{equation*}
            (i(r),j(r))\in\arg\min\{\overline p^{r}_{i} +\overline q^{r}_{j}+ c_{ij}^r: i\in I_{k(r)}(r, \overline p),j\in J_{l(r)}(r, \overline q)\}.
        \end{equation*}
        That is interhub arc $(i, j)\in A^r$ can be used for routing $r$ only when $i\in I_{k(r)}(r, \overline p)$, and $j\in J_{l(r)}(r, \overline q)$. Thus, the preprocessed routing costs are defined as
        \begin{align}
            \widehat{C}^r_{kl}= \begin{cases}
                \min\{\overline p^r_i +\overline q^{r}_{j}+ c_{ij}^r: (i,j)\in A^r,\, \text{s.t. }i\in I_k(r, \overline p), j\in J_k(r, \overline q)\}\quad & \text{if $I_k(r, \overline p)\ne \emptyset$ and $J_l(r, \overline q)\ne \emptyset$}\\ + \infty & \text{otherwise.}
            \end{cases}\label{Crkl1-F}
        \end{align}
        \begin{prop}[\SDOD$_F$ with preprocessed routing costs]\label{prop:defSD-F}
            Let $\overline p$ and $\overline q$ be given outsourcing fees. Let also $\widehat C^r_{kl}$ be the costs for routing commodity $r \in R$ if its origin is served by carrier $k$ and its destination is served by carrier $l$, defined as in \eqref{Crkl1-F}. Then, an optimal \SDOD$_F$ solution is obtained by finding an allocation of each non-hub node to at most one carrier,  together with a subset of commodities $R^* \subseteq R$, that maximize \[\sum_{r \in R^*} \left[ b^r - \widehat{C}^r_{k(r)l(r)} \right]\] where for each $r \in R$, its origin and destination are allocated to carriers $k(r)$  and $l(r)$, respectively.       
        \end{prop}
        \subsubsection{Analysis for $r$\SDOD}
            {One of the consequences of Proposition \ref{propo2} is that for any given solution $\overline f$, $\overline t$, and $\overline a$, satisfying \eqref{const:allocate-fixed}-\eqref{const:dr_allocation_MINLP}, the outsourcing fees $\overline p$, $\overline q$ as determined by \eqref{best_p} and \eqref{best_q} are optimal for both \SDOD\ and $r$\SDOD. That is, \SDOD\ and $r$\SDOD\ are, in fact, the same problem. The same conclusion can be reached for $r$\SDOD\, as the preprocessed routing costs \eqref{Crkl1} are exactly the same.
            
            On the contrary, \SDOD$_F$ and $r$\SDOD$_F$ are different problems, as can be seen by comparing the preprocessed routing costs of \SDOD$_F$ and $r$\SDOD$_F$. For the latter, the preprocessed routing costs are defined as
            \begin{align} \label{eq:tildeC}
            \widetilde{C}^r_{kl}=
                \begin{cases}
                    \min\{\overline p^r_i +\overline q^{r}_{j}+ c_{ij}^r: (i,j)\in A^r,\, \text{s.t. }i\in rI_k(r, \overline p), j\in rJ_k(r, \overline q)\}\quad & \text{if $rI_k(r, \overline p)\ne \emptyset$ and $rJ_l(r, \overline q)\ne \emptyset$}\\ + \infty & \text{otherwise,}
                \end{cases}
            \end{align}
            so an optimal $r$-\SDOD$_F$ solution can be obtained by by maximizing  $\sum_{r \in R^*} \left[ b^r - \widetilde{C}^r_{k(r)l(r)} \right]$.
    \subsection{Problem complexity}
        We first prove NP-hardness of \SDOD. We then show that the problem remains NP-hard, even when the outsourcing fees are fixed.
        \begin{prop} 
            The \SDOD\ is NP-hard, even with a single hub node.
        \end{prop}
        \textbf{Proof}: We show this result by a polynomial reduction from the NP-hard Quadratic Semi-Assignment Problem (QSAP) \citep[see, e.g.,][]{Greenberg1969, Loiola2007}. The QSAP is defined as follows. We are given a set $F$ of facilities and a set $L$ of locations. For each pair of facilities $(i,\,j) \in F\times F$, $i \neq j$, we are given a weight $w_{ij} >0$ and for each pair of locations  $(k,\,l) \in L \times L$ we are given a distance $d_{kl} \ge 0$. The goal is to find a mapping $f: F \mapsto L$ such that each facility $i\in F$ is allocated to one location $f(i) \in L$, so that \[ \sum_{i, j \in F}  w_{ij} d_{f(i) f(j)} \]
        is minimized.
        
        Given an input instance of QSAP, we transform it into an \SDOD\ instance (according to Proposition \ref{prop:defSD}) as follows. Let $K:=L$, $H:=\{0\}$, $V\setminus H:=F$,  $R := \{ (o,d) : o,d \in F, o \neq d \}$,  $C^{r}_{kl} := w_{r} d_{kl}$, $b^r := M$, for all $r \in R$ where $M:= \max_{r \in R, k,l \in L} C^r_{kl} + 1$. In this transformation, the set of hub nodes is a singleton, and contains an auxiliary hub node denoted by $0$.
        
        Then, according to Proposition \ref{prop:defSD}, and by definition of $b^r$, in an optimal solution of such constructed \SDOD\ instance, all the commodities will be routed. The profit of such obtained solution is equal to
        $$ \sum_{r \in R} \left [ M - w_r d_{k(r)l(r)} \right] = M \cdot |F| \cdot (|F| - 1) - \sum_{r \in R}w_r d_{k(r)l(r)} =  M \cdot |F| \cdot (|F| - 1) - OPT,$$ where $OPT$ refers to the optimal solution value of QSAP. Indeed, the optimal allocation of non-hubs to carriers corresponds to the optimal allocation of facilities $F$ to locations $L$, and hence, we have a polynomial time reduction from QSAP to \SDOD. Since in this reduction the size and structure of the hub network play no role, we conclude that the result holds even when the hub network contains a single node. \hfill$\blacksquare$
        
        \begin{prop}
            The \SDOD\ can be solved in polynomial time if the number of carriers is constant.   
        \end{prop}
        \textbf{Proof}: For $|K|=\kappa=const$, there are $|V\setminus H|^{\kappa+1}$ possible allocation decisions that could be enumerated.  For each such allocation, calculating $C^r_{kl}$ and finding the optimal subset of commodities to be routed can be done in $O(|H|^2 \cdot |R|)$ time, and hence,  the optimal solution could be found in $O(\kappa^2 \cdot |H|^2 \cdot |R|)$ time. \hfill$\blacksquare$
    
        \noindent We now turn our attention to the complexity of \SDOD$_F$ and {$r$\SDOD$_F$}.
        \begin{prop} 
            The \SDOD$_F$ and the {$r$\SDOD$_F$} are NP-hard.
        \end{prop}
        \textbf{Proof}: Again we show this result by reduction from the QSAP. Given an input instance of QSAP, we will transform it into an \SDOD$_F$ instance as follows. Let $K:=L$ and $R := \{ (o,d) : o,d \in F, o \neq d \}$.

        Consider a backbone network consisting of $|K|^2$ hubs (as many as possible combinations for the carriers allocation to the origins and destinations of the commodities). Let $H:=\{v_{kl}\}_{kl\in K\times K}$ denote the set of hubs and let $V=H\cup F$. The leader's routing costs for the non-loop arcs of the backbone network are set to an arbitrarily large constant $D$, i.e., $c_{ij}^r=D$, for all $r\in R$, $i,j\in H$, $i\ne j$. With these interhub costs, the only routing paths for commodity $r$ that do not have arbitrarily large costs are of the form $o(r)-v_{kl}-d(r)$.
        
        The outsourcing fees are defined as follows: 
        \begin{equation*}
            \overline p^{r}_{i}=
            \begin{cases}
                \frac{1}{2\,}w_{o(r)d(r)}\,d_{kl}& \text{if } i=v_{kl}\\ 0 & \text{otherwise,}
            \end{cases}
            \quad \text{ and }\quad \overline q^{r}_{i}=
            \begin{cases}
                \frac{1}{2\,}w_{o(r)d(r)}\,d_{kl}& \text{if } i=v_{kl}\\ 0 & \text{otherwise.}
            \end{cases}
        \end{equation*}
        The carriers reservation prices are defined as follows:
        \begin{equation*}
            \overline c^{rk}_{o(r)i}=
            \begin{cases}
                \overline p^{r}_{i}& \text{if }  i=v_{kl} \text{ for some } l\in K\\ \overline p^{r}_{i}+1 & \text{otherwise,}
            \end{cases} \quad \text{ and }\quad \overline c^{rl}_{id(r)}=
            \begin{cases}
                \overline q^{r}_{i}& \text{if }  i=v_{kl} \text{ for some } k\in K\\ \overline q^{r}_{i}+1 & \text{otherwise.}
            \end{cases}
        \end{equation*}
        With the above reservation prices, for a given commodity $r\in R$ and a given pair of carriers ($k$ and $l$, respectively) $\overline c^{rk}_{o(r)i}=\overline p^{r}_{i}$ and $\overline c^{rl}_{id(r)}=\overline q^{r}_{i}$ for $i=v_{kl}$, but when $i\ne v_{kl}$, either $\overline c^{rk}_{o(r)i}>\overline p^{r}_{i}$ or $\overline c^{rl}_{id(r)}>\overline q^{r}_{i}$ (or both). That is, with the above routing costs, the only offer that would be accepted by the pair of carriers $k$ and $l$, is the one with $i=v_{kl}$. Moreover, the overall cost of such offer for the leader is precisely  $\overline p^{r}_{i}+\overline q^{r}_{i}=w_{r} d_{kl}$. Hence, following \eqref{Crkl1-F}, we define,  $\widehat C^{r}_{kl} := w_{r} d_{kl}$ for all $r \in R$, $k,l\in K$. In addition, $b^r := M$, for all $r \in R$ where $M:= \max_{r \in R, k,l \in L} \widehat C^r_{kl} + 1$. Then, by Proposition \ref{prop:defSD-F} and the definition of $b^r$,  in an optimal solution of such constructed \SDOD$_F$ instance, all the commodities will be routed. The profit of such solution is: 
        \begin{equation*}
            \sum_{r \in R} \left [ M - w_r d_{k(r)l(r)} \right] = M \cdot |F| \cdot (|F| - 1) - \sum_{r \in R}w_r d_{k(r)l(r)} =  M \cdot |F| \cdot (|F| - 1) - OPT,
        \end{equation*}
        where $OPT$ refers to the optimal solution value of QSAP. Indeed, the optimal allocation of non-hubs to carriers corresponds to the optimal allocation of facilities $F$ to locations $L$, and hence, we have a polynomial time reduction from QSAP to \SDOD$_F$.
        
        {Finally, we observe that the above construction can be also used to transform the QSAP instance into an $r$\SDOD$_F$ instance. Indeed, by defining an $r$\SDOD$_F$ instance on the same graph, and with the above outsourcing fees and reservation prices as above, the only offers that would be accepted by both carriers (non-negative net profit for both involved carriers) are those with both $\overline c^{rk}_{o(r)i}>\overline p^{r}_{i}$ and $\overline c^{rl}_{id(r)}>\overline q^{r}_{i}$. Hence, the preprocessed routing costs of the $r$\SDOD$_F$ instance would be precisely $\widetilde C^{r}_{kl} =\widehat C^{r}_{kl}$. Thus, an optimal solution of the $r$\SDOD$_F$ instance would produce an optimal solution to QSAP.}\hfill$\blacksquare$
    
\section{Single-level MILP reformulations for \SDOD} \label{sec:single-levelformu}
    Model \eqref{subeq:MINLPreformulation} given in Section \ref{sec:bilevelModel} is a single-level reformulation for \SDOD\ as an MINLP. Besides the bilinear terms in the objective function, this model also introduces additional bilinear terms in constraints enforcing lower-level optimality conditions \eqref{primal-dual-first}-\eqref{primal-dual-third}. As such, this model is computationally highly intractable. Instead of using standard linearization techniques to get rid of these bilinear terms, in this section we will exploit solution properties presented in Section \ref{sec:properties} to obtain more efficient MILP formulations for \SDOD.
    
    The first group of reformulations exploits Proposition \ref{propo2} and postpones the decisions on the optimal outsourcing fees, which, as we have seen, can be determined \emph{ex post}. In the first among these formulations (that we call the \emph{explicit paths formulation}), the follower's responses are modeled explicitly and hence a large number of four-index variables is needed. In the remainder of Section \ref{sec:EP}, we then show alternative flow-based formulations, that allow us to reduce the number of four-index variables. In Section \ref{sec:IP} we then provide another model (\emph{implicit paths formulation}) which exploits Proposition \ref{prop:defSD} and is derived from the problem definition obtained after preprocessing the routing costs. 
    
    As to what concerns \SDOD$_F$, we will explain how these models need to be adapted for this case. 
    
    \subsection{Explicit Paths (EP) Formulation}\label{sec:EP}
        The explicit Paths (EP) formulation considers the \SDOD\ as the problem of finding an optimal set of commodities to be served and an allocation of carriers to non-hubs that maximizes the overall net profit of the routed commodities (total revenue minus overall {routing cost}).
        Hence, a formulation for this model that uses the same decision variables as \eqref{subeq:MINLPreformulation}, except for the $p$ and $q$ variables is: 
        \begin{subequations}
            \begin{align}
                (EP)\qquad  \max &  \sum_{r\in R} b^r s^r-\sum_{r\in R}\sum_{(i,\,j)\in A^r}c_{ij}^rx_{ij}^r - && \sum_{r\in R}\sum_{k\in K}\sum_{i\in H}\left[\overline c^{rk}_{o(r)i} f^{rk}_i+\overline c^{rk}_{id(r)} t^{rk}_i\right] \label{of_fixed-EP}\\
                 \mbox{s.t. }
                & \sum_{k\in K}a_{ik} \leq 1 && \qquad\qquad\qquad i\in V\setminus H\label{const:allocate-fixed-EP}\\
                &  \sum_{i\in H}f_{i}^{rk} \leq a_{o(r)k} && \qquad\qquad\qquad r\in R, k\in K \label{const:or_allocation_explicitpaths}\\
                &  \sum_{j\in H}t_{j}^{rk} \leq a_{d(r)k} && \qquad\qquad\qquad r\in R, k\in K \label{const:dr_allocation_explicitpaths}\\
                &  \sum_{k\in K}f_{i}^{rk} {\ge}  \sum_{(i,\,j)\in A^r}x^r_{ij}&& \qquad\qquad\qquad r\in R, i\in H \text{ s.t. } o(r)\ne i\label{access-route-fixed}\\
                &  \sum_{k\in K}t_{i}^{rk} {\ge}  \sum_{(j,\,i)\in A^r}x^r_{ji}&& \qquad\qquad\qquad r\in R, i\in H \text{ s.t. } d(r)\ne i\label{distribution-route-fixed}\\
                &\sum_{(i,\,j)\in A^r} x^r_{ij}= s^r && \qquad\qquad\qquad r\in R \label{one-interhub-fixed-EP} \\
                & s^r\in\{0, 1\} &&\qquad\qquad\qquad r\in R, \\
                & f_{i}^{rk}, t_{i}^{rk}  \ge 0 && \qquad\qquad\qquad i\in H, r\in R\\
                & a_{ik} \in\{0, 1\}&& \qquad\qquad\qquad i\in V\setminus H, k\in K\\
                & x_{ij}^r\geq 0&& \qquad\qquad\qquad (i,\,j)\in A^r, r\in R.
            \end{align}
        \end{subequations}
        Recall that Proposition \ref{propo2} states that the optimal outsourcing fees for commodity $r$ should be set to $\overline{c}^{rk}_{o(r)i}$ if it is to be routed through hub $i$ and if the non-hub $o(r)$ is allocated to $k$, whereas all other values of $p^r$ should be set to zero. This allows to state the objective function as in \eqref{of_fixed-EP} and to get rid of $p$ and $q$ variables from  model \eqref{subeq:MINLPreformulation}. The meaning of Constraints \eqref{const:allocate-fixed-EP}-\eqref{one-interhub-fixed-EP} is similar to that of Constraints \eqref{const:allocate-fixed}--\eqref{one-interhub-fixed}, so the main difference of (\emph{EP}) with respect to the \textcolor{black}{bilevel} formulation of \SDOD\ affects the \textcolor{black}{removal of the lower level optimality conditions and corresponding reformulation of the} objective function.
        
        \paragraph{Optimality condition.} There is an optimal solution to $(EP)$ that satisfies the following inequalities:
        \begin{itemize}
            \item[]2 $x_{ij}^r+a_{(o(r)k)}+a_{(d(r)l)}\leq 2$ for all $r\in R$, $(i,j)\in A^r$, $k,l\in K$ such that $\overline c^{rk}_{o(r)i}+\overline c^{rl}_{d(r)j}+c^r_{ij}\geq b^r$.
        \end{itemize}
        We notice that the explicit path MILP formulation for \SDOD$_F$ is given by model \eqref{subeq:MINLPreformulation}, where $p^r_i$ and $q^r_i$ are fixed to their given values.
    
    \subsection{Reducing the number of 4-index variables}
        Since $f$ and $t$ variables used in (EP) involve four indices, next we develop alternative formulations with fewer first and third leg decision variables.
        \subsubsection{Explicit Flow formulation}
            We replace the original $f$ and $t$ variables while still modelling the carriers' decisions by considering the following continuous cost variables for the first and third travel legs:
            \begin{itemize}
                \item $\mathcal{F}^r_i \geq 0,\, r\in R,\,i\in H$, the cost of routing the first leg of commodity $r$ through hub $i$.
                \item $\mathcal{T}^r_i \geq 0,\, r\in R,\,i\in H$, the cost of routing the third leg of commodity $r$ through hub $i$.
            \end{itemize}
            The resulting Explicit Flow (\textit{EF}) formulation is:
            \begin{subequations} \label{eq:subeq:EF}
                \begin{align}
                    (EF)\qquad  \max &  \sum_{r\in R} b^r s^r-\sum_{r\in R}\sum_{(i,\,j)\in A^r}c_{ij}^rx_{ij}^r -\sum_{r\in R}\sum_{i\in H}(\mathcal{F}^r_i + \mathcal{T}^r_i) \label{ahlpf_of}\\
                     \mbox{s.t. }
                    & \sum_{k\in K}a_{ik} \leq 1 &&   i\in V\setminus H\label{const:allocate-ahlpf}\\
                    & \sum_{(i,\,j)\in A^r}x^r_{ij} \leq \sum_{k\in K}a_{o(r)k} && i\,\in H,\, r \in R: o(r)\notin H\label{const:or_x_ahlpf}\\
                    & \sum_{(i,\,j)\in A^r}x^r_{ij} \leq \sum_{k\in K}a_{d(r)k} && j\,\in H,\, r \in R: d(r)\notin H\label{const:dr_x_ahlpf}\\
                    &\mathcal{F}^r_i \geq \sum_{k\in K}\bar{c}^{rk}_{o(r)i}a_{o(r)k} - M^F_{ir}(1-\sum_{(i, j)\in A^r}x_{ij}^r) && i \in H,\, r\in R: o(r)\notin H\label{const:new_F_ahlpf}\\
                    &\mathcal{T}^r_i \geq \sum_{k\in K}\bar{c}^{rk}_{id(r)}a_{d(r)k} - M^T_{ir}(1-\sum_{(j, i)\in A^r}x_{ji}^r) && i \in H,\, r\in R: d(r)\notin H\label{const:new_T_ahlpf}\\
                    &\sum_{(i,\,j)\in A^r} x^r_{ij}= s^r && r\in R \label{const:one-interhub-ahlpf}\\
                    & s^r\in\{0, 1\} &&r\in R\\
                    &x_{ij}^r\in \{0, 1\}&& r\in R,\, (i,\,j)\in A^r\\
                    & a_{ik} \in\{0, 1\}&& i\in V\setminus H,\, k\in K\\
                    & \mathcal{F}^r_i, \mathcal{T}^r_i \geq 0&& i\in H,\, r\in R. \label{EF-last}
                \end{align}
            \end{subequations}
            The objective function has been updated to include the total overall profit and the routing costs of the first and third legs, which are now adjusted by subtracting the new cost variables $\mathcal{F}$ and $\mathcal{T}$.
    
            Similarly to previous Constraints \eqref{const:or_allocation_explicitpaths}--\eqref{const:dr_allocation_explicitpaths}, Constraints \eqref{const:or_x_ahlpf} and \eqref{const:dr_x_ahlpf} indicate that, to route a commodity, there must be a carrier allocated to its origin and/or destination in case both/any of them are non-hub nodes. Constraints \eqref{const:new_F_ahlpf} restrict the minimum value of the first leg fee. Given that nodes can only be allocated at most to one carrier, the right hand side coefficient considers that, if the origin of a commodity $r$ is not a hub and it is served through interhub arc $(i, j)$, then the value of the   fee for the the first leg must be the carrier's reservation price for connecting $o(r)$ to $i$. Similarly, Constraints \eqref{const:new_T_ahlpf} provide the lower bound value  for the  cost of the third leg.
    
            The big-$M$ coefficients used are $M^F_{ir} = \max_{k\in K}\left\{\bar{c}^{rk}_{o(r)i}\right\}$ for each node $i \in H$ and commodity $r \in R$ in the context of constraint \eqref{const:new_F_ahlpf}, and $M^T_{ir} = \max_{k\in K}\left\{\bar{c}^{rk}_{id(r)}\right\}$ for each node $i \in H$ and commodity $r \in R$ in the context of Constraints \eqref{const:new_T_ahlpf}.
            
        \subsubsection{Implicit-Flows formulation}
            We introduce an aggregated variant of the previous formulation that redefines the cost variables $\mathcal{F}^r_i$ and $\mathcal{T}^r_i$ for all potential origins/destinations of the first and third legs of commodity $r\in R$. These newly defined variables are denoted as $\overline{\mathcal{F}}^r$ and $\overline{\mathcal{T}}^r$. The model (\emph{IF}) is obtained from (\emph{EF})  by replacing Constraints  \eqref{const:new_F_ahlpf} and \eqref{const:new_T_ahlpf} with the following ones:
            \begin{subequations}
                \begin{align}
                    &\overline{\mathcal{F}}^r \geq \sum_{i\in H}\bar{c}^{rk}_{o(r)i}\sum_{(i, j)\in A^r}x_{ij}^r- M^F_{rk}(1-a_{o(r)k}) &&  k \in K,\, r\in R:\, o(r)\notin H\label{const:new2_F_ahlpf}\\
                    &\overline{\mathcal{T}}^r \geq \sum_{i\in H}\bar{c}^{rk}_{id(r)}\sum_{(j, i)\in A^r}x_{ji}^r- M^T_{rk}(1-a_{d(r)k}) && k \in K,\, r\in R:\, d(r)\notin H.\label{const:new2_T_ahlpf}
                \end{align}
            \end{subequations}
            The new big-$M$ coefficients used are $M^F_{rk} = \max_{i\in H}\left\{\bar{c}^{rk}_{o(r)i}\right\}$ for each commodity $r \in R$ and carrier $k\in K$ for Constraints \eqref{const:new2_F_ahlpf}, and $M^T_{rk} = \max_{i\in H}\left\{\bar{c}^{rk}_{id(r)}\right\}$ for each commodity $r \in R$ and carrier $k\in K$ for Constraints \eqref{const:new2_T_ahlpf}.
            The corresponding term in the objective function is replaced with:
            $-\sum_{r\in R}(\overline{\mathcal{F}}^r + \overline{\mathcal{T}}^r)$.

        \subsubsection{Removing big-$M$ coefficients}
            To accurately model first and third leg outsourcing fees in (\emph{EF}) and (\emph{IF}),  we further exploit conditions  that must hold when routing a commodity $r\in R$ whose origin and/or destination is a non-hub node. If $\sum_{i\in H\setminus{o(r)}} x^r_{ij}=1$ and $a_{o(r)k}=1$, then the first leg outsourcing fee for $r$ must be $\overline{\mathcal{F}}^r \ge \bar{c}^{rk}_{o(r)i}$. Similarly, for the third leg, if $\sum_{j\in H\setminus{d(r)}} x^r_{ji}=1$ and $a_{d(r)k}=1$, then $\overline{\mathcal{T}}^r \ge \bar{c}^{rk}_{id(r)}$.
            
            To represent these conditions for (\emph{EF}), we can replace \eqref{const:new_F_ahlpf} and \eqref{const:new_T_ahlpf} with the following constraints:
            \begin{subequations}
                \begin{align}
                    &\overline{\mathcal{F}}^r_i \ge \bar{c}^{rk}_{o(r)i} ( \sum_{(i,\,j)\in A^r}x^r_{ij} + a_{o(r)k} - 1) && i\in H,\, k \in K,\, r \in R: o(r)\notin H \label{2-index-F-i}\\
                    &\overline{\mathcal{T}}^r_i \ge \bar{c}^{rk}_{id(r)} ( \sum_{(j,\,i)\in A^r}x^r_{ji} + a_{d(r)k} - 1) && i\in H,\, k \in K,\, r \in R: d(r)\notin H. \label{2-index-T-i}
                \end{align}
            \end{subequations}
            For \emph{(IF)}, constraints \eqref{const:new2_F_ahlpf} and \eqref{const:new2_T_ahlpf} can be replaced with the following ones:
            \begin{subequations}
                \begin{align}
                    &\overline{\mathcal{F}}^r \ge \bar{c}^{rk}_{o(r)i} ( \sum_{(i,\,j)\in A^r}x^r_{ij} + a_{o(r)k} - 1) && i\in H,\, k \in K,\, r \in R: o(r)\notin H \label{const:new_F_ahlpf-IF}\\
                    &\overline{\mathcal{T}}^r \ge \bar{c}^{rk}_{id(r)} ( \sum_{(j,\,i)\in A^r}x^r_{ji} + a_{d(r)k} - 1) && i\in H,\, k \in K,\, r \in R: d(r)\notin H. \label{const:new_T_ahlpf-IF}
                \end{align}
            \end{subequations}

    \subsection{Adapting \emph{(}EP\emph{)}, \emph{(}EF\emph{)},  and \emph{(}IF\emph{)}  to \SDOD$_F$ and $r$\SDOD$_F$}
         Next we explain how to adapt formulations (\emph{EP}), (\emph{EF}),  and (\emph{IF})  for \SDOD$_F$ and $r$\SDOD$_F$.  The adaptations for \SDOD$_F$ will be referred to as (\emph{EP}$_F$), (\emph{EF}$_F$), and (\emph{IF}$_F$), respectively.  For $r$\SDOD$_F$ they will be referred to as ($r$\emph{EF}$_F$), ($r$\emph{EF}$_F$), and ($r$\emph{IF}$_F$), respectively.
         
        \paragraph{Formulations ({EP}$_F$) and ($r${EP}$_F$).} Formulation (\emph{EP})  can be easily adapted to \SDOD$_F$, by considering the objective function \eqref{of:MINLP} with the given outsourcing fees $\overline p^r_i$, $\overline q^r_i$, $i\in H, r\in R$, and adding the Constraints \eqref{subeq:primaldual} (note that the coefficients $\left[\max_{i\in H}\rho_{i}^{rk}\right]^+$ and $\left[\max_{i\in H}\gamma_{i}^{rk}\right]^+$ are now constants). When dealing with $r$\SDOD$_F$, the Constraints \eqref{subeq:primaldual} should be substituted by their relaxed counterpart \eqref{subeq:primaldual-relax}.
        
        \paragraph{Formulations ({EF}$_F$) and ($r${EF}$_F$).} In order to adapt (\emph{EF})  for \SDOD$_F$, we note that when the outsourcing fees are given, then $\mathcal{F}^r_i$ will take the value $p^r_i$, provided that the first leg of commodity $r$ is routed through hub $i$. Similarly, $\mathcal{T}^r_i$  will take the value  $q^r_i$, provided that the third leg of commodity $r$ is routed through hub $i$. However, the constraints
        \begin{subequations}
            \begin{align}
                &\overline{\mathcal{F}}^r_i \geq p^{r}_{i}\sum_{(i, j)\in A^r}x_{ij}^r &&   i\in H,\, r\in R:\, o(r)\notin H\label{const:new_F_ahlpf-fixed-rel}\\
                &\overline{\mathcal{T}}^r_i \geq q^{r}_{i}\sum_{(j, i)\in A^r}x_{ji}^r &&  i\in H,\, r\in R:\, d(r)\notin H,\label{const:new_T_ahlpf-fixed-rel}
            \end{align}
        \end{subequations}
        are not enough to guarantee that feasible solutions correspond to optimal carriers decisions. Indeed, when the outsourcing fees are given, optimal carriers responses require that the first leg of a commodity  served by carrier $k\in K$ is routed through some node $i\in I_{k}(r,p)$ and the third leg of a served commodity served by carrier $k\in K$ is routed through some node $i\in J_{k}(r,q)$. Note that such a choice of the \emph{connecting hub} vertices already guarantees that $\mathcal{F}^r_i\geq p^r_i\geq \overline c^{rk}_i$ and $\mathcal{T}^r_i\geq q^r_i\geq \overline c^{rk}_i$. Hence, for  (\emph{EF}$_F$) in formulation \eqref{eq:subeq:EF} we substitute Constraints \eqref{const:new_F_ahlpf}-\eqref{const:new_T_ahlpf} by
        \begin{subequations}
            \begin{align}
                & \overline{\mathcal{F}}^r_i \ge p^r_i ( \sum_{(i,\,j)\in A^r}x^r_{ij} + a_{o(r)k} - 1) && i\in I_{k}(r,p),\, k \in K,\, r\in R:\, o(r)\notin H\label{const:new_F_ahlpf-fixed}\\
                & \overline{\mathcal{T}}^r_i \ge q^r_i ( \sum_{(j,\,i)\in A^r}x^r_{ji} + a_{d(r)k} - 1) && i\in J_{k}(r,q),\, k \in K,\, r\in R:\, d(r)\notin H,\label{const:new_T_ahlpf-fixed}
             \end{align}
         \end{subequations} 
        \noindent whereas for ($r$\emph{EF}$_F$) it is enough to consider Constraints \eqref{const:new_F_ahlpf-fixed-rel}-\eqref{const:new_T_ahlpf-fixed-rel} instead.
        
        Note that \eqref{const:new_F_ahlpf-fixed}-\eqref{const:new_T_ahlpf-fixed} are, in fact a re-statement of Constraints \eqref{2-index-F-i}-\eqref{2-index-T-i} for the case of fixed fees which, in addition, guarantee the optimality of the carriers' responses.
        
        \paragraph{Formulations ({IF}$_F$) and ($r${IF}$_F$).} In order to adapt (\emph{IF}) for \SDOD$_F$, we proceed quite similarly to the case of (\emph{EF}). Now we observe that when the oursourcing fees are given, then $\overline{\mathcal{F}}^r$ will take the value $\sum_{i\in H}p^{r}_{i}\sum_{(i, j)\in A^r}x_{ij}^r$ and $\overline{\mathcal{T}}^r$ will take the value $\sum_{i\in H}q^{r}_{i}\sum_{(j, i)\in A^r}x_{ji}^r$, independently of the carriers to which $o(r)$ and $d(r)$ are assigned, respectively. Now, the constraints \vspace{-5.pt}
        \begin{subequations}
            \begin{align}
                &\overline{\mathcal{F}}^r \geq \sum_{i\in H}p^{r}_{i}\sum_{(i, j)\in A^r}x_{ij}^r &&   r\in R:\, o(r)\notin H\label{const:new2_F_ahlpf_fixed}\\
                &\overline{\mathcal{T}}^r \geq \sum_{i\in H}q^{r}_{i}\sum_{(j, i)\in A^r}x_{ji}^r &&  r\in R:\, d(r)\notin H\label{const:new2_T_ahlpf_fixed}
            \end{align}
        \end{subequations}
        \noindent are not enough to ensure that the obtained solutions correspond to optimal carriers decisions. Arguments similar to the above ones lead to use the following sets of constraints in (\emph{IF}$_F$):\vspace{-5.pt}
        \begin{subequations}
            \begin{align}
                & \overline{\mathcal{F}}^r \ge p^r_i ( \sum_{(i,\,j)\in A^r}x^r_{ij} + a_{o(r)k} - 1) && i\in I_{k}(r,p),\, k \in K,\, r\in R:\, o(r)\notin H\label{const:new_F_ahlpf-fixed-IF}\\
                & \overline{\mathcal{T}}^r \ge q^r_i ( \sum_{(j,\,i)\in A^r}x^r_{ji} + a_{d(r)k} - 1)&& i\in I_{k}(r,q),\, k \in K,\, r\in R:\, d(r)\notin H.\vspace{-10.pt}\label{const:new_T_ahlpf-fixed-IF}
             \end{align}
         \end{subequations}      
        Similarly, for ($r$\emph{IF}$_F$) it is enough to consider constraints \eqref{const:new2_F_ahlpf_fixed}-\eqref{const:new2_T_ahlpf_fixed} instead.
        
        The reader may again appreciate that \eqref{const:new_F_ahlpf-fixed-IF}-\eqref{const:new_T_ahlpf-fixed-IF} are, in fact a re-statement of constraints \eqref{const:new_F_ahlpf-IF}-\eqref{const:new_T_ahlpf-IF} for the case of fixed fees which, in addition, guarantee the optimality of the carriers' responses. 

    \subsection{Implicit Paths \emph{(IP)} formulation}\label{sec:IP}
        We now provide alternative single-level reformulations for \SDOD\  and \SDOD$_F$ based on Propositions \ref{prop:defSD} and  \ref{prop:defSD-F}, respectively. That is, we model the problems after preprocessing the routing costs. We use decision variables associated with potential implicit paths for serving the commodities, with preprocessed routing costs $C^r_{kl}$ and $\widehat C^r_{kl}$, as indicated in \eqref{Crkl1} and  \eqref{Crkl1-F}, respectively.  In addition to the service and allocation variables used in the bilevel formulation \eqref{mod:leader}, $s$, $a$, respectively, 
        we define binary decision variables $\pi^r_{kl}$, $r\in R$, $k,l\in K$, which take the value one if and only if first and third legs of commodity $r\in R$ are ``served'' by carriers $k,l\in K$, respectively. Then, the formulation for \SDOD\ is:
        \begin{subequations}
            \begin{align}
                (IP)\qquad  \max & \sum_{r\in R}b^rs^r-\sum_{r\in R}\sum_{k, l\in K} C^r_{kl} \pi_{kl}^r \label{of:RCAS3}\\
                 \mbox{s.t. }
                & \sum_{k\in K}a_{ik}\leq 1 &&   i\in V\setminus H \label{IP:allocate}\\
                & \sum_{\substack{k\in K\\l\in K}}\pi_{kl}^r =  s^r &&   r\in R\label{IP:serve}\\
                & \sum_{l\in K}\pi_{kl}^r\leq a_{o(r)k} && k\in K,\, r\in R:\, o(r)\notin H \label{IP:allocate-origin}\\
                & \sum_{k\in K}\pi_{kl}^r\leq a_{d(r)l}\, && l\in K,\, r\in R:\, d(r)\notin H \label{IP:allocate-destination}\\
                & \pi^r_{kl}\in\{0, 1\}&& k\in K,\, l\in K,\, r\in R\label{IP:pi}\\
                & s^r\in\{0, 1\} && r\in R\vspace{-10.pt}\\
                & a_{ik}\in\{0, 1\} && i\in V\setminus H,\, k\in K,\vspace{-10.pt}
            \end{align}
        \end{subequations}
        where constraints \eqref{IP:allocate} model the allocation decision to non-hub nodes by the leader. To serve a commodity, Constraints \eqref{IP:serve} indicate that routing costs must be incurred. For commodities whose origin are non-hub nodes, Constraints \eqref{IP:allocate-origin} imply that a carrier must be allocated to the origin node if the commodity is to be served. Analogously, Constraints \eqref{IP:allocate-destination} apply the same functionality for commodities with non-hub destination nodes. One could easily remove variables $s$ from the model, but we prefer to keep them.
        
        To obtain a valid formulation for \SDOD$_F$, the cost coefficients in the objective function must be set to $\widehat C^r_{kl}$ as defined in \eqref{Crkl1-F}. The resulting formulation will be referred to as (\emph{IP}$_F$). Formulation (\emph{IP}) can also be easily adapted to $r$\SDOD$_F$ by using cost coefficients $\widetilde C^r_{kl}$, as defined in \eqref{eq:tildeC}.
\section{Computational Results}\label{sec:compu}
    In this section we report results from computational experiments with benchmark instances adapted from the hub location literature. The main purpose of this computational study is to: $1)$ empirically evaluate the computational performance and scalability of the proposed models and formulations; and $2)$ derive some managerial insights concerning the sensitivity of decisions for \SDOD\ and \SDOD$_F$.
   
    Experiments have been performed on a PC equipped with a Ryzen 7 5700G CPU and 32Gb of RAM. Models have been implementated in Python 3.11 and solved with Gurobi 10.0.1. To provide reproducible results, we set Gurobi's \emph{Threads} parameter to 1 and turned off the \emph{Presolve} option.
    
    \subsection{Datasets, testing methodology, and parameter computing} \label{sec:data}
        For the computational experiments we have used benchmark instances generated from two well-known datasets from the hub location literature:
        \begin{itemize}
            \item The Civil Aeronautics Board (CAB) dataset, contains data from 100 cities in the United States of America \citep[see][]{Okelly1987}. It provides symmetric commodities demands, $w^{r}$ and unit arc costs $\hat c_{ij}$. Instances have been generated from this data set with a number of nodes $|V|\in \{20,\,30,\,40,\,50,\,60,\,70\}$. For every instance size  $n=|V|$ ten subsets of $n$ nodes have been randomly selected from the original data set, and their associated data adopted.  For each subset of nodes we have obtained two instances, one for every number of carriers $|K|\in\{3,\,4\}$. In total, we have generated 120 CAB instances. 
            \item The Australian Post (AP) dataset, first published by \citet{Ernst1996}, contains data from 200 nodes with non-integer asymmetric demands for the commodities $w^{r}$ and unit arc costs $\hat c_{ij}$. In these instances self flows $w_{ii}\ne 0$, $i\in V$. Still, we ignore them since we assume that  $o(r)\ne d(r)$ , $r\in R$. Instances from this data set with a number of nodes $n\in \{100,\, 120,\, 140,\, 160,\, 180,\, 200\}$ have been generated. For every instance size $n\ne 200$, ten subsets of $n$ nodes have been randomly selected from the original data set and their associated data adopted. One single set with all the nodes of the original data set is considered for $n=200$. For each subset of nodes we have obtained two instances, one for every number of carriers $|K|\in \{3,\,4,\,5\}$. In total, we have generated 153 AP instances.  
        \end{itemize}
        The parameters $n$, $|K|$, and $|H|$, are displayed in the first column of Tables  \ref{tab:cexp_01_CAB1050_pathbased&implicitpath}-\ref{tab:cexp_02_AP100200_bis}.
        
        Appendix \ref{app:data} gives the details on how we have generated the additional data that was not included in the original instances: $i)$ set of hubs, $ii)$ routing costs and reservation prices, $iii)$ commodities revenues, and $iv)$ outsourcing first and third leg fees for  \SDOD$_F$ and $r$\SDOD$_F$.  
        
        For all the experiments the time limit was set to 600 seconds for the smaller CAB instances and to 3600 seconds for the larger AP instances. 
    \subsection{Comparison of formulations for \SDOD}\label{sec:performance}
        Next we compare the computational performance of the formulations proposed  for \SDOD\,: (\emph{EP}), (\emph{EF}), (\emph{IF})  and (\emph{IP}) with the 120 instances generated from the CAB data set as explained above. 
        
        First, we carried out some preliminary tests for finding out the best modeling and algorithmic settings for (EF) and (IF).
        Detailed explanations on the considered alternatives and their performance can be found in Appendix \ref{Apdx:Tables}.  The best results were obtained with a \emph{lazy-callback} strategy according to which violated Big-$M$ constraints  \eqref{const:new_F_ahlpf}-\eqref{const:new_T_ahlpf} (or \eqref{const:new_F_ahlpf-fixed}-\eqref{const:new_T_ahlpf-fixed} for the fixed fees versions) are dynamically added. Thus, this is the variant used for both (EF) and (IF) in the experiments reported in the remainder of this section.

        \paragraph{Comparison of (EP), (EF), ({IF})  and ({IP}).} Results are summarized in Table \ref{tab:cexp_01_CAB1050_pathbased&implicitpath} where, for each tested formulation and group of instances, we show its runtime in seconds ($t(s)$), the amount of explored nodes in the enumeration tree (\emph{Nodes}) and the percent optimality gap at termination (\emph{GAP}($\%$)). The rows of the table display average results over the ten instances of the corresponding dimensions.
        
        Consistently, (\emph{IF}) is the formulation performing worse, followed by  (\emph{EF}).  (\emph{IF}) reaches the time limit already for instances with $n=30$ and $|K|=3$.   (\emph{EF})  produces optimal solutions for all instances although, on average, more than four minutes are needed when $n=70$. Therefore, we no longer consider (\emph{EF}) or (\emph{IF}), as they are significantly outperformed by (\emph{EP}) as well as by (\emph{IP}), both of which  solve all instances without branching and within a few seconds, even for the largest CAB instances with $n=70$ and $|K|=4$.
        
        The computing times of (\emph{EP}) and (\emph{IP}) are similar for small size instances of up to 30 nodes. Still, (\emph{IP}) clearly outperforms (\emph{EP}) as $n$ and $|K|$ increase. This is due to the smaller number of variables and constraints required by (\emph{IP}). For the largest instance referenced in Table \ref{tab:cexp_01_CAB1050_pathbased&implicitpath}, with $n=70$ and $|K|=4$, (\emph{EP}) has $440,214$ variables and $108,843$ constraints, whereas (\emph{IP}) has $72,261$ variables and $44,170$ constraints.

        \begin{table}[H]
            \caption{Computational performance comparison for CAB[20-70] instances with $|K| \in \{3,4\}$. \label{tab:cexp_01_CAB1050_pathbased&implicitpath}}
            \resizebox{\textwidth}{!}{\begin{tabular}{c|rrr|rrr|lrr|rrr}
            \multicolumn{1}{c|}{Instance} & \multicolumn{3}{c|}{(EP)} & \multicolumn{3}{c|}{(EF)} & \multicolumn{3}{c|}{(IF)} & \multicolumn{3}{c}{(IP)} \\ \hline
            \multicolumn{1}{c|}{$|V|.|K|.|H|$} & \multicolumn{1}{c}{$t(s)$} & \multicolumn{1}{c}{Nodes} & \multicolumn{1}{c|}{GAP(\%)} & \multicolumn{1}{c}{$t(s)$} & \multicolumn{1}{c}{Nodes} & \multicolumn{1}{c|}{GAP(\%)} & \multicolumn{1}{c}{$t(s)$} & \multicolumn{1}{c}{Nodes} & \multicolumn{1}{c|}{GAP(\%)} & \multicolumn{1}{c}{$t(s)$} & \multicolumn{1}{c}{Nodes} & \multicolumn{1}{c}{GAP(\%)} \\ \hline
            20.3.1 & 0.03 & 1 & 0.00 & 0.35 & 1 & 0.00 & \multicolumn{1}{r}{67.04} & 50485 & 0.00 & 0.01 & 1 & 0.00 \\
            20.4.1 & 0.03 & 1 & 0.00 & 0.51 & 1 & 0.00 & \multicolumn{1}{r}{244.50} & 137555 & 0.01 & 0.02 & 1 & 0.00 \\
            30.3.2 & 0.10 & 1 & 0.00 & 2.00 & 1 & 0.00 & $\quad$-- & 33854 & 0.25 & 0.03 & 1 & 0.00 \\
            30.4.2 & 0.11 & 1 & 0.00 & 2.02 & 1 & 0.00 & $\quad$-- & 33237 & 0.82 & 0.04 & 1 & 0.00 \\
            40.3.2 & 0.29 & 1 & 0.00 & 7.81 & 1 & 0.00 & $\quad$-- & 22599 & 0.67 & 0.06 & 1 & 0.00 \\
            40.4.2 & 0.32 & 1 & 0.00 & 8.21 & 1 & 0.00 & $\quad$-- & 18005 & 1.01 & 0.10 & 1 & 0.00 \\
            50.3.3 & 0.68 & 1 & 0.00 & 25.20 & 2 & 0.00 & $\quad$-- & 6157 & 0.76 & 0.10 & 1 & 0.00 \\
            50.4.3 & 0.74 & 1 & 0.00 & 30.78 & 7 & 0.00 & $\quad$-- & 3484 & 1.36 & 0.16 & 1 & 0.00 \\
            60.3.3 & 1.57 & 1 & 0.00 & 92.05 & 1 & 0.00 & $\quad$-- & 1796 & 1.51 & 0.16 & 1 & 0.00 \\
            60.4.3 & 1.83 & 1 & 0.00 & 87.15 & 5 & 0.00 & $\quad$-- & 1115 & 2.07 & 0.31 & 1 & 0.00 \\
            70.3.4 & 2.95 & 1 & 0.00 & 248.80 & 5 & 0.00 & $\quad$-- & 381 & 1.39 & 0.23 & 1 & 0.00 \\
            70.4.4 & 3.23 & 1 & 0.00 & 268.26 & 35 & 0.00 & $\quad$-- & 241 & 2.56 & 0.48 & 1 & 0.00 \\ \hline
            \end{tabular}}\\
            \footnotesize{Time limit reached ($-$).}
        \end{table}
        
        \paragraph{Scalability of (\emph{EP}) and (\emph{IP}).} To test the scalability of  (\emph{EP}) and (\emph{IP}), we use the 153 larger instances generated from the AP data set  with $n\in[100, 200]$ and $|K|\in\{3, 4, 5\}$. 
        Table \ref{tab:cexp_02_AP100200_bis}  
        reports computing times in seconds for each phase of the solution process: preprocessing  (\textit{Preprocess}), loading the model  (\emph{Loading}), and solving the respective MILP (\emph{Runtime}). All instances are solved to optimality at the root node, so the number of nodes explored in the enumeration tree is not shown. Rows corresponding to instances with $n<200$ display average values over the ten instances in the corresponding group, whereas the last three rows give results for the only instance with $n=200$, for the different values of $|K|$.

        Since the preprocessing phase of the two formulations is quite different, it is not surprising that the time consumed at this phase is quite different as well. The preprocessing of (\emph{EP}) removes all infeasible and/or unnecessary  \emph{routing} variables $x$, $f$ and $t$. Preprocessing (\emph{IP}) involves computing the total routing cost $C^r_{kl}$ for each commodity $r\in R$ and pair of carriers $k,\, l\in K$. For a fair comparison of formulations, we report the time needed for computing such costs, even if it depends on the data structures and not on the solver.
     
        The loading phase includes the time needed for generating and loading the model. Note that for both (\emph{EP}) and (\emph{IP}), loading times exceed solution times, which can be explained by the fact that both formulations use four-index variables and several families of constraints that involve three indices. Hence, the time needed to load the formulations increases very fast with the input size.
        
        The results shown in Table \ref{tab:cexp_02_AP100200_bis} confirm the overall superiority of (\emph{IP}), even if preprocessing times are larger and increase faster for (\emph{IP})  than for  (\emph{EP}).  Observe that, on average, the computing time needed to load and solve (\emph{IP}) is 13 and 20 times smaller than that of (\emph{EP}), respectively.

        \begin{table}[H]
            \caption{Computing times of the different phases for solving (\emph{EP}) and (\emph{IP}) for AP[100-200] with $|K| \in \{3,4,5\}$. \label{tab:cexp_02_AP100200_bis}}
            \resizebox{\textwidth}{!}{\begin{tabular}{lrrrrrr}
            \multicolumn{1}{c}{Instance} & \multicolumn{3}{c}{(EP)} & \multicolumn{3}{c}{(IP)} \\ \hline
            \multicolumn{1}{c}{$|V|.|K|.|H|$} & \multicolumn{1}{c}{Preprocess (s)} & \multicolumn{1}{c}{Loading (s)} & \multicolumn{1}{c}{Runtime (s)} & \multicolumn{1}{c}{Preprocess (s)} & \multicolumn{1}{c}{Loading (s)} & \multicolumn{1}{c}{Runtime (s)} \\ \hline
            100.3.5$\qquad$ & 0.57$\quad$ & 20.55$\quad$ & 9.32$\quad$ & 4.85$\quad$ & 2.62$\quad$ & 0.51$\quad$ \\
            100.4.5$\qquad$ & 0.59$\quad$ & 24.81$\quad$ & 9.80$\quad$ & 8.36$\quad$ & 4.27$\quad$ & 1.24$\quad$ \\
            100.5.5$\qquad$ & 0.61$\quad$ & 29.08$\quad$ & 10.83$\quad$ & 12.87$\quad$ & 6.15$\quad$ & 1.97$\quad$ \\
            120.3.6$\qquad$ & 1.24$\quad$ & 38.64$\quad$ & 21.44$\quad$ & 10.05$\quad$ & 3.68$\quad$ & 0.99$\quad$ \\
            120.4.6$\qquad$ & 1.29$\quad$ & 46.05$\quad$ & 22.36$\quad$ & 17.48$\quad$ & 5.95$\quad$ & 1.81$\quad$ \\
            120.5.6$\qquad$ & 1.30$\quad$ & 54.61$\quad$ & 24.22$\quad$ & 26.25$\quad$ & 8.89$\quad$ & 4.30$\quad$ \\
            140.3.7$\qquad$ & 2.32$\quad$ & 65.73$\quad$ & 41.26$\quad$ & 18.34$\quad$ & 5.31$\quad$ & 1.44$\quad$ \\
            140.4.7$\qquad$ & 2.34$\quad$ & 78.51$\quad$ & 43.45$\quad$ & 31.49$\quad$ & 8.62$\quad$ & 3.14$\quad$ \\
            140.5.7$\qquad$ & 2.29$\quad$ & 90.72$\quad$ & 45.68$\quad$ & 49.94$\quad$ & 12.25$\quad$ & 6.25$\quad$ \\
            160.3.8$\qquad$ & 3.83$\quad$ & 105.32$\quad$ & 70.59$\quad$ & 31.79$\quad$ & 6.81$\quad$ & 2.31$\quad$ \\
            160.4.8$\qquad$ & 3.87$\quad$ & 120.31$\quad$ & 72.37$\quad$ & 53.90$\quad$ & 10.99$\quad$ & 4.69$\quad$ \\
            160.5.8$\qquad$ & 3.93$\quad$ & 137.23$\quad$ & 77.17$\quad$ & 81.64$\quad$ & 15.81$\quad$ & 9.29$\quad$ \\
            180.3.9$\qquad$ & 6.25$\quad$ & 158.42$\quad$ & 142.13$\quad$ & 48.57$\quad$ & 8.78$\quad$ & 2.96$\quad$ \\
            180.4.9$\qquad$ & 6.48$\quad$ & 187.64$\quad$ & 135.34$\quad$ & 85.52$\quad$ & 13.98$\quad$ & 7.34$\quad$ \\
            180.5.9$\qquad$ & 6.38$\quad$ & 213.35$\quad$ & 136.69$\quad$ & 129.93$\quad$ & 20.79$\quad$ & 12.32$\quad$ \\
            200.3.10$\qquad$ & 9.43$\quad$ & 314.43$\quad$ & 191.25$\quad$ & 78.47$\quad$ & 10.95$\quad$ & 5.70$\quad$ \\
            200.4.10$\qquad$ & 12.63$\quad$ & 365.26$\quad$ & 202.15$\quad$ & 141.74$\quad$ & 17.86$\quad$ & 10.36$\quad$ \\
            200.5.10$\qquad$ & 9.56$\quad$ & 333.90$\quad$ & 243.67$\quad$ & 212.63$\quad$ & 26.30$\quad$ & 15.75$\quad$ \\ \hline
            \end{tabular}}
        \end{table}
        
    \subsection{Comparison of formulations for \SDOD$_F$ and $r$\SDOD$_F$}\label{sec:performance-fixed}
        First we compare the effect of the preprocessing on the formulations for the models with fixed outsourcing fees, \SDOD$_F$ and $r$\SDOD$_F$. For these experiments two different sets of fixed outsourcing fees for first and third legs have been generated for each instance, as explained in Appendix \ref{ap:fixed-out}. The first one considers as outsourcing fees the carriers maximum reservation prices, whereas the second one considers as outsourcing fees the average values over the reservation prices of all the carriers.
        
        Table \ref{tab:FFees_EP_VarRem} of Appendix \ref{ap:fixed-out} displays average percentages of the number of variables fixed in the preprocessing phase of (\emph{EP}$_F$) and (\emph{IP}$_F$). Since no noticeable differences have been observed when varying the number of carriers or when considering \emph{maximum} or \emph{average} outsourcing fees, for a given input size $n$, the average is computed over all the instances of size $n$ for varying values of $|K|$ and the two considered alternatives of outsourcing fees. Disaggregated results are also shown in Tables \ref{tab:Granular_max}-- \ref{tab:Granular_avg} of Appendix \ref{ap:fixed-out}.
        
        As can be seen, for (\emph{EP}$_F$) the preprocessing  removes, on average, more than $90\%$ of the routing variables $x$ and about $50\%$ of the first and third leg variables $f$ and $t$, respectively. Note that, even if the percentage of fixed variables is high, the remaining formulation still has a considerable size. For example, on the AP instance with 180 nodes and 4 carriers the remaining number of variables ($s,\,x,\,f$ and $t$) is  691,360.  
        The reduction in the number of $x$ variables is aligned with the results reported in \citet{Alibeyg2018}, where more than 85$\%$ routing variables are fixed in the preprocess of a prize-collecting hub location model.

        Tables \ref{tab:AvgFFees_CAB_COMP}--\ref{tab:MaxFFees_CAB_COMP} of Appendix \ref{ap:fixed-out} display results of the different formulations for \SDOD$_F$ and  $r$\SDOD$_F$ for the (preprocessed) benchmark instances from the CAB data set, using the \emph{maximum} and \emph{average} fixed outsourcing fees, respectively. We report runtimes in seconds, ($t(s)$), number of explored nodes (\textit{Nodes}), and percent optimality gaps at termination ($\%$GAP). Each row displays average values over the 10 instances with the corresponding number of nodes and carriers. Broadly speaking, the obtained results show again the superiority of implicit-path formulations over explicit-path formulations for both  \SDOD$_F$ and $r$\SDOD$_F$. 
        
        As can be observed,  the choice  of the fixed outsourcing fees has a different impact on the difficulty for solving the selected formulation.  
        In particular, for ($r$\emph{EP}$_F$) the increase of computing times 
        with the instance size remains very moderate for \emph{maximum} outsourcing fees but is remarkable for \emph{average} outsourcing fees. 
        This is clearly related to the algorithmic burden for the separation of  Constraints \eqref{const:new_F_ahlpf-fixed}-\eqref{const:new_T_ahlpf-fixed}, which become redundant for \emph{maximum} fees (any offer would be accepted by the carriers), but are dynamically separated for \emph{average} outsourcing fees. Nevertheless, the effect of the outsourcing fee, is not so relevant for (\emph{EP}$_F$), were notably fewer offers would be accepted by the carriers independently of the choice of outsourcing fees, so the separation of the inequalities has a minor impact on the overall performance of the formulation. 

        Instead, neither implicit path formulation is noticeably affected by the choice of outsourcing fees, even if a slight influence can be appreciated in the case of ($r$\emph{IP}$_F$).

        The joint effect of the above remarks is that, while all tested instances could be solved to proven optimality both with (\emph{EP}$_F$) and (\emph{IP}$_F$), the largest $70$ node CAB instances could not be solved within the five minutes time limit by the  worse performing formulation, ($r$\emph{EP}$_F$), when \emph{average} outsourcing fees are considered.

        The above observations are also supported by the results on the larger instances from the AP data set with $n\in \{100,\,120,\,140,\,160,\,180,\,200\}$ and $|K|\in\{3,\,4,\,5\}$, which are summarized in Tables \ref{tab:AvgFFees_AP_COMP} -\ref{tab:MAXFFees_AP_COMP} of Appendix \ref{ap:fixed-out}. For each tested formulation we include the computing times consumed in the preprocessing phase (\emph{Prep.}), instance loading phase (\emph{Load.}) and solution phase (\emph{Run.}). The rows corresponding to instances with less than 200 nodes display average values over the ten instances in the corresponding group, whereas the last three rows give the results for the only instance with 200 nodes and the corresponding number of carriers. Note that for the largest instances with 200 nodes, the solver runs out of memory in the optimization phase of both (\emph{EP}$_F$) and ($r$\emph{EP}$_F$) (already when solving the linear programming relaxation) both with \emph{maximum} and \emph{average} outsourcing fees, whereas with both (\emph{IP}$_F$) and ($r$\emph{IP}$_F$) all the instances can be solved to proven optimality within the one hour time limit for both versions of the  outsourcing fees.

        Further insight on the scalability of the implicit-path formulations, which are the best performing ones, can be obtained by increasing the number carriers on the instances generated from the largest AP dataset. Table \ref{tab:FFees_Stressing}  displays computing times for instances with $n=200$ nodes  and a number of carriers $|K|\in\{6,\,7,\,8\}$. As can be observed, for both outsourcing fees variants, the preprocessing times of ($r$\emph{IP}$_F$) are more than one order of magnitude higher than those of (\emph{IP}$_F$), and the differences increase with the value of $|K|$. This is clearly due to the increase of the computational burden required to compute the coefficients $\widehat C^r_{km}$, which, for the largest instance with eight carriers, takes more than one hour.
        
        \begin{table}[H]
            \caption{Time comparison for fixed outsourcing fees when considering (\emph{IP}$_F$) and $|K|\in\{6,\,7,\,8\}$ on AP dataset with 200 nodes.\label{tab:FFees_Stressing}}
            \resizebox{\textwidth}{!}{\begin{tabular}{lrrrrrr|rrrrrr}
            \multicolumn{1}{c}{} & \multicolumn{6}{c|}{(\emph{IP}$_F$)} & \multicolumn{6}{c}{($r$\emph{IP}$_F$)} \\
            \multicolumn{1}{c}{\multirow{2}{*}{$|V|.|K|.|H|$}} & \multicolumn{3}{c}{Average} & \multicolumn{3}{c|}{Maximum} & \multicolumn{3}{c}{Average} & \multicolumn{3}{c}{Maximum} \\
            \multicolumn{1}{c}{} & \multicolumn{1}{c}{Prep.} & \multicolumn{1}{c}{Load.} & \multicolumn{1}{c}{Run.} & \multicolumn{1}{c}{Prep.} & \multicolumn{1}{c}{Load.} & \multicolumn{1}{c|}{Run.} & \multicolumn{1}{c}{Prep.} & \multicolumn{1}{c}{Load.} & \multicolumn{1}{c}{Run.} & \multicolumn{1}{c}{Prep.} & \multicolumn{1}{c}{Load.} & \multicolumn{1}{c}{Run.} \\ \hline\hline
            200.6.10 & 186.13 & 35.69 & 3.14 & 306.85 & 35.82 & 2.72 & 1160.50 & 35.95 & 25.31 & 2491.58 & 35.73 & 9.29 \\
            200.7.10 & 273.30 & 47.66 & 4.43 & 409.89 & 47.85 & 3.81 & 1617.29 & 47.26 & 40.87 & 3389.77 & 48.02 & 13.52 \\
            200.8.10 & 525.58 & 60.01 & 7.49 & 526.06 & 60.41 & 5.20 & 2234.48 & 60.45 & 145.44 & 4424.94 & 60.69 & 17.72\\
            \hline
            \end{tabular}}
        \end{table}

    \subsection{Analysis of optimal solutions}\label{sec:sensitivity}
        In this section we analyze some characteristics of optimal solutions to \SDOD\,, \SDOD$_F$, and $r$\SDOD$_F$  obtained with the best performing formulations, (\emph{IP}), (\emph{IP}$_F$), and ($r$\emph{IP}$_F$), respectively. For this analysis we have used 24 instances generated from the CAB dataset, with a number of nodes ranging in 30-100, and a number of carriers $|K|\in\{3,\,4,\,5\}$. These instances have been generated so they share the backbone network with the same 5 hub nodes. Moreover, when increasing the number of nodes, we ensure that the smaller instance is contained in the larger one, so the $\widehat C^r_{km}$ values corresponding to pairs of nodes of the smaller network remain unchanged.
        
        Figure \ref{fig:CABCommEvol} depicts the number of served commodities for \SDOD\, and compares it with \SDOD$_F$ for both \emph{average} outsourcing fees (Figure \ref{fig:CABCommEvol}(a), left) and \emph{maximum} outsourcing fees (Figure \ref{fig:CABCommEvol}(b), right). Both figures also depict the total number of commodities (which increases quadratically in $n$).
        
        Indeed, the number of commodities served 
        is always notably larger for \SDOD\, than for \SDOD$_F$. Even if, for both models, the number of served commodities increases with $|V|$, the increase on the number of served commodities is smaller than the increase on the number of commodities, resulting in an overall decrease of the service rate.  Still, the service rate of \SDOD\, remains close to $\frac{1}{2}$ for all instance sizes, whereas it is considerably smaller for \SDOD$_F$ with both \emph{average} and and \emph{maximum} outsourcing fees.

        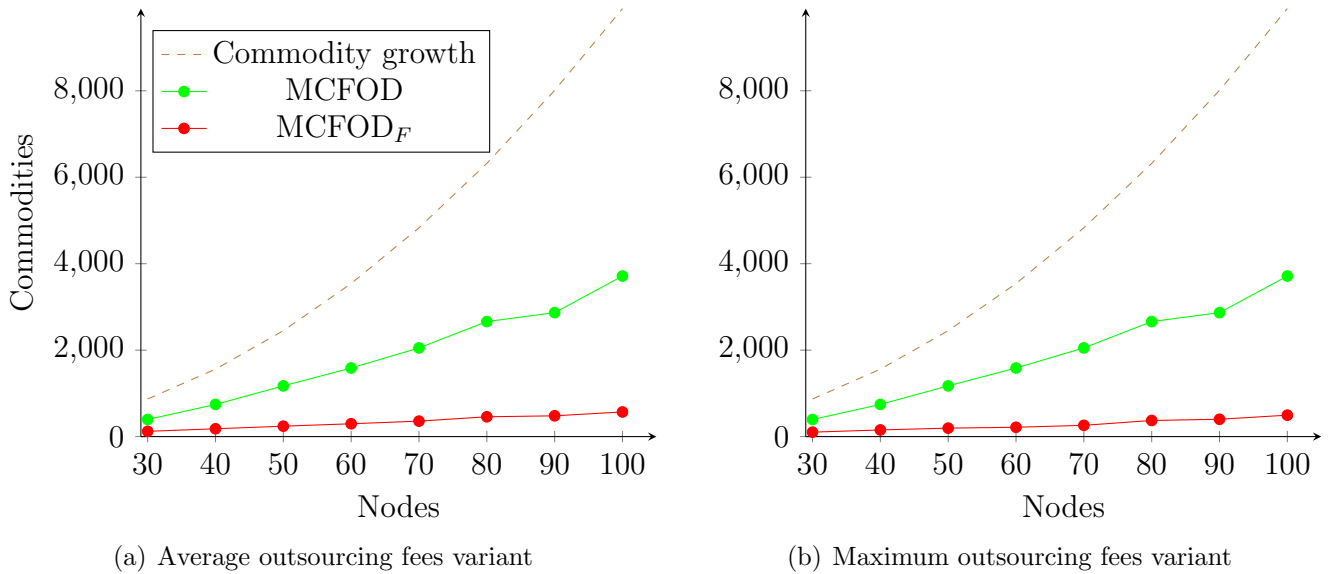
\begin{figure}[H]
            {\subfigure[Average outsourcing fees variant]{
                \begin{tikzpicture}
                    \begin{axis}
                        [
                        xlabel={Nodes}, 
                        ylabel={Commodities}, 
                        axis lines=left, 
                        xmin = 29,
                        xmax = 105,
                        ymin = 0,
                        xtick = data,
                        legend style={at={(0.35,0.95)}, anchor=north,legend columns=1},
                        ]
                        \addplot[brown, dashed] table [col sep=comma, x=Nodes, y=Total_Commodities]{01_ComEvol.csv};
                        \addplot[green, mark=*, solid] table [col sep=comma, x=Nodes, y=Qty_Opt]{01_ComEvol.csv};
                        \addplot[red, mark=*, solid] table [col sep=comma, x=Nodes, y=Qty_Served]{01_ComEvol.csv};
                        \legend{Commodity growth, \SDOD, \SDOD$_F$}
                    \end{axis}
                \end{tikzpicture}}
                \quad 
            \subfigure[Maximum outsourcing fees variant]{
                \begin{tikzpicture}
                        \begin{axis}
                            [
                            xlabel={Nodes}, 
                            ylabel={}, 
                            axis lines=left, 
                            xmin = 29,
                            xmax = 105,
                            ymin = 0,
                            xtick = data,
                            ]
                            \addplot[brown, dashed] table [col sep=comma, x=Nodes, y=Total_Commodities]{01_ComEvolMax.csv};
                            \addplot[green, mark=*, solid] table [col sep=comma, x=Nodes, y=Qty_Opt]{01_ComEvolMax.csv};
                            \addplot[red, mark=*, solid] table [col sep=comma, x=Nodes, y=Qty_Served]{01_ComEvolMax.csv};
                        \end{axis}
                    \end{tikzpicture}}
                }
            \caption{Comparison between optimal \SDOD\, and \SDOD$_F$ solutions on served commodities and commodity evolution for instances from the CAB dataset. \label{fig:CABCommEvol}}
        \end{figure}
        
        This decrease of the service rate can be better appreciated in Table \ref{tab:FFees_CAB_COMCOMP} of the Appendix \ref{app:service-rate}, which displays the average percentage of served commodities in optimal solutions to (\emph{IP}), (\emph{IP}$_F$), and ($r$\emph{IP}$_F$). The row corresponding to each  value of $n$  shows average results over the three instances with $n$ nodes and a number of carriers $|K|\in\{3,\,4,\,5\}$.      
        As expected, the percentage of served commodities, is higher for \SDOD\,, where the leader chooses the best  outsourcing fees, than for any of the models with fixed fees. Within the two models with fixed fees, the highest service rates are attained by the \emph{relaxed}  model where carriers  accept any profitable offer. 
        The service rates attained by $r$\SDOD$_F$ are only slightly smaller than those of \SDOD, but notably larger than those of \SDOD$_F$. Still, as $n$ increases the percentage of served commodities decreases for all three models.

        The choice of outsourcing fees affects differently \SDOD$_F$ and $r$\SDOD$_F$. For \SDOD$_F$, \emph{average}  outsourcing fees lead, on average, to a $9.7\%$ service rate, outperforming the $7.2\%$ average service rate obtained with  \emph{maximum} outsourcing fees. Very few offers are acceptable for the carriers with \SDOD$_F$ for both types of outsourcing fees,  but the cost for the leader of the accepted offers  is cheaper when using \emph{average} outsourcing fees. On the contrary, for $r$\SDOD$_F$,  the \emph{maximum} outsourcing fees produce an average service rate of around $40\%$, which is higher than the $38\%$  obtained with the \emph{average}  outsourcing fees. This is natural since with $r$\SDOD$_F$ all \emph{maximum} outsourcing fees offers are acceptable (profitable) for the carriers so the leader has a wider range of available options among which to decide. 

        We finally analyze the total profit produced by the different models. Figure \ref{fig:CABProfitAvg} in the Appendix \ref{app:service-rate} compares the net profit for optimal solutions to \SDOD, \SDOD$_F$, and $r$\SDOD$_F$. The figure shows that the highest profits are obtained with  \SDOD, even if $r$\SDOD$_F$ compares very closely. The superiority of \SDOD\, is due to its flexibility, as it allows the leader to adapt the values of the outsourcing fees according to its own interest. The difference between \SDOD\, and $r$\SDOD$_F$ is due, not only to the fact that in the former commodities may become profitable by decreasing the outsourcing fees, but also to the fact offers that in $r$\SDOD$_F$ are rejected by the carriers, may become acceptable for the carriers (and also for the leader) by increasing the outsourcing fees in \SDOD. This is also consistent with the results of Table \ref{tab:FFees_CAB_COMCOMP}, where it is shown that \SDOD\, always display the largest percentage of served commodities.

\section{Conclusions}\label{sec:conclu}
    We have addressed a problem of routing a subset of commodities within an existing network so as to maximize the collected profit. Contrary to the standard setting in which a single agent controls all the decisions, we introduced a new problem in which interactions between several actors are modeled {within a Stackelberg game}. The major decision maker (the leader) selects a subset of commodities to serve, and also seeks to outsource the first- and the last-leg of these commodities to external carriers. These {multiple independent} carriers receive the outsourcing offers of the leader and decide whether to accept them or not, based on their  individual profits. In this context, the leader makes critical decisions regarding carrier {allocation to non-hub nodes} and corresponding payments for outsourcing the service, while anticipating the optimal responses of the carriers. Initially, we formulated this problem as a bilevel MINLP model. Two additional problem settings were considered: one in which the leader only decides how to allocate carriers to non-hub nodes, whereas the outsourcing fees are fixed, and the other one in which the followers refuse the offer only if their resulting profit is negative. We showed that all considered problem variants are NP-hard and proposed several MILP formulations for solving them to optimality. In a computational study, we studied computational scalability of these formulations and compared them based on different criteria such as the number of nodes explored, runtimes and gaps at termination. Among the formulations considered, the implicit path formulation, which simultaneously encodes allocation and routing decisions, exhibits the best computational performance. It requires a more computationally demanding preprocessing phase, but it allows to significantly reduce the size of the model, by removing non-profitable carrier-to-non-hub allocations.

     Furthermore, we conducted a comprehensive analysis of the effects of various parameters on problem resolution and on the structure of the obtained solution, notably the profit of the leader and the service rate. The high reservation prices of the carriers tend to have a larger influence on runtimes when the number of available carriers and the number of nodes is larger, justifying a more complex decision for allocating the carriers and deciding the commodities to serve. Indeed, when these values increase, the overall service rate tends to decrease. The overall profit of the leader tends to decrease as the carriers apply a larger reservation price, meanwhile the carriers profit tends to increase. More importantly, our computational experiments indicate that including decisions about outsourcing fees can significantly improve service rates and the profit of the major firm, compared to the setting in which the outsourcing fees are regulated. Finally, we also investigated how the profit-maximizing vs. reservation-price-oriented strategies of the carriers affect the overall profit of the leader and the service rate. We observed that  profit-maximizing carriers tend to reduce the overall service rate and profit of the leader as compared to reservation-price-oriented ones. This could be expected, since profit-maximizing carriers will accept fewer transportation offers, thus reducing the profitable alternatives for the leader.

    As future work, we identify the following several relevant extensions of the \SDOD\,:
    \begin{itemize}
        \item Variations concerning the backbone network or the service mode of commodities, which may fit better some potential applications. These include single-allocation of origins/destinations, where all commodities with the same origin will be routed through the same hub node, as well as extensions that incorporate hub location decisions for the leader.
        \item Arc capacities. These could be considered both on hub arcs  and on  access/distribution arcs. Concerning the former, in the MCFP with arc capacities, the commodity flow may be split, and imposing the unsplittable flow already renders the problem NP-hard \citep{Even:1976}. Hence, assuming that the hub network is complete can no longer be done without loss of generality, and different bilevel models need to be considered. Capacity restrictions on access/distribution arcs, or hub or non-hub nodes will affect the optimality conditions of the carriers' subproblems, and would therefore lead to different problem formulations.
        \item Carriers' capacity or fairness constraints: Imposing a maximum service capacity on the carriers would no longer allow to consider each commodity and each carrier independently. Moreover, the regulatory authorities may force the leader to consider some additional fairness allocation constraints, by e.g., minimizing the gap between the maximum and minimum demand allocated to each carrier. 
        \item Dealing with imperfect information: A common criticism in bilevel optimization is the assumption of perfect information. In our context, it may be very difficult for the leader to acquire full knowledge about the reservation prices of the carriers, or the demand, or the routing costs. To overcome these difficulties, one could exploit bilevel optimization under uncertainty \citep{Beck-et-al:2023}, in which the leader has to set the optimal outsourcing fees and allocate non-hubs to carriers, while hedging against the uncertain parameters. Stochastic bilevel models would be considered in case the leader has knowledge about the distribution of uncertain parameters,  otherwise, robust bilevel optimization \citep{Beck-SIAM} would be a model of choice. 
    \end{itemize}
    All the above problem variants render the problem even more difficult, and we consider them worth exploring for future work. They would require the development of innovative solution algorithms to address large-scale instances effectively. It would also be interesting to explore extensions of the applied techniques (discretization and linearization) to other types of problems (e.g., in the context of last mile delivery).
    
    \bibliography{Bibliography_NEW}  
    \newpage
\renewcommand{\thesection}{A-\arabic{section}}
 \setcounter{section}{0}
\begin{section}{Summary of notation}\label{Apdx:Notation}
    \begin{table}[H]
    \caption{General notation used. \label{tab:NotSummary}}
    \resizebox{\textwidth}{!}{\begin{tabular}{lll}
        \textbf{Sets} &  & \multicolumn{1}{c}{} \\
         & $V$ & Node set, $\{1,\,2,\,3,\dots,\,n\}$ \\
         & $H$ & Hub nodes set, $H\subset V$. \\
         & $V\setminus H$ & Non-nub nodes set. \\
         & $A_H$ & Arcs of the backbone network, $\{(i,\,j): i,j\in H\}$ \\
         & $A$ & Arc set, $A_H\cup \{(i,\,j):\left(i\in V\setminus H, j\in H\right) \text{ or } \left(i\in H, j\in V\setminus H\right)\}$ \\
         & $R$ & Index set for commodities. \\
         & $A^r$ & Potential interhub connections for serving commodity $r\in R$ \\
         & $K$ & Index set for carriers $\{1,\,2,\,3,\dots,\,|K|\}$ \\
        \textbf{Additional data} &  &  \\
         & $(o(r),\,d(r),\,w^r,\,b^r)$ & Origin, destination, amount of demand and revenue of commodity $r\in R$ \\
         & $\bar{c}^{rk}_{ij}$ & Reservation price of carrier $k\in K$ to serve commodity $r\in R$ through arc $(i,\,j)\in A$ \\
         & $c^r_{ij}$ & Leader's routing cost for serving commodity $r\in R$ through interhub $(i,\,j)\in A_H$ \\
         & $p^r_i$ & First leg outsourcing fee to serve commodity $r\in R$ through hub $i\in H$ \\
         & $q^r_i$ & Third leg outsourcing fee to serve commodity $r\in R$ through hub $i\in H$
    \end{tabular}}
    \end{table}%
\end{section}

\begin{section}{Data generation for benchmark instances}\label{app:data}
    Below we describe how we have generated the data that was not included in the original benchmark CAB and AP instances.
    \paragraph{Setting the set of hubs $H$.} For every instance, the backbone network is established through the following steps. First, the number of hubs is set to $ |H|=\left\lceil\tau|V|\right\rceil$, where we use the fraction $\tau=0.05$. Then, the center of mass of the input graph is computed, based on the node positions and their respective total demand. Then, a radius is incrementally expanded around this center until the cumulative demand covered reaches a specified threshold $\mu W$, where $W=\sum_{r\in R}w^r$ and $\mu=0.6$. Finally, the $|H|$ nodes with the highest total demand within the current coverage radius are designated as hubs, so $A_H$ as well as $V\setminus H$ and the access/distribution arcs are well defined. The resulting number of hubs, together with the parameters $n$ and $|K|$, are displayed in the first column of Table  \ref{tab:cexp_01_CAB1050_pathbased&implicitpath}.
        
    \paragraph{Setting the routing costs and reservation prices.}
    The leader's routing cost of interhub arc $(i,j)\in A_H$ is set to $c^r_{ij}=w^r\alpha\hat{c}_{ij}$ where $\alpha = 0.5$ is an interhub discount factor.
    We set the carriers reservation prices based on their routing costs through access and distribution arcs. For this, for each non-hub node $i\in V\setminus H$ and carrier $k\in K$, a coefficient $\theta_{ki}$ is drawn from a uniform distribution $U[0.6, 1.2]$, which represents a perturbation associated with the allocation of origin $i$ to carrier $k$. In addition, for each access/distribution arc $(i,j)\in A\setminus A_H$, a coefficient $\chi_{ij}$ is drawn from a uniform distribution $U[0.89, 0.99]$, which represents a perturbation associated with the connection between origin/destination node $i$ and hub $j$. Then, for every commodity $r\in R$ the carriers' routing costs for access arc $(o(r), i)$ are set to  $w^r\theta_{ko(r)}\chi_{o(r)i}\hat{c}_{o(r)i}$, $k\in K$, whereas the carriers' routing costs for distribution arc $(i, d(r))$ are set to  $w^r\theta_{ko(r)}\chi_{id(r)}\hat{c}_{id(r)}$, $k\in K$. 
    Finally, the reservation prices are set to
    \begin{align*}
        \bar{c}^{rk}_{o(r)i} = w^r\theta_{ko(r)}\chi_{o(r)i}\hat{c}_{o(r)i} (1 + \epsilon),\, \qquad 
        \bar{c}^{rk}_{id(r)} = w^r\theta_{kd(r)}\chi_{id(r)}\hat{c}_{id(r)} (1 + \epsilon),
    \end{align*}
    where $\epsilon =0.01$ indicates the profit margin for a carrier in order to accept the service.
    
    \paragraph{Setting the revenues.}
    Following \citet{Alibeyg2018}, the revenues, $b^r$, $r\in R$, are generated as
    \begin{align*}
        b^r=w^r\varphi_r\sum_{(i,j)\in A_H}\frac{\Lambda^r_{ij}}{|A_H|},
    \end{align*}
    where, $\Lambda^r_{ij}=\hat{c}_{o(r)i}+\alpha\hat{c}_{ij}+\hat{c}_{jd(r)}$ is the total unit routing cost of path $o(r)-i-j-d(r)$ relative to the original data set unit costs $\hat c_{ij}$, when applying a discount factor $\alpha = 0.5$ to the interhub arc, and $\varphi_r$ is drawn from a uniform distribution $U[0.25, 0.35]$.
        
    \paragraph{Setting the outsourcing fees for \SDOD$_F$.}
    Two sets of fixed outsourcing fees for first and third legs have been generated for each instance. The first one considers the carriers maximum reservation prices, i.e.,
    \begin{subequations} \label{eq:MaxFees}
        \begin{align}
            &\bar p^r_{i} = \max_{k\in K}\bar c^{rk}_{o(r)i} && i\in H,\,r\in R:o(r)\notin H\\
            &\bar q^r_{i} = \max_{k\in K}\bar c^{rk}_{id(r)} && i\in H,\,r\in R:d(r)\notin H.
        \end{align}
    \end{subequations}
    From the carriers' point of view, these offers would always produce some profit as their values are at least those of the reservation prices. Thus, this variant represents the worst-case scenario for the leader for \SDOD$_F$. 
    The second set considers the average values over the reservation prices of all the carriers, i.e., 
    \begin{subequations} \label{eq:AvgFees}
        \begin{align}
            &\bar p^r_{i} = \frac{\sum_{k\in K}\bar c^{rk}_{o(r)i}}{|K|} && i\in H,\,r\in R:o(r)\notin H\\
            &\bar q^r_{i} = \frac{\sum_{k\in K}\bar c^{rk}_{id(r)}}{|K|} && i\in H,\,r\in R:d(r)\notin H.
        \end{align}    
    \end{subequations}
    This variant represents an intermediate scenario in which offers could be rejected by some carriers even if they could be profitable for the leader.
\end{section}

\begin{section}{Best modeling and algorithmic settings for (EF) and (IF)}\label{Apdx:Tables}

    We first compared the variants of (\emph{EF}) and (\emph{IF}), which use the \emph{BigM} constraints \eqref{const:new_F_ahlpf} and \eqref{const:new_T_ahlpf}, respectively, against the \emph{NoBigM} variants, which use \eqref{2-index-F-i} and \eqref{2-index-T-i}, respectively. Each variant, was tested under three alternative algorithmic settings.  
    The first one (\emph{No-Lazy}) is to enumerate all of the constraints at the beginning; the second setting (\emph{Lazy-Attr}) is to apply Gurobi's \emph{lazy} attribute; and the third approach (\emph{Lazy-Callback}) is to dynamically add violated constraints as lazy constraints using callbacks.

    For both (\emph{EF}) and (\emph{IF})  there are two different families of constraints that can be applied. The \emph{BigM} variant uses \eqref{const:new_F_ahlpf} and \eqref{const:new_T_ahlpf}, whereas the \emph{NoBigM} variant uses \eqref{2-index-F-i} and \eqref{2-index-T-i}.

    \begin{table}[H]
        \caption{Performance comparison for $(EF)$ considering CAB[20-70] and $|K|\in\{3,4\}$.\label{tab:cexp_01_CAB1050_3-Index}}
        \resizebox{\textwidth}{!}{\begin{tabular}{clcccccc}
            Dataset &  & \multicolumn{3}{c}{BigM} & \multicolumn{3}{c}{NoBigM} \\
            $|V|.|K|$ & Approach & t(s) & Nodes & GAP(\%) & t(s) & Nodes & GAP(\%) \\ \hline \hline
            \multirow{3}{*}{20.3} & No-Lazy & 0.48 & 1 & 0.00 & -- & 177129 & -- \\
             & Lazy-Attr. & -- & 176977 & -- &  &  &  \\
             & Lazy-Callback & 0.35 & 1 & 0.00 & -- & 400554 & 0.10 \\ \hline
            \multirow{3}{*}{20.4} & No-Lazy & 0.51 & 1 & 0.00 & -- & 185590 & -- \\
             & Lazy-Attr. & -- & 162302 & -- &  &  &  \\
             & Lazy-Callback & 0.37 & 1 & 0.00 & -- & 230196 & 0.20 \\ \hline
            \multirow{3}{*}{30.3} & No-Lazy & 3.30 & 1 & 0.00 & -- & 84574 & -- \\
             & Lazy-Attr. & -- & 84174 & -- &  &  &  \\
             & Lazy-Callback & 2.00 & 1 & 0.00 & -- & 36822 & 0.50 \\ \hline
            \multirow{3}{*}{30.4} & No-Lazy & 4.49 & 7 & 0.00 & -- & 81681 & -- \\
             & Lazy-Attr. & -- & 69691 & -- &  &  &  \\
             & Lazy-Callback & 2.02 & 1 & 0.00 & -- & 33314 & 0.63 \\ \hline
            \multirow{3}{*}{40.3} & No-Lazy & 19.87 & 1 & 0.00 & -- & 35454 & -- \\
             & Lazy-Attr. & -- & 26614 & -- &  &  &  \\
             & Lazy-Callback & 7.81 & 1 & 0.00 & -- & 6791 & 1.31 \\ \hline
            \multirow{3}{*}{40.4} & No-Lazy & 30.03 & 14 & 0.00 & -- & 41814 & -- \\
             & Lazy-Attr. & -- & 21009 & -- &  &  &  \\
             & Lazy-Callback & 8.21 & 1 & 0.00 & -- & 5937 & 1.80 \\ \hline
            \multirow{3}{*}{50.3} & No-Lazy & 41.58 & 10 & 0.00 & -- & 13882 & -- \\
             & Lazy-Attr. & -- & 6763 & -- &  &  &  \\
             & Lazy-Callback & 25.20 & 2 & 0.00 & -- & 3295 & -- \\ \hline
            \multirow{3}{*}{50.4} & No-Lazy & 62.62 & 64 & 0.00 & -- & 17607 & -- \\
             & Lazy-Attr. & -- & 6823 & -- &  &  &  \\
             & Lazy-Callback & 30.78 & 7 & 0.00 & -- & 3704 & -- \\ \hline
            \multirow{3}{*}{60.3} & No-Lazy & 194.19 & 371 & 0.00 & -- & 5208 & -- \\
             & Lazy-Attr. & -- & 1754 & -- &  &  &  \\
             & Lazy-Callback & 92.05 & 1 & 0.00 & -- & 2568 & -- \\ \hline
            \multirow{3}{*}{60.4} & No-Lazy & 273.44 & 272 & 0.00 & -- & 6500 & -- \\
             & Lazy-Attr. & -- & 2136 & -- &  &  &  \\
             & Lazy-Callback & 87.15 & 5 & 0.00 & -- & 2122 & -- \\ \hline
            \multirow{3}{*}{70.3} & No-Lazy & -- & 1 & 0.10 & -- & 4384 & -- \\
             & Lazy-Attr. & -- & 1924 & -- &  &  &  \\
             & Lazy-Callback & 248.80 & 5 & 0.00 & -- & 1128 & -- \\ \hline
            \multirow{3}{*}{70.4} & No-Lazy & -- & 1 & 0.17 & -- & 3317 & -- \\
             & Lazy-Attr. & -- & 1655 & -- &  &  &  \\
             & Lazy-Callback & 268.26 & 35 & 0.00 & -- & 1128 & --
        \end{tabular}}\\
        \footnotesize{--: Time limit reached or no feasible solution found.}
    \end{table} 
    \begin{table}[H]
        \caption{Performance comparison for $(IF)$ considering CAB[20-70] and $|K|\in\{3,4\}$.\label{tab:cexp_01_CAB1050_2-Index}}
        \resizebox{\textwidth}{!}{\begin{tabular}{clcccccc}
                Dataset &  & \multicolumn{3}{c}{BigM} & \multicolumn{3}{c}{NoBigM} \\
                $|V|.|K|$ & Approach & t(s) & Nodes & GAP(\%) & t(s) & Nodes & GAP(\%) \\ \hline \hline
                \multirow{3}{*}{20.3} & No-Lazy & 67.04 & 50485 & 0.00 & -- & 151244 & -- \\
                 & Lazy-Attr. & -- & 126740 & -- & \multicolumn{1}{l}{} & \multicolumn{1}{l}{} & \multicolumn{1}{l}{} \\
                 & Lazy-Callback & -- & 652558 & 0.11 & -- & 450107 & 0.14 \\ \hline
                \multirow{3}{*}{20.4} & No-Lazy & 244.50 & 137555 & 0.01 & -- & 190940 & -- \\
                 & Lazy-Attr. & -- & 146133 & -- & \multicolumn{1}{l}{} & \multicolumn{1}{l}{} & \multicolumn{1}{l}{} \\
                 & Lazy-Callback & -- & 471470 & 0.44 & -- & 322586 & 0.79 \\ \hline
                \multirow{3}{*}{30.3} & No-Lazy & -- & 33854 & 0.25 & -- & 63200 & -- \\
                 & Lazy-Attr. & -- & 56553 & -- & \multicolumn{1}{l}{} & \multicolumn{1}{l}{} & \multicolumn{1}{l}{} \\
                 & Lazy-Callback & -- & 121391 & 0.41 & -- & 54444 & 1.41 \\ \hline
                \multirow{3}{*}{30.4} & No-Lazy & -- & 33237 & 0.82 & -- & 63200 & -- \\
                 & Lazy-Attr. & -- & 68018 & -- & \multicolumn{1}{l}{} & \multicolumn{1}{l}{} & \multicolumn{1}{l}{} \\
                 & Lazy-Callback & -- & 103822 & 1.04 & -- & 36313 & 1.42 \\ \hline
                \multirow{3}{*}{40.3} & No-Lazy & -- & 15634 & 0.87 & -- & 53648 & -- \\
                 & Lazy-Attr. & -- & 40093 & -- & \multicolumn{1}{l}{} & \multicolumn{1}{l}{} & \multicolumn{1}{l}{} \\
                 & Lazy-Callback & -- & 22599 & 0.67 & -- & 6783 & -- \\ \hline
                \multirow{3}{*}{40.4} & No-Lazy & -- & 11364 & 1.34 & -- & 46875 & -- \\
                 & Lazy-Attr. & -- & 46184 & -- & \multicolumn{1}{l}{} & \multicolumn{1}{l}{} & \multicolumn{1}{l}{} \\
                 & Lazy-Callback & -- & 18005 & 1.01 & -- & 7635 & 6.30 \\ \hline
                \multirow{3}{*}{50.3} & No-Lazy & -- & 4032 & 1.65 & -- & 16583 & -- \\
                 & Lazy-Attr. & -- & 6999 & -- & \multicolumn{1}{l}{} & \multicolumn{1}{l}{} & \multicolumn{1}{l}{} \\
                 & Lazy-Callback & -- & 6157 & 0.76 & -- & 4251 & -- \\ \hline
                \multirow{3}{*}{50.4} & No-Lazy & -- & 3474 & 1.95 & -- & 12381 & -- \\
                 & Lazy-Attr. & -- & 9572 & -- & \multicolumn{1}{l}{} & \multicolumn{1}{l}{} & \multicolumn{1}{l}{} \\
                 & Lazy-Callback & -- & 3484 & 1.36 & -- & 5172 & -- \\ \hline
                \multirow{3}{*}{60.3} & No-Lazy & -- & 1796 & 1.51 & -- & 5507 & -- \\
                 & Lazy-Attr. & -- & 2870 & -- & \multicolumn{1}{l}{} & \multicolumn{1}{l}{} & \multicolumn{1}{l}{} \\
                 & Lazy-Callback & -- & 1712 & 1.71 & -- & 1653 & -- \\ \hline
                \multirow{3}{*}{60.4} & No-Lazy & -- & 1893 & 2.23 & -- & 3261 & -- \\
                 & Lazy-Attr. & -- & 3311 & -- & \multicolumn{1}{l}{} & \multicolumn{1}{l}{} & \multicolumn{1}{l}{} \\
                 & Lazy-Callback & -- & 1115 & 2.07 & -- & 1526 & -- \\ \hline
                \multirow{3}{*}{70.3} & No-Lazy & -- & 381 & 1.39 & -- & 2823 & -- \\
                 & Lazy-Attr. & -- & 1481 & -- & \multicolumn{1}{l}{} & \multicolumn{1}{l}{} & \multicolumn{1}{l}{} \\
                 & Lazy-Callback & -- & 325 & 56.09 & -- & 1101 & -- \\ \hline
                \multirow{3}{*}{70.4} & No-Lazy & -- & 241 & 2.56 & -- & 2543 & -- \\
                 & Lazy-Attr. & -- & 1274 & -- & \multicolumn{1}{l}{} & \multicolumn{1}{l}{} & \multicolumn{1}{l}{} \\
                 & Lazy-Callback & -- & 98 & 43.95 & -- & 1022 & --
            \end{tabular}}\\
            \footnotesize{--: Time limit reached or no feasible solution found.}
    \end{table}
    Tables \ref{tab:cexp_01_CAB1050_3-Index}-\ref{tab:cexp_01_CAB1050_2-Index} summarize the obtained results for the CAB instances for (\emph{EF}) and (\emph{IF}), respectively. The tables have a block of three columns for each modeling variant (\emph{BigM}, \emph{NoBigM}). Columns $t(s)$ give computing times in seconds (or a dash when the maximum time limit was reached), columns \textit{Nodes} show the number of nodes explored in the enumeration tree, and  columns \textit{GAP$\%$} the percent optimality gap at termination, calculated as $|\frac{UB-LB}{UB}|\times 100$. The rows of the tables display average results over the ten instances of the corresponding dimensions for the respective algorithmic approach. As can be seen, the best results are obtained with the \emph{BigM} variant when solved with the \emph{Lazy-Callback} strategy. Thus, this is the variant that is used for both (\emph{EF}) and (\emph{IF}) in the remainder of this paper.
\end{section}
\begin{section}{Comparison of \SDOD$F$ and $r$\SDOD$_F$ for \emph{average} and \emph{maximum} fixed outsourcing fees}\label{ap:fixed-out}
    Table \ref{tab:FFees_EP_VarRem} displays average percentages of the number of variables fixed in the preprocessing phase of (\emph{EP}$_F$) and (\emph{IP}$_F$).
    \begin{table}[H]
        \caption{Average percentage of preprocessed variables for \SDOD$_F$ and $r$\SDOD$_F$ \rojooscuro{with} fixed outsourcing fees. \label{tab:FFees_EP_VarRem}}
        {\resizebox{\textwidth}{!}{\begin{tabular}{cl|rrrrrrrrr|rrr}
    \multirow{3}{*}{Dataset} & \multicolumn{1}{c|}{\multirow{3}{*}{$|V|.|H|$}} & \multicolumn{9}{c|}{(EP$_F$)} & \multicolumn{3}{c}{(IP$_F$)} \\
     & \multicolumn{1}{c|}{} & \multicolumn{3}{c|}{$x$ variables} & \multicolumn{3}{c|}{$f$ variables} & \multicolumn{3}{c|}{$t$ variables} & \multicolumn{3}{c}{$\pi$ variables} \\
     & \multicolumn{1}{c|}{} & Initial & Remaining & \multicolumn{1}{r|}{$\%$} & Initial & Remaining & \multicolumn{1}{r|}{$\%$} & Initial & Remaining & $\%$ & Initial & Remaining & $\%$ \\ \hline \hline
    \multirow{6}{*}{CAB} & 20.1 & 10000 & 1417 & \multicolumn{1}{r|}{85.83} & 7000 & 4025 & \multicolumn{1}{r|}{42.83} & 7000 & 4025 & 42.83 & 4818 & 877 & 83.71 \\
     & 30.2 & 22500 & 2329 & \multicolumn{1}{r|}{89.65} & 15750 & 9128 & \multicolumn{1}{r|}{42.38} & 15750 & 9128 & 42.38 & 10459 & 1253 & 86.97 \\
     & 40.2 & 40000 & 4884 & \multicolumn{1}{r|}{87.79} & 28000 & 16370 & \multicolumn{1}{r|}{41.83} & 28000 & 16370 & 41.83 & 18929 & 1752 & 90.80 \\
     & 50.3 & 62500 & 6798 & \multicolumn{1}{r|}{89.12} & 43750 & 25788 & \multicolumn{1}{r|}{41.37} & 43750 & 25788 & 41.37 & 29043 & 2836 & 90.41 \\
     & 60.3 & 90000 & 8804 & \multicolumn{1}{r|}{90.22} & 63000 & 37530 & \multicolumn{1}{r|}{40.73} & 63000 & 37530 & 40.73 & 42333 & 4866 & 88.59 \\
     & 70.4 & 122500 & 10835 & \multicolumn{1}{r|}{91.15} & 85750 & 51030 & \multicolumn{1}{r|}{40.77} & 85750 & 51030 & 40.77 & 56920 & 7062 & 87.18 \\
    \hdashline\multirow{5}{*}{AP} & 100.5 & 250000 & 1221 & \multicolumn{1}{r|}{99.51} & 200000 & 123233 & \multicolumn{1}{r|}{39.07} & 200000 & 123233 & 39.07 & 155220 & 22857 & 85.30 \\
     & 120.6 & 360000 & 1375 & \multicolumn{1}{r|}{99.62} & 288000 & 179120 & \multicolumn{1}{r|}{38.48} & 288000 & 179120 & 38.48 & 223282 & 33524 & 84.41 \\
     & 140.7 & 490000 & 1549 & \multicolumn{1}{r|}{99.68} & 392000 & 245653 & \multicolumn{1}{r|}{38.01} & 392000 & 245653 & 38.01 & 303683 & 62491 & 79.71 \\
     & 160.8 & 640000 & 1913 & \multicolumn{1}{r|}{99.70} & 512000 & 320533 & \multicolumn{1}{r|}{38.07} & 512000 & 320533 & 38.07 & 396424 & 71252 & 82.46 \\
     & 180.9 & 810000 & 2081 & \multicolumn{1}{r|}{99.74} & 648000 & 405540 & \multicolumn{1}{r|}{38.09} & 648000 & 405540 & 38.09 & 501504 & 99156 & 79.75 \\ \hline
    \end{tabular}}}
    \end{table}
    Tables \ref{tab:Granular_max}--\ref{tab:Granular_avg} show  disaggregated results on the variables eliminated in $(EP_F)$ and $(IP_F)$  due to the preprocessing for fixed outsourcing fees used to compute Table \ref{tab:FFees_EP_VarRem}. 
    \begin{table}[H]
        \caption{Number of preprocessed variables for \SDOD$_F$ and $r$\SDOD$_F$  with maximum fixed outsourcing fees.\label{tab:Granular_max}}     
        {\resizebox{\textwidth}{!}{\begin{tabular}{cc|rrrrrrrrr|rrr|rrr}
     \multirow{3}{*}{Dataset} & \multirow{3}{*}{$|V|.|K|$} & \multicolumn{9}{c|}{$(EP)$ \& $(rEP)$} & \multicolumn{3}{c|}{($IP_F$)} & \multicolumn{3}{c}{($rIP_F$) } \\
     &  & \multicolumn{3}{c|}{$x$ variables} & \multicolumn{3}{c|}{$f$ variables} & \multicolumn{3}{c|}{$t$ variables} & \multicolumn{3}{c|}{$\pi$ variables} & \multicolumn{3}{c}{$\pi$ variables} \\
     &  & Initial & Remaining & \multicolumn{1}{r|}{$\%$} & Initial & Remaining & \multicolumn{1}{r|}{$\%$} & Initial & Remaining & $\%$ & Initial & Remaining & $\%$ & Initial & Remaining & $\%$ \\ \hline \hline
    \multirow{12}{*}{CAB} & 20.3 & 10000 & 1419 & \multicolumn{1}{r|}{85.81} & 6000 & 3800 & \multicolumn{1}{r|}{36.67} & 6000 & 3800 & 36.67 & 3512 & 510 & 85.48 & 3512 & 510 & 85.48 \\
     & 20.4 & 10000 & 1420 & \multicolumn{1}{r|}{85.80} & 8000 & 5700 & \multicolumn{1}{r|}{28.75} & 8000 & 5700 & 28.75 & 6124 & 2156 & 64.79 & 6124 & 2156 & 64.79 \\
     & 30.3 & 22500 & 2333 & \multicolumn{1}{r|}{89.63} & 13500 & 8700 & \multicolumn{1}{r|}{35.56} & 13500 & 8700 & 35.56 & 7622 & 1189 & 84.40 & 7622 & 2393 & 68.60 \\
     & 30.4 & 22500 & 2334 & \multicolumn{1}{r|}{89.63} & 18000 & 13050 & \multicolumn{1}{r|}{27.50} & 18000 & 13050 & 27.50 & 13296 & 671 & 94.95 & 13296 & 3106 & 76.64 \\
     & 40.3 & 40000 & 4899 & \multicolumn{1}{r|}{87.75} & 24000 & 15600 & \multicolumn{1}{r|}{35.00} & 24000 & 15600 & 35.00 & 13752 & 1246 & 90.94 & 13752 & 2006 & 85.41 \\
     & 40.4 & 40000 & 4904 & \multicolumn{1}{r|}{87.74} & 32000 & 23400 & \multicolumn{1}{r|}{26.88} & 32000 & 23400 & 26.88 & 24106 & 2215 & 90.81 & 24106 & 3818 & 84.16 \\
     & 50.3 & 62500 & 6818 & \multicolumn{1}{r|}{89.09} & 37500 & 24500 & \multicolumn{1}{r|}{34.67} & 37500 & 24500 & 34.67 & 21110 & 471 & 97.77 & 21110 & 4978 & 76.42 \\
     & 50.4 & 62500 & 6825 & \multicolumn{1}{r|}{89.08} & 50000 & 36750 & \multicolumn{1}{r|}{26.50} & 50000 & 36750 & 26.50 & 36976 & 1035 & 97.20 & 36976 & 9617 & 73.99 \\
     & 60.3 & 90000 & 8828 & \multicolumn{1}{r|}{90.19} & 54000 & 35400 & \multicolumn{1}{r|}{34.44} & 54000 & 35400 & 34.44 & 30720 & 974 & 96.83 & 30720 & 9191 & 70.08 \\
     & 60.4 & 90000 & 8837 & \multicolumn{1}{r|}{90.18} & 72000 & 53100 & \multicolumn{1}{r|}{26.25} & 72000 & 53100 & 26.25 & 53946 & 1349 & 97.50 & 53946 & 17122 & 68.26 \\
     & 70.3 & 122500 & 10864 & \multicolumn{1}{r|}{91.13} & 73500 & 48300 & \multicolumn{1}{r|}{34.29} & 73500 & 48300 & 34.29 & 41326 & 1091 & 97.36 & 41326 & 14848 & 64.07 \\
     & 70.4 & 122500 & 10875 & \multicolumn{1}{r|}{91.12} & 98000 & 72450 & \multicolumn{1}{r|}{26.07} & 98000 & 72450 & 26.07 & 72514 & 2531 & 96.51 & 72514 & 17229 & 76.24 \\
    \hdashline\multirow{15}{*}{AP} & 100.3 & 250000 & 1200 & \multicolumn{1}{r|}{99.52} & 150000 & 99000 & \multicolumn{1}{r|}{34.00} & 150000 & 99000 & 34.00 & 84840 & 3359 & 96.04 & 84840 & 25136 & 70.37 \\
     & 100.4 & 250000 & 1196 & \multicolumn{1}{r|}{99.52} & 200000 & 148500 & \multicolumn{1}{r|}{25.75} & 200000 & 148500 & 25.75 & 149200 & 4612 & 96.91 & 149200 & 46780 & 68.65 \\
     & 100.5 & 250000 & 1195 & \multicolumn{1}{r|}{99.52} & 250000 & 198000 & \multicolumn{1}{r|}{20.80} & 250000 & 198000 & 20.80 & 231620 & 9210 & 96.02 & 231620 & 64610 & 72.11 \\
     & 120.3 & 360000 & 1351 & \multicolumn{1}{r|}{99.62} & 216000 & 142800 & \multicolumn{1}{r|}{33.89} & 216000 & 142800 & 33.89 & 121992 & 6251 & 94.88 & 121992 & 42543 & 65.13 \\
     & 120.4 & 360000 & 1346 & \multicolumn{1}{r|}{99.63} & 288000 & 214200 & \multicolumn{1}{r|}{25.63} & 288000 & 214200 & 25.63 & 214614 & 5311 & 97.53 & 214614 & 63918 & 70.22 \\
     & 120.5 & 360000 & 1344 & \multicolumn{1}{r|}{99.63} & 360000 & 285600 & \multicolumn{1}{r|}{20.67} & 360000 & 285600 & 20.67 & 333240 & 14130 & 95.76 & 333240 & 87150 & 73.85 \\
     & 140.3 & 490000 & 1523 & \multicolumn{1}{r|}{99.69} & 294000 & 194600 & \multicolumn{1}{r|}{33.81} & 294000 & 194600 & 33.81 & 165872 & 7697 & 95.36 & 165872 & 65526 & 60.50 \\
     & 140.4 & 490000 & 1515 & \multicolumn{1}{r|}{99.69} & 392000 & 291900 & \multicolumn{1}{r|}{25.54} & 392000 & 291900 & 25.54 & 291886 & 11507 & 96.06 & 291886 & 113706 & 61.04 \\
     & 140.5 & 490000 & 1513 & \multicolumn{1}{r|}{99.69} & 490000 & 389200 & \multicolumn{1}{r|}{20.57} & 490000 & 389200 & 20.57 & 453292 & 20215 & 95.54 & 453292 & 188417 & 58.43 \\
     & 160.3 & 640000 & 1884 & \multicolumn{1}{r|}{99.71} & 384000 & 254400 & \multicolumn{1}{r|}{33.75} & 384000 & 254400 & 33.75 & 216480 & 6254 & 97.11 & 216480 & 76451 & 64.68 \\
     & 160.4 & 640000 & 1868 & \multicolumn{1}{r|}{99.71} & 512000 & 381600 & \multicolumn{1}{r|}{25.47} & 512000 & 381600 & 25.47 & 381016 & 8521 & 97.76 & 381016 & 123632 & 67.55 \\
     & 160.5 & 640000 & 1861 & \multicolumn{1}{r|}{99.71} & 640000 & 508800 & \multicolumn{1}{r|}{20.50} & 640000 & 508800 & 20.50 & 591776 & 11065 & 98.13 & 591776 & 242826 & 58.97 \\
     & 180.3 & 810000 & 2048 & \multicolumn{1}{r|}{99.75} & 486000 & 322200 & \multicolumn{1}{r|}{33.70} & 486000 & 322200 & 33.70 & 273816 & 8652 & 96.84 & 273816 & 121125 & 55.76 \\
     & 180.4 & 810000 & 2032 & \multicolumn{1}{r|}{99.75} & 648000 & 483300 & \multicolumn{1}{r|}{25.42} & 648000 & 483300 & 25.42 & 482004 & 13700 & 97.16 & 482004 & 205400 & 57.39 \\
     & 180.5 & 810000 & 2024 & \multicolumn{1}{r|}{99.75} & 810000 & 644400 & \multicolumn{1}{r|}{20.44} & 810000 & 644400 & 20.44 & 748692 & 22388 & 97.01 & 748692 & 278762 & 62.77 \\ \hline
    \end{tabular}}}
    \end{table}
    \begin{table}[H]
        \caption{Number of preprocessed variables for \SDOD$_F$ and $r$\SDOD$_F$ with average fixed outsourcing fees.\label{tab:Granular_avg}}
        {\resizebox{\textwidth}{!}{\begin{tabular}{cc|rrrrrrrrr|rrr|rrr}
         \multirow{3}{*}{Dataset} & \multirow{3}{*}{$|V|.|K|$} & \multicolumn{9}{c|}{$(EP)$ \& $(rEP)$} & \multicolumn{3}{c|}{($IP_F$)} & \multicolumn{3}{c}{($rIP_F$) } \\
     &  & \multicolumn{3}{c|}{$x$ variables} & \multicolumn{3}{c|}{$f$ variables} & \multicolumn{3}{c|}{$t$ variables} & \multicolumn{3}{c|}{$\pi$ variables} & \multicolumn{3}{c}{$\pi$ variables} \\
     &  & Initial & Remaining & \multicolumn{1}{r|}{$\%$} & Initial & Remaining & \multicolumn{1}{r|}{$\%$} & Initial & Remaining & $\%$ & Initial & Remaining & $\%$ & Initial & Remaining & $\%$ \\ \hline \hline
    \multirow{12}{*}{CAB} & 20.3 & 10000 & 1415 & \multicolumn{1}{r|}{85.85} & 6000 & 2780 & \multicolumn{1}{r|}{53.67} & 6000 & 2780 & 53.67 & 3512 & 143 & 95.93 & 3512 & 143 & 95.93 \\
     & 20.4 & 10000 & 1415 & \multicolumn{1}{r|}{85.85} & 8000 & 3820 & \multicolumn{1}{r|}{52.25} & 8000 & 3820 & 52.25 & 6124 & 697 & 88.62 & 6124 & 697 & 88.62 \\
     & 30.3 & 22500 & 2325 & \multicolumn{1}{r|}{89.67} & 13500 & 6240 & \multicolumn{1}{r|}{53.78} & 13500 & 6240 & 53.78 & 7622 & 659 & 91.35 & 7622 & 911 & 88.05 \\
     & 30.4 & 22500 & 2324 & \multicolumn{1}{r|}{89.67} & 18000 & 8520 & \multicolumn{1}{r|}{52.67} & 18000 & 8520 & 52.67 & 13296 & 289 & 97.83 & 13296 & 806 & 93.94 \\
     & 40.3 & 40000 & 4868 & \multicolumn{1}{r|}{87.83} & 24000 & 11320 & \multicolumn{1}{r|}{52.83} & 24000 & 11320 & 52.83 & 13752 & 747 & 94.57 & 13752 & 953 & 93.07 \\
     & 40.4 & 40000 & 4865 & \multicolumn{1}{r|}{87.84} & 32000 & 15160 & \multicolumn{1}{r|}{52.63} & 32000 & 15160 & 52.63 & 24106 & 1341 & 94.44 & 24106 & 1691 & 92.99 \\
     & 50.3 & 62500 & 6776 & \multicolumn{1}{r|}{89.16} & 37500 & 17850 & \multicolumn{1}{r|}{52.40} & 37500 & 17850 & 52.40 & 21110 & 408 & 98.07 & 21110 & 1695 & 91.97 \\
     & 50.4 & 62500 & 6772 & \multicolumn{1}{r|}{89.16} & 50000 & 24050 & \multicolumn{1}{r|}{51.90} & 50000 & 24050 & 51.90 & 36976 & 876 & 97.63 & 36976 & 3606 & 90.25 \\
     & 60.3 & 90000 & 8777 & \multicolumn{1}{r|}{90.25} & 54000 & 26340 & \multicolumn{1}{r|}{51.22} & 54000 & 26340 & 51.22 & 30720 & 648 & 97.89 & 30720 & 2821 & 90.82 \\
     & 60.4 & 90000 & 8773 & \multicolumn{1}{r|}{90.25} & 72000 & 35280 & \multicolumn{1}{r|}{51.00} & 72000 & 35280 & 51.00 & 53946 & 1237 & 97.71 & 53946 & 5584 & 89.65 \\
     & 70.3 & 122500 & 10802 & \multicolumn{1}{r|}{91.18} & 73500 & 35910 & \multicolumn{1}{r|}{51.14} & 73500 & 35910 & 51.14 & 41326 & 1187 & 97.13 & 41326 & 6573 & 84.09 \\
     & 70.4 & 122500 & 10800 & \multicolumn{1}{r|}{91.18} & 98000 & 47460 & \multicolumn{1}{r|}{51.57} & 98000 & 47460 & 51.57 & 72514 & 3331 & 95.41 & 72514 & 9703 & 86.62 \\
    \hdashline\multirow{15}{*}{AP} & 100.3 & 250000 & 1245 & \multicolumn{1}{r|}{99.50} & 150000 & 69700 & \multicolumn{1}{r|}{53.53} & 150000 & 69700 & 53.53 & 84840 & 4013 & 95.27 & 84840 & 16162 & 80.95 \\
     & 100.4 & 250000 & 1241 & \multicolumn{1}{r|}{99.50} & 200000 & 99700 & \multicolumn{1}{r|}{50.15} & 200000 & 99700 & 50.15 & 149200 & 7494 & 94.98 & 149200 & 31624 & 78.80 \\
     & 100.5 & 250000 & 1247 & \multicolumn{1}{r|}{99.50} & 250000 & 124500 & \multicolumn{1}{r|}{50.20} & 250000 & 124500 & 50.20 & 231620 & 14295 & 93.83 & 231620 & 46983 & 79.72 \\
     & 120.3 & 360000 & 1404 & \multicolumn{1}{r|}{99.61} & 216000 & 103080 & \multicolumn{1}{r|}{52.28} & 216000 & 103080 & 52.28 & 121992 & 7874 & 93.55 & 121992 & 31171 & 74.45 \\
     & 120.4 & 360000 & 1400 & \multicolumn{1}{r|}{99.61} & 288000 & 146520 & \multicolumn{1}{r|}{49.13} & 288000 & 146520 & 49.13 & 214614 & 9027 & 95.79 & 214614 & 46804 & 78.19 \\
     & 120.5 & 360000 & 1405 & \multicolumn{1}{r|}{99.61} & 360000 & 182520 & \multicolumn{1}{r|}{49.30} & 360000 & 182520 & 49.30 & 333240 & 18274 & 94.52 & 333240 & 69839 & 79.04 \\
     & 140.3 & 490000 & 1583 & \multicolumn{1}{r|}{99.68} & 294000 & 142380 & \multicolumn{1}{r|}{51.57} & 294000 & 142380 & 51.57 & 165872 & 7609 & 95.41 & 165872 & 48634 & 70.68 \\
     & 140.4 & 490000 & 1576 & \multicolumn{1}{r|}{99.68} & 392000 & 204260 & \multicolumn{1}{r|}{47.89} & 392000 & 204260 & 47.89 & 291886 & 19589 & 93.29 & 291886 & 88901 & 69.54 \\
     & 140.5 & 490000 & 1582 & \multicolumn{1}{r|}{99.68} & 490000 & 251580 & \multicolumn{1}{r|}{48.66} & 490000 & 251580 & 48.66 & 453292 & 26919 & 94.06 & 453292 & 151172 & 66.65 \\
     & 160.3 & 640000 & 1958 & \multicolumn{1}{r|}{99.69} & 384000 & 185600 & \multicolumn{1}{r|}{51.67} & 384000 & 185600 & 51.67 & 216480 & 6307 & 97.09 & 216480 & 57672 & 73.36 \\
     & 160.4 & 640000 & 1949 & \multicolumn{1}{r|}{99.70} & 512000 & 265440 & \multicolumn{1}{r|}{48.16} & 512000 & 265440 & 48.16 & 381016 & 15320 & 95.98 & 381016 & 98784 & 74.07 \\
     & 160.5 & 640000 & 1955 & \multicolumn{1}{r|}{99.69} & 640000 & 327360 & \multicolumn{1}{r|}{48.85} & 640000 & 327360 & 48.85 & 591776 & 20612 & 96.52 & 591776 & 187575 & 68.30 \\
     & 180.3 & 810000 & 2131 & \multicolumn{1}{r|}{99.74} & 486000 & 235440 & \multicolumn{1}{r|}{51.56} & 486000 & 235440 & 51.56 & 273816 & 8082 & 97.05 & 273816 & 95734 & 65.04 \\
     & 180.4 & 810000 & 2123 & \multicolumn{1}{r|}{99.74} & 648000 & 331200 & \multicolumn{1}{r|}{48.89} & 648000 & 331200 & 48.89 & 482004 & 21289 & 95.58 & 482004 & 165891 & 65.58 \\
     & 180.5 & 810000 & 2130 & \multicolumn{1}{r|}{99.74} & 810000 & 416700 & \multicolumn{1}{r|}{48.56} & 810000 & 416700 & 48.56 & 748692 & 30260 & 95.96 & 748692 & 218583 & 70.80 \\ \hline
    \end{tabular}}}%
    \end{table}%
    Tables \ref{tab:AvgFFees_CAB_COMP}--\ref{tab:MaxFFees_CAB_COMP} display results of \SDOD$_F$ and  $r$\SDOD$_F$ for the (preprocessed) benchmark instances from the CAB dataset, using the \emph{maximum} and \emph{average} fixed outsourcing fees, respectively.
    \begin{table}[H]
            \caption{Comparison of  \SDOD$_F$ and $r$\SDOD$_F$ formulations on CAB[20-70] with $|K| \in \{3,\,4\}$ with average fixed outsourcing fees.\label{tab:AvgFFees_CAB_COMP}}
            \resizebox{\textwidth}{!}{\begin{tabular}{lrrrrrr|rrrrrr}
                & \multicolumn{3}{c}{(\emph{EP}$_F$)} & \multicolumn{3}{c|}{($r$\emph{EP}$_F$)} & \multicolumn{3}{c}{(\emph{IP}$_F$)} & \multicolumn{3}{c}{($r$\emph{IP}$_F$)} \\
                \multicolumn{1}{c}{$|V|.|K|.|H|$} & \multicolumn{1}{c}{t(s)} & \multicolumn{1}{c}{Nodes} & \multicolumn{1}{c}{GAP} & \multicolumn{1}{c}{t(s)} & \multicolumn{1}{c}{Nodes} & \multicolumn{1}{c|}{GAP} & \multicolumn{1}{c}{t(s)} & \multicolumn{1}{c}{Nodes} & \multicolumn{1}{c}{GAP} & \multicolumn{1}{c}{t(s)} & \multicolumn{1}{c}{Nodes} & \multicolumn{1}{c}{GAP} \\ \hline \hline
                20.3.1 & 0.02 & 1.0 & 0.00 & 0.12 & 1.0 & 0.00 & 0.01 & 1.0 & 0.00 & 0.02 & 1.0 & 0.00 \\
                20.4.1 & 0.02 & 1.0 & 0.00 & 0.09 & 1.0 & 0.00 & 0.01 & 1.0 & 0.00 & 0.02 & 1.0 & 0.00 \\
                30.3.2 & 0.09 & 1.0 & 0.00 & 0.57 & 1.0 & 0.00 & 0.01 & 1.0 & 0.00 & 0.02 & 1.0 & 0.00 \\
                30.4.2 & 0.13 & 1.0 & 0.00 & 0.83 & 1.2 & 0.00 & 0.03 & 1.0 & 0.00 & 0.04 & 1.0 & 0.00 \\
                40.3.2 & 0.28 & 1.0 & 0.00 & 3.07 & 1.0 & 0.00 & 0.03 & 1.0 & 0.00 & 0.05 & 1.0 & 0.00 \\
                40.4.2 & 0.45 & 1.0 & 0.00 & 6.68 & 25.2 & 0.00 & 0.05 & 1.0 & 0.00 & 0.08 & 1.0 & 0.00 \\
                50.3.3 & 0.63 & 1.0 & 0.00 & 10.82 & 58.6 & 0.00 & 0.04 & 1.0 & 0.00 & 0.08 & 1.0 & 0.00 \\
                50.4.3 & 0.94 & 1.0 & 0.00 & 21.74 & 132.8 & 0.00 & 0.08 & 1.0 & 0.00 & 0.13 & 1.0 & 0.00 \\
                60.3.3 & 1.39 & 1.0 & 0.00 & 54.58 & 220.0 & 0.00 & 0.06 & 1.0 & 0.00 & 0.12 & 1.0 & 0.00 \\
                60.4.3 & 2.02 & 1.0 & 0.00 & 64.31 & 354.4 & 0.00 & 0.10 & 1.0 & 0.00 & 0.21 & 1.0 & 0.00 \\
                70.3.4 & 2.50 & 1.0 & 0.00 & --$\quad$ & 733.2 & 0.22 & 0.08 & 1.0 & 0.00 & 0.17 & 1.0 & 0.00 \\
                70.4.4 & 4.07 & 1.0 & 0.00 & --$\quad$ & 690.2 & 0.95 & 0.14 & 1.0 & 0.00 & 0.36 & 1.0 & 0.00 \\
                \hline
            \end{tabular}}
        \end{table}
    \begin{table}[H]
        \caption{Comparison of \SDOD$_F$ and $r$\SDOD$_F$ formulations on CAB[20-70] with $|K| \in \{3,\,4\}$ with maximum fixed outsourcing fees. \label{tab:MaxFFees_CAB_COMP}}
        \resizebox{\textwidth}{!}{\begin{tabular}{lrrrrrr|rrrrrr}
                 & \multicolumn{3}{c}{(\emph{EP}$_F$)} & \multicolumn{3}{c|}{($r$\emph{EP}$_F$)} & \multicolumn{3}{c}{(\emph{IP}$_F$)} & \multicolumn{3}{c}{($r$\emph{IP}$_F$)} \\
                \multicolumn{1}{c}{$|V|.|K|.|H|$} & \multicolumn{1}{c}{t(s)} & \multicolumn{1}{c}{Nodes} & \multicolumn{1}{c}{GAP} & \multicolumn{1}{c}{t(s)} & \multicolumn{1}{c}{Nodes} & \multicolumn{1}{c|}{GAP} & \multicolumn{1}{c}{t(s)} & \multicolumn{1}{c}{Nodes} & \multicolumn{1}{c}{GAP} & \multicolumn{1}{c}{t(s)} & \multicolumn{1}{c}{Nodes} & \multicolumn{1}{c}{GAP} \\ \hline \hline
                20.3.1 & 0.02 & 1.0 & 0.00 & 0.01 & 1.0 & 0.00 & 0.01 & 1.0 & 0.00 & 0.01 & 1.0 & 0.00 \\
                20.4.1 & 0.02 & 1.0 & 0.00 & 0.01 & 1.0 & 0.00 & 0.02 & 1.0 & 0.00 & 0.01 & 1.0 & 0.00 \\
                30.3.2 & 0.09 & 1.0 & 0.00 & 0.03 & 1.0 & 0.00 & 0.02 & 1.0 & 0.00 & 0.02 & 1.0 & 0.00 \\
                30.4.2 & 0.13 & 1.0 & 0.00 & 0.03 & 1.0 & 0.00 & 0.03 & 1.0 & 0.00 & 0.04 & 1.0 & 0.00 \\
                40.3.2 & 0.28 & 1.0 & 0.00 & 0.07 & 1.0 & 0.00 & 0.03 & 1.0 & 0.00 & 0.04 & 1.0 & 0.00 \\
                40.4.2 & 0.46 & 1.0 & 0.00 & 0.08 & 1.0 & 0.00 & 0.05 & 1.0 & 0.00 & 0.07 & 1.0 & 0.00 \\
                50.3.3 & 0.70 & 1.0 & 0.00 & 0.14 & 1.0 & 0.00 & 0.04 & 1.0 & 0.00 & 0.07 & 1.0 & 0.00 \\
                50.4.3 & 1.10 & 1.0 & 0.00 & 0.17 & 1.0 & 0.00 & 0.07 & 1.0 & 0.00 & 0.12 & 1.0 & 0.00 \\
                60.3.3 & 1.56 & 1.0 & 0.00 & 0.28 & 1.0 & 0.00 & 0.05 & 1.0 & 0.00 & 0.11 & 1.0 & 0.00 \\
                60.4.3 & 2.35 & 1.0 & 0.00 & 0.33 & 1.0 & 0.00 & 0.10 & 1.0 & 0.00 & 0.18 & 1.0 & 0.00 \\
                70.3.4 & 3.00 & 1.0 & 0.00 & 0.47 & 1.0 & 0.00 & 0.08 & 1.0 & 0.00 & 0.15 & 1.0 & 0.00 \\
                70.4.4 & 4.81 & 1.0 & 0.00 & 0.55 & 1.0 & 0.00 & 0.13 & 1.0 & 0.00 & 0.32 & 1.0 & 0.00 \\
                \hline
            \end{tabular}}
    \end{table}
Tables \ref{tab:AvgFFees_AP_COMP} -\ref{tab:MAXFFees_AP_COMP} give the results on the larger instances from the AP data set. We report runtimes, in seconds, ($t(s)$), number of explored nodes (\textit{Nodes}), and percent optimality gaps at termination ($\%$GAP). The values displayed in each row are averages over the 10 instances with the corresponding number of nodes and carriers.
        \begin{table}[H]
            \caption{Comparison of   \SDOD$_F$ and $r$\SDOD$_F$ on AP[100-200] with $|K| \in \{3,\,4,\,5\}$ with average fixed outsourcing fees. \label{tab:AvgFFees_AP_COMP}}
            \resizebox{\textwidth}{!}{\begin{tabular}{lrrrrrr|rrrrrr}
                                      & \multicolumn{6}{c|}{(EP)}                                                                                                                                            & \multicolumn{6}{c}{(IP)}                                                                                                                                            \\
                                      & \multicolumn{3}{c}{$SDOD_F$}                                                     & \multicolumn{3}{c|}{$rSDOD_F$}                                                   & \multicolumn{3}{c}{$SDOD_F$}                                                     & \multicolumn{3}{c}{$rSDOD_F$}                                                   \\
        \multicolumn{1}{c}{$|V|.|K|.|H|$} & \multicolumn{1}{c}{Prep.} & \multicolumn{1}{c}{Load.} & \multicolumn{1}{c}{Run.} & \multicolumn{1}{c}{Prep.} & \multicolumn{1}{c}{Load.} & \multicolumn{1}{c|}{Run.} & \multicolumn{1}{c}{Prep.} & \multicolumn{1}{c}{Load.} & \multicolumn{1}{c}{Run.} & \multicolumn{1}{c}{Prep.} & \multicolumn{1}{c}{Load.} & \multicolumn{1}{c}{Run.} \\ \hline \hline
        100.3.5                         & 29.39                     & 77.00                     & 6.26                     & 29.45                     & 71.50                     & 152.49                    & 13.74                     & 2.55                      & 0.15                     & 69.58                     & 2.54                      & 0.36                     \\
        100.4.5                         & 34.05                     & 88.86                     & 9.50                     & 34.24                     & 81.84                     & $-\quad$                        & 33.98                     & 4.13                      & 0.34                     & 129.75                    & 4.13                      & 0.97                     \\
        100.5.5                         & 38.67                     & 100.63                    & 12.99                    & 38.58                     & 90.64                     & $-\quad$                        & 33.40                     & 6.06                      & 0.44                     & 188.67                    & 6.01                      & 1.61                     \\
        120.3.6                         & 43.03                     & 148.23                    & 13.90                    & 42.85                     & 137.75                    & $-\quad$                        & 14.50                     & 3.78                      & 0.22                     & 42.23                     & 3.69                      & 0.54                     \\
        120.4.6                         & 49.58                     & 169.75                    & 19.43                    & 49.59                     & 156.01                    & $-\quad$                        & 31.82                     & 6.12                      & 0.50                     & 78.32                     & 6.19                      & 1.30                     \\
        120.5.6                         & 58.10                     & 189.60                    & 29.92                    & 58.46                     & 171.87                    & $-\quad$                        & 36.73                     & 9.04                      & 0.84                     & 114.00                    & 9.10                      & 2.16                     \\
        140.3.7                         & 59.61                     & 256.28                    & 26.90                    & 59.84                     & 240.77                    & $-\quad$                        & 23.30                     & 5.30                      & 0.33                     & 77.24                     & 5.27                      & 0.96                     \\
        140.4.7                         & 70.94                     & 290.38                    & 35.97                    & 71.40                     & 267.80                    & $-\quad$                        & 51.80                     & 8.47                      & 0.77                     & 145.72                    & 8.60                      & 2.41                     \\
        140.5.7                         & 84.70                     & 320.05                    & 57.66                    & 85.01                     & 294.41                    & $-\quad$                        & 57.15                     & 12.60                     & 1.11                     & 213.67                    & 12.60                     & 4.56                     \\
        160.3.8                         & 82.59                     & 423.94                    & 45.66                    & 83.16                     & 397.52                    & $-\quad$                        & 33.67                     & 7.13                      & 0.43                     & 132.37                    & 7.12                      & 1.41                     \\
        160.4.8                         & 102.69                    & 476.20                    & 64.66                    & 102.69                    & 436.70                    & $-\quad$                        & 78.96                     & 11.45                     & 1.04                     & 251.59                    & 11.31                     & 3.55                     \\
        160.5.8                         & 122.39                    & 529.62                    & 110.33                   & 122.71                    & 478.08                    & $-\quad$                        & 82.01                     & 17.03                     & 1.37                     & 360.29                    & 16.99                     & 7.50                     \\
        180.3.9                         & 147.48                    & 660.11                    & 129.65                   & 146.56                    & 625.00                    & $-\quad$                        & 46.55                     & 8.99                      & 0.54                     & 211.10                    & 9.18                      & 2.14                     \\
        180.4.9                         & 164.66                    & 738.98                    & 177.13                   & 165.47                    & 673.10                    & $-\quad$                        & 112.77                    & 14.43                     & 1.35                     & 392.49                    & 14.40                     & 5.85                     \\
        180.5.9                         & 201.36                    & 807.43                    & 210.13                   & 199.61                    & 731.78                    & $-\quad$                        & 116.41                    & 21.54                     & 1.99                     & 572.64                    & 21.18                     & 10.21                    \\
        200.3.10                         & \multicolumn{6}{c|}{\multirow{3}{*}{Out of memory.}}                                                                                             & 60.04                     & 10.86                     & 0.57                     & 308.86                    & 11.22                     & 2.83                     \\
        200.4.10                         & \multicolumn{6}{c|}{}                                                                                                                                                & 153.03                    & 17.54                     & 1.54                     & 587.91                    & 17.54                     & 8.31                     \\
        200.5.10                         & \multicolumn{6}{c|}{}                                                                                                                                                & 152.83                    & 26.13                     & 2.28                     & 872.72                    & 25.93                     & 16.04                    \\ \hline
        \end{tabular}}
        \end{table}

        \begin{table}[H]
            \caption{Comparison of \SDOD$_F$ and $r$\SDOD$_F$ on AP[100-200] with $|K| \in \{3,\,4,\,5\}$ with maximum fixed outsourcing fees. \label{tab:MAXFFees_AP_COMP}}
            \resizebox{\textwidth}{!}{\begin{tabular}{lrrrrrr|rrrrrr}
         & \multicolumn{3}{c}{(\emph{EP}$_F$)} & \multicolumn{3}{c|}{($r$\emph{EP}$_F$)} & \multicolumn{3}{c}{(\emph{IP}$_F$)} & \multicolumn{3}{c}{($r$\emph{IP}$_F$)} \\
        \multicolumn{1}{c}{$|V|.|K|.|H|$} & \multicolumn{1}{c}{Prep.} & \multicolumn{1}{c}{Load.} & \multicolumn{1}{c}{Run.} & \multicolumn{1}{c}{Prep.} & \multicolumn{1}{c}{Load.} & \multicolumn{1}{c|}{Run.} & \multicolumn{1}{c}{Prep.} & \multicolumn{1}{c}{Load.} & \multicolumn{1}{c}{Run.} & \multicolumn{1}{c}{Prep.} & \multicolumn{1}{c}{Load.} & \multicolumn{1}{c}{Run.} \\ \hline \hline
        100.3.5 & 29.25 & 77.31 & 7.99 & 29.24 & 70.86 & 1.64 & 20.57 & 2.54 & 0.14 & 141.58 & 2.53 & 0.31 \\
        100.4.5 & 34.27 & 88.74 & 13.58 & 34.49 & 81.89 & 2.05 & 33.93 & 4.09 & 0.29 & 250.63 & 4.12 & 0.74 \\
        100.5.5 & 38.71 & 100.72 & 17.97 & 38.70 & 90.81 & 2.48 & 50.28 & 5.95 & 0.42 & 389.05 & 5.99 & 1.20 \\
        120.3.6 & 42.43 & 148.11 & 15.97 & 42.37 & 138.44 & 3.28 & 19.71 & 3.83 & 0.21 & 85.44 & 3.73 & 0.48 \\
        120.4.6 & 49.40 & 169.03 & 28.33 & 49.47 & 156.22 & 4.12 & 32.67 & 6.12 & 0.45 & 150.76 & 6.17 & 1.01 \\
        120.5.6 & 58.36 & 190.46 & 41.35 & 58.71 & 173.06 & 5.36 & 48.46 & 9.12 & 0.77 & 235.08 & 9.18 & 1.52 \\
        140.3.7 & 58.80 & 256.07 & 29.99 & 59.24 & 240.48 & 5.82 & 31.45 & 5.24 & 0.38 & 159.87 & 5.22 & 0.81 \\
        140.4.7 & 71.43 & 290.61 & 50.15 & 71.68 & 267.02 & 7.36 & 52.20 & 8.41 & 0.67 & 282.21 & 8.50 & 1.76 \\
        140.5.7 & 84.24 & 320.12 & 78.05 & 84.97 & 294.77 & 9.06 & 78.38 & 12.72 & 0.96 & 436.44 & 12.54 & 2.70 \\
        160.3.8 & 83.82 & 419.23 & 50.73 & 84.17 & 398.48 & 9.77 & 47.53 & 6.95 & 0.42 & 274.61 & 6.99 & 1.08 \\
        160.4.8 & 103.99 & 469.25 & 87.62 & 104.99 & 439.15 & 12.36 & 79.03 & 11.25 & 0.89 & 480.59 & 11.13 & 2.62 \\
        160.5.8 & 122.05 & 519.08 & 167.23 & 122.25 & 480.06 & 15.60 & 116.45 & 16.79 & 1.26 & 743.08 & 16.79 & 3.64 \\
        180.3.9 & 120.13 & 644.85 & 99.53 & 120.00 & 615.96 & 54.89 & 67.52 & 8.86 & 0.61 & 436.18 & 9.06 & 1.79 \\
        180.4.9 & 149.42 & 719.53 & 155.15 & 149.47 & 679.01 & 83.81 & 112.79 & 14.31 & 1.22 & 773.14 & 14.32 & 3.41 \\
        180.5.9 & 174.64 & 792.17 & 225.16 & 175.66 & 740.28 & 74.82 & 167.91 & 21.34 & 1.79 & 1197.13 & 21.25 & 5.24 \\
        200.3.10 & \multicolumn{6}{c|}{\multirow{3}{*}{Out of memory}} & 90.78 & 10.83 & 0.61 & 641.81 & 10.81 & 2.26 \\
        200.4.10 & \multicolumn{6}{l|}{} & 153.22 & 17.61 & 1.45 & 1157.00 & 17.75 & 4.56 \\
        200.5.10 & \multicolumn{6}{l|}{} & 229.71 & 26.31 & 1.97 & 1820.99 & 25.91 & 6.81 \\ \hline
        \end{tabular}}
        \end{table}
\end{section}

\begin{section}{Percentage of served commodities and net profit comparison}\label{app:service-rate}
 Table \ref{tab:FFees_CAB_COMCOMP} displays the average percentage of served commodities in optimal solutions to (\emph{IP}), (\emph{IP}$_F$), and ($r$\emph{IP}$_F$). The row corresponding to each number of nodes ($n$)  shows average results over the three instances with $n$ nodes and varying number of carriers $|K|\in\{3,\,4,\,5\}$.      
        \begin{table}[H]
            \caption{Service level comparison among \SDOD\,, \SDOD$_F$, and $r$\SDOD$_F$ on CAB[30-100].\label{tab:FFees_CAB_COMCOMP}}
            \resizebox{\textwidth}{!}{\begin{tabular}{lrrrrr}
                \multicolumn{1}{c}{} & \multicolumn{4}{c}{} & \multicolumn{1}{l}{} \\  
                \multicolumn{1}{c}{} & \multicolumn{2}{c|}{Average} & \multicolumn{2}{c|}{Maximum} & \multicolumn{1}{c}{Optimal} \\
                \multicolumn{1}{c}{Instance} & \multicolumn{1}{c}{\SDOD$_F$ ($\%$)} & \multicolumn{1}{c|}{$r$\SDOD$_F$ ($\%$)} & \multicolumn{1}{c}{\SDOD$_F$ ($\%$)} & \multicolumn{1}{c|}{$r$\SDOD$_F$ ($\%$)} & \multicolumn{1}{c}{\SDOD\,($\%$)} \\ \hline \hline
                30 & 13.86 & \multicolumn{1}{r|}{41.38} & 11.60 & \multicolumn{1}{r|}{42.96} & 44.00 \\
                40 & 12.88 & \multicolumn{1}{r|}{43.37} & 10.15 & \multicolumn{1}{r|}{44.98} & 46.37 \\
                50 & 11.00 & \multicolumn{1}{r|}{42.84} & 8.19 & \multicolumn{1}{r|}{44.65} & 46.91 \\
                60 & 9.18 & \multicolumn{1}{r|}{39.98} & 6.31 & \multicolumn{1}{r|}{41.97} & 44.06 \\
                70 & 8.21 & \multicolumn{1}{r|}{37.22} & 5.63 & \multicolumn{1}{r|}{39.98} & 41.97 \\
                80 & 8.63 & \multicolumn{1}{r|}{37.08} & 5.97 & \multicolumn{1}{r|}{39.76} & 41.64 \\
                90 & 7.19 & \multicolumn{1}{r|}{31.39} & 5.05 & \multicolumn{1}{r|}{33.77} & 35.45 \\
                100 & 6.81 & \multicolumn{1}{r|}{32.59} & 4.92 & \multicolumn{1}{r|}{35.50} & 37.19 \\ \hline
                Average& 9.72 &	38.23	& 7.23 &	40.45 &	42.20
            \end{tabular}}
        \end{table}

Figure \ref{fig:CABProfitAvg} compares the net profit associated with 
        for optimal solutions to \SDOD, \SDOD$_F$, and $r$\SDOD$_F$.
        
      \begin{figure}[H]
            \centering
            \subfigure[Maximum fixed outsourcing fees]{
                \begin{tikzpicture}
                    \begin{axis}
                        [
                        xlabel={Nodes}, 
                        ylabel={Total Profit}, 
                        axis lines=left, 
                        xmin = 29, xmax = 105, ymin = 0,
                        legend style={at={(0.5,-0.2)}, anchor=north,legend columns=3},
                        xtick= data,
                        legend to name = {mylegend2}
                        ]
                        \addplot[blue, mark=square, solid] table[col sep=comma, x=Nodes, y=Mean MAX Relaxed]{01_Profit.csv};
                        
                        \addplot[red, mark=square, solid] table[col sep=comma, x=Nodes, y=Mean MAX BilvlOpt]{01_Profit.csv};
                        
                        \addplot[green, mark=square, solid] table[col sep=comma, x=Nodes, y=Mean IP]{01_Profit.csv};

                        \legend{$r$\SDOD, \SDOD$_F$, \SDOD}
                        
                    \end{axis}
                \end{tikzpicture}
            }\quad
            \subfigure[Average fixed outsourcing fees]{
                \begin{tikzpicture}
                    \begin{axis}
                        [
                        xlabel={Nodes}, 
                        ylabel={}, 
                        axis lines=left, 
                        xmin = 29,
                        xmax = 105,
                        ymin = 0,
                        xtick= data
                        ]
                        \addplot[blue, mark=square, solid] table[col sep=comma, x=Nodes, y=Mean AVG Relaxed]{01_Profit.csv};
                        
                        \addplot[red, mark=square, solid] table[col sep=comma, x=Nodes, y=Mean AVG Bilvl Opt]{01_Profit.csv};
                        
                        \addplot[green, mark=square, solid] table[col sep=comma, x=Nodes, y=Mean IP]{01_Profit.csv};

                    \end{axis}
                \end{tikzpicture}
                }
            \caption{Average total profit comparison for CAB dataset}
            \label{fig:CABProfitAvg}
        \end{figure}
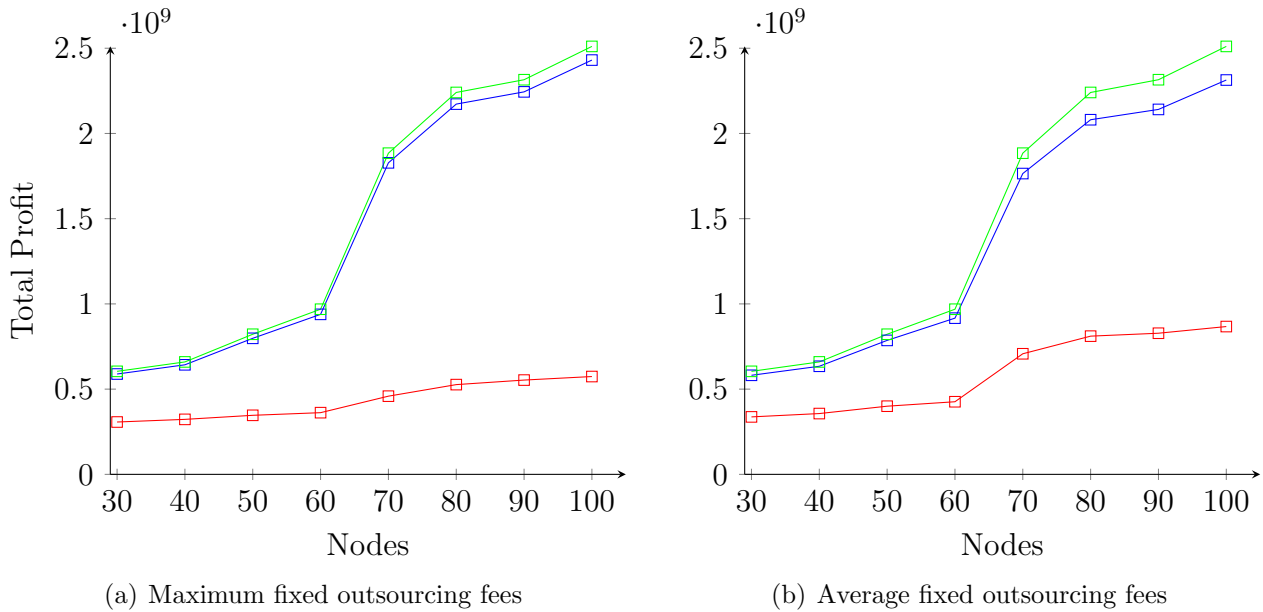
\end{section}   
\end{document}